# Dynamic systems with quantum behaviour

A.P. Alexandrov

*Research Institute*
*Nikolayev Astronomical Observatory, Ukraine*

Abstract

It is argued that the world is a dissipative dynamic system, a phase flow of which is formed by conformally-symplectic mapping. The key assumption is that the concept of energy in microcosm makes sense only for the steady motions corresponding to quantum eigenstates. The constant, which determines the exponential phase volume contraction, is supposed to be a new universal constant, in addition to the speed of light and Planck constant. It is shown that statistical treatment of quantum objects as the ensembles concentrated on smooth connected attractors provides a simple explanation of stochastic behaviour of these objects as well as leads to a natural interpretation of the wave function, stationary Schrödinger equation, and scattering matrix.

To validate the general hypotheses stated in the work, some physical models are presented. In particular, the models support the view that the inertial motion and quantum properties are basically determined by the vacuum as a dynamic subsystem. The matter-vacuum interaction is described formally by means of multivalued Hamilton function. A drawback of all these models is the non-locality of equations of motion stemming from the non-locality of variation of multivalued Hamilton function. This drawback can be overcome in the field theory using multivalued functionals with a local first variation opened by S. P. Novikov. Hence, a prospect of construction of the unified field theory in Einstein spirit, without a procedure of quantization is opened.

*E-mail:* a-lex@mksat.net

# Contents

# Introduction.

This work suggests a conjecture about possible application of the *conformally-dissipative systems* (further – CD-systems, or simply CD) in theoretical physics.

This term is used by the author for denotation of dynamic systems the phase flow of which conformally contracts a symplectic form[1]. In some respects such systems looks more fundamental, than Hamiltonian ones. Indeed, these systems give the simplest *dissipative* generalization of Hamiltonian systems, are conceptually *simpler*, and at the same time admit a formal Hamiltonian limit. Conversely, it is impossible to obtain the CD from the Hamiltonian system as some kind of limit, obviously. In general case the CD-systems are *irreversible* and admit non-trivial attractors, in particular limit cycles.

The assumption is that conception of a CD-system, in its field variant, can become the basic *dynamic principle of future fundamental physical theory*. Here is implied that dynamics of the system will be described by the *classic* fields[2], but the conception of quantized field will connected with the next, statistical floor of theoretical building.

To form this point of view *the Newton's dynamic postulate* (NDP) served as a starting point. The exact formulation of this postulate will be set below. In addition to common considerations the dynamic models demonstrating the possibility of quantum behavior on the dynamic basis confirm this postulate. Among them are the model of quantum harmonic oscillator, model of relativistic massless particle, model of spin, and model of fermion type.

Obviously, using of the dissipative dynamic systems at the fundamental level supposes *changes in the physical foundations* of existent theory. Among of these necessary changes I will note the following.

Firstly, the restrictions on the conceptions *of energy* and Hamiltonian are necessary. Here the *classic* (not quantum) conceptions are meant. Secondly, non-existence of the «relativistic CD-systems» results in approximate nature of the known relativistic dependences. Thirdly, new frameworks inevitably lead to the conception of *vacuum* as a dynamic subsystem which provides quantization of the material fields and, presumably, generalizes gravitation.

---

[1] Such systems could be called as conformally-symplectic ones, but it seems more important to emphasise their dissipative character.

[2] Here a classic field is meant in wide sense: a) configurations of the field are locally set by the finite number of *numerical* functions on a finite-dimensional basic manifold; b) these functions admit *local* variations.

It should be emphasized once more, that NDP and all the mathematics of the CD-systems supposes the *dynamic* nature of theory and, in particular, the possibility in principle of *unstatistical*, «deterministic» description in a quantum domain. It is possible to say in this sense, that here offer a variant of mathematics for a «deterministic quantum theory»[3].

Another physical hypotheses being automatically introduced by the formalism of CD-systems, is the assumption that there exist *a new fundamental constant* $(\kappa)$ in nature. This constant has a frequency dimensionality and is substantially different (to a less side) from Planckian frequency $\nu_{Pl} \sim 10^{43} \sec^{-1}$. Further we will call it *the dissipative constant*. Together with Planck constant, this constant is responsible for formation of the quantum eigenstates. Thus, the supposed theory must include three fundamental dimensional constants: $c, \hbar, \kappa$.

Judging from the presented models, a value $\kappa$ must be small in comparison with typical atomic frequencies, but large enough on macroscopic measures, for example, belonging to an interval of $10^{12} - 10^3 \sec^{-1}$. Such a value is necessary for the accordance with quantum theory. The point is that the quantum spectrum of harmonic oscillator is reproduced only at $\kappa \to 0$, with deviation of the second order $O(\hbar\kappa^2/\omega_0)$, where $\omega_0$ is the oscillator frequency. On the other hand, a relaxation time to the quantum eigenstates must not be too large, that bounds $\kappa$ from below.

Of course, it would be desirable immediately to guess the phenomena in which so unusual scale of the dissipative constant can exhibits directly. In this respect the next observations can be presented. Firstly, new fundamental mass $m = \hbar\kappa c^{-2}$ can have an order of now unknown a mass of neutrino. Secondly, it is noticed (see, for example, [23]) that to *the density of space vacuum* measured recently

$$\rho_V = (4 \pm 0.3) \cdot 10^{-30} g \cdot sm^{-3} = \hbar c^{-1} \lambda_V^{-4},$$

the length of $\lambda_V \approx 1 mm$ corresponds. This length $\lambda_V$ is close to the average wave-length of *cosmic microwave background*, and the corresponding value $\kappa = c\lambda_V^{-1} \approx 3 \cdot 10^{11} \sec^{-1}$ seems acceptable.

In spite of their dissipativity, the CD-systems completely keep within the standard frameworks of Lagrange formalism. If the unified field CD-theory is possible, then, as will be clear further, the proper Lagrangian must look like

$$L' = e^{\kappa t}\left(L(\mathbf{q},\dot{\mathbf{q}}) - \hbar \cdot \frac{d}{dt}\Phi(\mathbf{q})\right),$$

---

[3] As an example of the mathematics, pursuing a similar purpose, but based on other initial ideas, point out to known works of G. `t Hooft (see, for example, [15]).

where $t$ is *the absolute time*, $\mathbf{q}$ – the fields, in the conditional notation, $L(\mathbf{q},\dot{\mathbf{q}})$ – some analogue of a classic (relativistic?) Lagrangian, $\Phi(\mathbf{q})$ – the *multivalued* dimensionless functional, having multiple $2\pi$ periods[4]. A multiplier $\hbar$ specifies that exactly this functional is responsible for the quantum behaviour of the system. It is natural to accept, that a total derivative $d\Phi(\mathbf{q})/dt$ will have the form of space integral of some local density that will lead to equations in partials. General construction of such multivalued functionals is given in [8].

To avoid misunderstanding it should be noted that Hamiltonian description in variables $(\mathbf{p}' = \partial L'/\partial \dot{\mathbf{q}},\ \mathbf{q})$ corresponding to Lagrangian $L'$ is not compatible with attractors and that is why it has no physical sense. In this case the true «Hamiltonian variables» are $(\mathbf{p} = \partial L/\partial \dot{\mathbf{q}},\ \mathbf{q})$, but in these variables system is not Hamiltonian, but a conformally-dissipative one.

The structure of the work is as follows. There are three chapters and two appendixes at the end of the main text.

Chapter 1 is devoted to general questions. It begins with criticism of the known objections against the dynamic ideology at fundamental level, then makes an attempt to comprehend the «dynamic essence» of the existent theory, and then proposes the hypotheses and heuristic considerations about of a future theory, which lead, in particular, to the conception of a CD-system.

Chapter 2 is the preliminary mathematical considerations of the CD-systems from the point of view indicated higher, and also the search of the connection with quantum formalism[5]. The proofs of the majority of mathematical statements (further called *theorems*, for the sake of uniformity) are omitted. Basically they have computational character and it is easy to prove them.

In Chapter 3 the primary heuristic model of harmonic oscillator is generalized to the model of coherent states (the CS-model). Further the different special cases of it are examined that make possible to judge about the principle possibilities of the CD-systems for description of real physics.

The physically most interesting local field models are not considered in the given work. At first it will be enough to demonstrate such systems *do exist*, at least, in formal sense[6]. I hope that results of this work will find continuation also in the field models of CD-type.

---

[4] It is possible to suppose existence other periods. In particular, such periods can determine the *type of statistics* of particles. It should be noted also, that relations of different periods are *mathematical* constants.

[5] It seems this class of dynamic systems does not especially be analyzed now.

[6] A problem is in checking a regularity of a system.

# Chapter I. General reasonings and hypotheses.

## §1. Newton's dynamic postulate (NDP).

I suggest to call so the basic Newton idea that physical world is (in modern terminology) a continuous *dynamic system*, with the equations of motion

(0.1) $$dx/dt = V(x).$$

It is supposed, that these equations define a semi-flow $g_t$ (continuous one-parametrical semi-group of transformations) on some manifold $M$ which points represent the *real states* of system being understood regardless to any «observations» or «measurements». Newton, however, *identified* a dynamic parameter $t$ *(absolute, truly mathematical time)* with empirical time in any inertial system of reference though recognized the existence of «relative time». Taking account of the theory of relativity it is necessary, of course, *to refuse* from that identification.

Let's emphasize, that *contrary to* the modern quantum-relativistic views, NDP assumes the following:

− the world is objectively real, irrespectively of *measurements, observers, consciousness*;
− time course *is not an illusion* (and of vague origin) in consciousness of observers;
− the probability should not be as a fundamental principle but should *be explained*.

These attractive (however, probably, not for all) consequences force to think once again on, whether there was a refusal from NDP in the last century only a provisional measure, and whether it is necessary nevertheless to accept it, as *the basic* physical hypothesis.

### 1.1. Quantum mechanics versus NDP.

The main objections against NDP are put forward by the quantum theory.

The first objection consists that the uncertainty principle is incompatible with the usual dynamic description the states of objects. Among the quantities describing the states *of quantum objects*, cannot be the canonically conjugated pairs of dynamic variables.

The attempts to declare this principle erroneous, to prove that it is only the property reflecting incompleteness of quantum formalism has led to the opposite result, as is known. Later on all the doubts that the uncertainty principle concerns an essence of quantum objects have disappeared.

However it would be an error to consider NDP as incompatible with this negative principle. Actually it is quite possible *to agree* with it. Eventually, the uncertainty principle asserts, that the quantum objects are not the objects in usual sense, and that the states of quantum objects is not the same as the states of a dynamic system. So, in such cautious formulation, the contradiction with NDP is absent.

The second objection against NDP is given by known theorems of impossibility of theories of hidden variables, beginning from the von Neumann's theorem, and finishing at Bell's inequalities and their generalizations.

The answer to this objection is close to the answer to the first objection and, in essence, is already known. For example, in article [13] written in connection with discussion concerning the interpretation of quantum mechanics, it is noted:

«The system for which the Bell's inequalities are valid, can and should be exhaustively presented as actual *set* of some objects-elements which are characterized *itself* by actually inherent properties».

In other words, similar theorems ignores the dynamic character of a quantum measurement procedure which are not reduced actually to ascertaining of some «preexisting» properties as it is meant, for example, in the classic mechanics. Properties of quantum objects grow out *of interaction* of the device with some deeper essence, and in this sense they are similar *to relations*, i.e. have *relational* character, as well as the author of article confirms further.

The third objection against NDP consists in obvious reference to probabilistic character of the quantum theory assuming experimentally confirmed «primary probability» at fundamental level.

In opposition to similar views, however, it is easy to show that *the visible probabilistic behaviour of quantum objects can be a simple consequence of their special dynamic treatment*. The important help is given in this respect by multi-world interpretation of quantum mechanics of Hugh Everett. Recall that Everett has identified the formal unitary Schrödinger dynamics with the real world dynamics, and then has rather logically come to conclusion about Universe branching, occurring at each quantum measurement with potentially non-unique outcome.

It is natural to think, what could it mean within the frameworks of usual dynamic approach. Here it is possible to begin with the remark, that the parallel Universes actually are thought up long before Everett, under the name *of phase space* of dynamic system. Mentally enough to actualize this space and then, at due change of terminology, it is possible to consider, that the phase points which are carried away by a phase flow, describe the «parallel existence of Universes».

In such an interpretation to the idea about branching of some Universe there corresponds the thought about branching of a corresponding integral curve in the expanded phase space (including time). But it is impossible in usual dynamic theory, and the idea is necessary for modifying: the bundle of trajectories can branch only. After this minor editing the Everett concept becomes quite dynamic, and admits the translation into usual language, without a mention of the Universes.

The basic heuristic statement obtained by this way, is that the unitary quantum evolution should be connected with *a bundle* (or with an ensemble) of usual *dynamic* trajectories. In particular, it means, that *the state* of a quantum system (or a quantum object) should be connected not with a phase point but with some *subset* in phase space. Obviously, the converse statement is also valid, i.e. from the last assumption follows also the possibility of branching. In initial physical terms it is equivalent *to an explanation* of «probabilistic» behaviour of the quantum objects, owing to their treatment as *the statistical ensembles*.

Note, that the given approach turns the correlation of results of remote measurements to obvious possibility, in spite of the visible accidental outcome of each measurement (that forms the EPR paradox essence). The correlation becomes clear in view of presence *of a dynamic basis* of all measurements. Also, automatically is solved «a choice problem» at quantum measurement, because ultimately we always deal with *single* phase trajectory.

### 1.2. Theory of relativity versus NDP.

Let's consider now the objections against NDP which can be put forward in connection with the theory of relativity. Obviously, the basic objections here cannot be, as, in effect, the relativity proclaims the existence *of the set of dynamic interpretations* of reality *equal in rights*, instead of the one fundamental interpretation. The basic objection, possibly, is that NDP separates again all that the relativity theory has already united: the space and time, energy and momentum, electric and magnetic fields, and so on. It looks as refusal from the symmetry, and from a beautiful geometrical picture (degeometrization). Besides, unlike the principle of relativity, NDP is represented purely speculative, in any way empirically not supported principle.

In reply here it is possible to note, that the beauty and unifying tendencies of the relativity have basically a group-theoretical origin. The source of unification is the Lorentz group. But, as the history of quarks shows, for example, the symmetry reasons and the group approach can be only *the first step* on a way *to deeper* dynamic picture, in which this symmetry are not present. Besides, the beauty of the theory depends on the mathematicians who developed it. Therefore it is possible to assume, that, in the long term, the dynamic theory will become not less graceful, than relativity.

Let me remind also, that cosmic microwave background sets a natural *absolute frame of reference*. It is possible to reject mentally, of course, a matter, and to declare, that it fills a relativistic-invariant vacuum, but at the unbiased approach the given fact is the direct experimental evidence in favour of NDP. The answer to a problem where and how to search for others, more essential experimental evidences in favour of NDP, only the future theory can give.

### 1.3. Philosophers versus NDP.

Perhaps, it is necessary to mention a *philosophical* objection against NDP, as the concept *of a determinism* excluding the «creative evolution».

Similar reasons, however, after opening of the possibility of dynamic chaos, and even in low-dimensional systems, have lost a sense. It is clear now, that the dynamic concept if not to absolutise it and to consider it «sensibly», *is not identical to the idea of determinism*. The impossibility of absolutely exact knowledge of initial data excludes the real possibility of exact predictions in the future. Besides, all models can be only approximate, and all real systems are open. Therefore, even the «deterministic» dynamic systems (including the classic mechanics) have physical sense actually only in *a statistical variant* that assumes, for example, replacement of point states by $\delta$-shaped phase distributions of probabilities.

### 1.4. The plan for the further.

So, nothing hinders us to accept NDP as a primary fundamental principle.

The general plan of the further researches (partially realized in the given work) consists in the beginning to analyze the dynamic content of the theory of relativity and the quantum theory, and to try to find their true dynamic interpretation. This analysis should help the search of bases a new, *adequate to NDP* mathematical apparatus. The author believes that last problem is solved by the concept of the CD-systems. Then, as the *deduction* of quantum and relativistic principles here is meant as approach to deeper dynamic laws, it is necessary to establish connection with the contemporary theory. As to connection with the quantum theory, at least, *at basic level*, the CD-systems *give* such a connection. It means, in particular, that CD-systems give one more, in addition to already existing, *interpretation of quantum mechanics*. As to con-nection with the relativity, probably, for solving of this problem we must consider only the *field* CD-systems.

## §2. Relativity and dynamics.

Let's consider the features of the dynamic approach to the theory of relativity.

### 2.1. The concept of inertial motion.

Obviously, the concept of «a frame of reference» has the exclusively empirical sense, and it does not be entered in abstract frameworks of NDP. Owing to this circumstance the principle of relativity does not admit *the equivalent* dynamic formulation. On the other hand, the relativity principle is based on the concept of inertial motion which in the essence already is dynamic. Therefore it is natural to search for dynamic definition not for a relativity, but for inertia. It also will be a dynamic substitute of a principle of relativity.

The only thing that is necessary for such definition is the concept of dynamic symmetry. The group $G$ of transformations of phase space is called a dynamic symmetry if the phase flow $g_t$ belongs to this group, and $g_t$ does not belong to the center of this group[7]. In the special theory of relativity (SR) the role of such group $G$ plays the Poincare group.

Suppose that some solution $x(t) = g_t x_0$ to the equations of motion describes the isolated material system which is in a «rest state». Then the solution, having the initial data which are obtained by action of group transformation on the initial data of a «rest» subsystem, *by definition* describes the *inertial «movement as a whole»* of the same subsystem[8]. In other words, the transition from the system in rest to the inertially moving one looks like

$$g_t x_0 \to g_t h x_0, \quad h \in G.$$

Such understanding of inertial motion is non-trivial for the case of the Poincare group. For example, from it Fitzgerald-Lorentz's reduction follows. It is important to note, that from the given definition it is possible to obtain the dynamic formulation *of all* relativistic effects[9]. In the standard interpretation these are considered to be purely «kinematic» or «geometrical». It allows to consider the concept of the dynamic symmetry as adequate expression of SR ideas in the dynamic frameworks.

It is necessary to note, that the given point of view does not use the concept of frame of reference (because it is supposed, that the states of a system are objective, «real» states) and, especially, does not use the space of events (i.e. space-time).

### 2.2. The vacuum concept.

However for the general theory of relativity (GR) the approach with dynamic symmetry any more does not work. Strictly speaking, the symmetry and «inertial motion» disappears. A metric tensor $g_{ik}$ becomes *a dynamic variable* owing to what the «free» Lagrangians of particles and fields turn, as a matter of fact, into *interaction Lagrangians* with $g_{ik}$ as with a usual physical field. The inertial dynamics is replaced with some interaction of universal character. The concept of inertial movement remains only in approach when $g_{ik}$ is considered as an external field and when the spin of particles is ignored.

---

[7] Otherwise, the symmetry can be called as usual, or *geometrical*. Dynamic symmetry also can be subdivided into types. For example, the time course in moving system depends on whether $g_t$ belongs to commutant $G'$, or not.

[8] Change of the initial moment leads to the equivalent definition.

[9] See Appendix 1.

It is natural to accept, that in the future dynamic theory we will have some fundamental field of deeper level instead of a metric field. It means that metric tensor $g_{ik}$ is necessary to understand as *the incomplete* empirical description of this more fundamental field in terms of characteristics of the approximate «inertial motion» of test particles, caused by interaction with the given field. Let's call this field *a vacuum*. The metrics can be only one of vacuum displays. Presumably, another display of the vacuum properties is the quantum effects[10], for example, the wave properties of particles which are present *also in flat space*.

*The universality* of interaction with vacuum demands some explanation.

**Hypothesis 1**. *Interaction with a vacuum concerns the nature of particles (e.g. is the cause of their existence), and causes the property of inertia of particles*.

The fact, that vacuum serves as «reference point» for *forces*, *energy*, and others «material characteristics» can be considered as consequence of this hypothesis.

## 2.3. Reconstruction of effective geometrical structure of space-time.

Generally speaking, NDP does not exclude the use of traditional geometrical description in the theory of relativity. At the empirical level this description has arisen as consequence of experimenting in various inertial frames of reference. But if we already have the equations of the dynamic theory, the laws of motion *of test particles* should become a sufficient *theoretical* basis for reconstruction of four-dimensional geometry at macroscopic level.

It is logically admissible, that reconstruction of the metrics or other structure will be realizable not always because the elementary or other test particles, under some physical conditions can simply not exist. The occurrence of geometrical singularities that will specify in necessity of search of a new suitable interpretation is in that case possible.

Thus, the future theory should not only give a dynamic explanation to the relativistic dependences, including the relativity principle, but also should establish the necessary *restrictions* for the geometrical approach.

## 2.4. Remarks on a four-dimensional formalism.

Probably, here and there we *already have been faced* with similar borders. Really, if GR will result from the approximate *geometrical interpretation* of a deeper dynamic theory, then the known singularities GR [5] can mean the *borders* of such interpretation. Possibly, the geometrical frameworks appear too narrow for a dynamics. That these general assumptions have not seemed too dubious, the toy model of a similar situation is presented in Appendix 2.

It is known also, that not in each GR world there are the global hypersurfaces of Cauchy data uniquely defining the future [5]. In such situations, on the contrary, the dynamic frameworks set by NDP are too narrow for geometry. So, because NDP is accepted, there is a doubt in a physical sense of such «worlds». Also, in a sense of any others geometrical worlds containing the loops in time, additional time-like dimensions, and similar *dynamic absurdities*.

---

[10] For this reason the metric tensor $g_{ik}$ is characterized here as the incomplete description of a vacuum.

Therefore geometrical and dynamic approaches lead, finally, to the *nonequivalent* theories. Thus, there is a potential danger of the absolutization of the existing geometrical formalism.

## §3. Quanta and dynamics.
### 3.1. The basic hypotheses.

The dynamic approach to the quantum theory, as well as any its other interpretation, should give some explanation to the basic non-classic features of this theory. In short they can be characterized as *discreteness*, *probability*, *coherence*.

Possible dynamic interpretation of likelihood behaviour of micro objects already was discussed above. Most likely, the states of these objects are necessary associate not with the points, but with the statistical ensembles in a phase space.

As to the nature of quantum coherence, this feature seems as a hint about the dynamic character of quantum background. Eventually, the phase flow can be characterized by some «phase». The existence of such «phase» is not something unusual to the dynamic systems. For example, in Hamilton mechanics to invariant Lagrange manifolds corresponds the action function $S$ satisfying the Hamilton-Jacobi equation, and the phase role play the function $\varphi = S/\hbar$.

Quantum coherence and, in particular, quantum correlations seem mysterious only within the frameworks of quantum mechanics when the «natural» postulate on accident of an outcome of each measurement is accepted. At the same time experimentally this postulate *is refuted* by correlation of outcomes of the simultaneous individual distant measurements. Actually experiment puts us before alternative: or quantum accident is non-local, or we deal with usual dynamic or «seeming» accident.

The most essential feature of the quantum world, certainly, property *of discreteness* is. At first sight, the problem consists in explaining from dynamic positions the discreteness of spectra of quantum observables, basic of which the energy is, first of all. However such formulation of a problem is not too fruitful because it is too adhered to an existing formalism. The matter is that this formulation implicitly assumes the «quantum jumps», which are incompatible with continuous dynamics. These jumps was entered already in the first Bohr model of hydrogen atom. We can say that Bohr had sacrificed a dynamic postulate in this model, but had preserved, however, *as a more fundamental*, the concept of the conservative energy. His relations $E_m - E_n = \hbar \omega_{mn}$, obviously, express the energy conservation law.

On the other hand, Bohr's idea about the *quantum steady-states* is quite dynamical one by the nature because Bohr described these states within the frameworks of classic mechanics. Therefore, if we wish to keep NDP, we should construct *a dynamic model* of these states, refusing at the same time from the quantum jumps.

The fundamental nature of the steady-states concept for the quantum physics follows from experiments on *measurement* of values of physical observables, because these experiments are reduced, finally, to detection of such states. It is clear from the theory. So, for example, the possible results of measurement of electron *spin projection* are deduced from consideration of the eigenstates of Pauli's Hamiltonian for electron in a magnetic field. If we could explain the steady-states, then the observable «discreteness» of energy and of all other observables thereby will be explained.

The core what is required from quantum steady-states in respect of dynamics is their *stability and clear splitting into the certain kinds*. It is known, that V. Weisskopf considered these properties as the cores in quantum physics, contrary to uncertainty and probability usually focusing attention. They provide the identity of atoms, molecules, elementary particles, and lead to definiteness of qualities of substance in the nature.

One more distinctive property of the quantum steady-states is their *flexibility*: the atomic levels of energy are capable to «adiabatic deformations» under influence of slowly changing external conditions. Contrary to it, for example, the topological charges in the field theory preserve the values.

If to look for something similar in mathematics the answer will be obvious. All specified properties evidently are inherent also for the limit cycles describing the auto-oscillations and for other attractors formed in the dissipative dynamic systems. And, what is especially important, any reasonable alternative to such connection of the continuous dynamics with discreteness, apparently, simply does not exist. It should be noted also that the smooth attractors are the submanifolds in a phase space, and are described *by the additional equations* which quite can play a role *of quantization conditions*.

On the other hand, it is necessary to recognize, that the idea of quantum auto-oscillations is absolutely not original and is rejected by the obvious reasons. Really, any auto-oscillatory process demands the *constant inflow of energy* from the outside. However, for example, for the isolated atom the similar process seems absolutely improbable.

The elementary way of overcoming this «no go» is to sacrifice, contrary to the Bohr approach, the fundamental nature of the concept of energy in favour of NDP, and to accept the following general point of view.

**Hypothesis 2.** *The dynamic law in quantum domain, unlike of the Hamilton dynamics, should admit the contraction of a phase volume[11] and the auto-oscillations. The concept of energy as secondary, arising on the basis of the dynamic law, should be modified.*

According to such approach, the «inflow and dissipation of energy» which is necessary for auto-oscillations, are only the admissible *formal* concepts which, however, should not be connected with usual physical understanding of energy. The most important in the accepted hypothesis is that it liquidates a taboo regarding use of the dissipative systems, and we immediately will take advantage of it.

**Hypothesis 3.** *To the quantum steady-states in the dynamic theory there correspond the limit cycles or others attractors*.

### 3.2. The qualitative consequences.

Now we should understand more in detail the physical sense of a last hypothesis. It supposes that we should operate with the limit cycles in the same manner as with the eigenstates in quantum mechanics (QM). In particular, having a limit cycle, it is necessary to be able to calculate a corresponding average $\overline{f}$ for any observable $f$, i.e. for a function. Obviously, such an average should be calculated as a usual time average. On the other hand, calculation of time averages leads to probability distributions. For example, if the point $x_0$ belongs to the basin of attraction of a limit cycle then the time average of $f(g_t x_0)$ looks as an average with respect to *the invariant distribution* $\rho = g_t \rho$ localized on a cycle

(1.1) $$\overline{f}(x_0) = \lim_{T \to \infty} \tfrac{1}{T} \int_0^T f(g_t x_0) dt = \langle \rho; f \rangle .$$

---

[11] In other words, the system should be *dissipative*.

Here $g_t$ denotes a phase flow, and flow action on distributions is defined by the formula

(1.2) $$\langle g_t \rho; f \rangle = \langle \rho; f \circ g_t \rangle .$$

Note, that a mix of distributions

(1.3) $$\rho = \alpha \rho_1 + \beta \rho_2, \quad \alpha + \beta = 1, \quad \alpha, \beta > 0,$$

corresponding to different limit cycles, cannot correspond to a QM eigenstate. It allows to improve the formulation of the basic hypothesis.

**Hypothesis 3'.** *The eigenstates of a quantum system will be described in the future dynamic theory as invariant distributions, concentrating on attractors, and indecomposable into a mix of the invariant distributions.*

Let's call such distributions *q-states*. Thus, there are states of two kinds – the *usual states, and q-states*. The usual state is the distribution localized at a point, but the support of the q-state is the attractor or its invariant part.

To all this we can attribute the following physical sense. Possibly, in quantum area we deal with the dynamic system the usual states of which, thanks to our ways of observation, and also thanks to the existence of «hidden» or fast dynamics, are empirically inaccessible. The «quantum measurements» are not instant and comprise some temporal averaging. It leads to that the usual states become «not physical» and the role of the true, «physical» states is transferred to the q-states. The theory which considers this circumstance, and operates exclusively with the q-states, should coincide with QM or to be close to it[12].

However for low frequencies when the averaging is not adequate, the description by means of q-states becomes empirically incomplete. In such cases the q-states exhibit themselves only as *statistical ensembles* that supposes a set of measurements. Unlike it, the high-frequency q-states are perceived in individual measurements.

Further, if aspire to keep the formal concept going from classic mechanics according to which behind a set of «states» there should be the certain *object* characterized by these states it is natural to find out, what the new objects (it is possible to call them *q-objects*) correspond to q-states. Obviously, owing to hypotheses accepted above, they need to be identified with the observable quantum objects. It is clear, that such approach provides formal correspondence between q-objects and usual objects. For this purpose it is enough to assume, that limit $\hbar \to 0$ leads to contraction of the extended q-states to usual point states.

---

[12] Actually the situation is not so simple. As it was noted above, the dynamic approach demands to include *a vacuum*. Therefore transition to that type of the description which is accepted in QM should be accompanied also by the inverse procedure *of an exclusion* of vacuum.

And still q-objects are rather specific. The q-object states (i.e. the q-states) can form a complicated set. It can be discrete, and also can be a smooth manifold. But the main difference is that q-states depend on the basic phase flow, contrary to the set of usual states, and, in particular, depend *on the external conditions* which define what and where are formed the limit cycles or others attractors.

Besides, q-states are concentrated on sets, instead of separate points, so the behaviour of q-objects is of generally *probabilistic* character. For a dynamics of q-object to have a direct physical sense, the corresponding q-states should be, certainly, the «high-frequency» ones, as it was already spoken.

Other obvious consequence of extent of q-states in phase space is *the uncertainty* of the observable values in these states: functions can accept the *various* values over attractors. Comparing this situation with that which is available in QM, we come to conclusion that the extent of q-states should be guaranteed by a nonzero value of Planck's constant. At the same time the possibility of exact measurement *of one* of two canonically conjugate variables in QM gives the hint, that q-states, probably, are placed on Lagrange submanifolds of symplectic phase space.

One more feature concerns the description of compound systems. If the subsystems interact, the attractor of compound system is not obliged to look like the direct product of attractors of subsystems. Therefore the q-state of compound system is not reduced to states of subsystems, and these subsystems are not in certain q-states. In this sense the compound q-object generally cannot be considered as «actually» consisting of q-objects which have originally formed it.

In general this picture is close to that what is available in QM. However, in the given interpretation of QM the Schrödinger equation can play the role of description of the steady motions only. But it is known, that this equation allows to predict the quantum transitions and to find their probabilities. The reasons of such discrepancy, obviously, that the Hypothesis 3 and corresponding dynamic interpretation, strictly speaking, concern only a stationary case.

### 3.3. The de Broglie's idea.

The general-dynamic picture cannot clear the questions connected with the evident, space-time interpretation of some paradoxical quantum phenomena. The phenomenon *of interference of single particles* concerns that, first of all.

For explanation of the given phenomenon it is necessary to involve *the additional* hypotheses including as well a postulate *of locality* of the basic equations. In this respect, the wave-pilot theory, or the de Broglie-Bohm mechanics seems to be a step in true direction.

In itself this theory, nevertheless, is unsatisfactory. First, in a case of many particles the given theory does not admit realization in three-dimensional space. Secondly, because the formal dynamics of Schrödinger equation is considered in this model as *a part of real dynamics*, the theory cannot describe a relaxation of arbitrary states to eigenstates. But the existence of similar property in the nature is supposed by the Hypothesis 3.

These drawbacks, nevertheless, are substantially connected with concrete *mathematical realization*. The *physical idea* that *interfere some field which is not transferring energy, and energy is transferred by the particles guiding by this field*, is remarkable simple and general, and it is necessary to take advantage of it. Obviously, the role of «a field which are not transferring energy» in our case can play a vacuum.

## §4. Heuristic reasonings about the mechanism of quantum auto-oscillations.

Relying on the above qualitative picture, let's try now to formulate more or less constructive requirements to the dynamic system. First of all, being guided by the Hypothesis 2, it is natural to reflect on search *of a new definition* of the energy more suitable to systems, admitting attractors. Really, if we admit this hypothesis and consider the energy to be a function in phase space, then we face with the following problems:

− energy cannot be the constant of motion defined *on all states* of system;

− not clearly also, why energy is constant of motion *over limit cycles* as it is dictated by conformity with QM.

It is interesting, that the idea about the necessity of special understanding of energy in quantum area was stated in 1926 by E. Schrödinger, in the letter to M. Planck [6]:

*Concept «energy» – it is something that we have introduced from macroscopic experience and only from it. I don't think that it is possible to transfer directly this concept to micromechanics and to speak about the energy of separate special oscillation. The energy property of separate special oscillation is frequency.*

Thus, Schrödinger suggested to identify energy with frequency on the basis of Planck's formula $E = h\nu$. Basically, such a solution *suits* us, because, at least, auto-oscillations can have a frequency. Eventually, it means, that energy has physical sense only for the considered above q-states. A lack is that such a solution in any way does not take account of the idea of auto-oscillations. Besides, not clearly, how (if not to use the quantum formalism) the property *of additivity* of energy can be provided.

Note, that nothing like the «own energy» exist for usual auto-oscillations. Instead, however, the *other* general characteristics are easily found out. Really, the physical auto-oscillatory process is characterized by such quantities:

− the frequency $\nu = 1/T$;

− the energy dissipation $\varepsilon$ during $T$, equal to its inflow to the system for the same time;

− the average energy dissipation rate $\varepsilon/T$, or $\varepsilon\nu$.

Considering these quantities, it is necessary to note the similarity of last expression with Planck's formula for energy:

$$\varepsilon\nu \sim h\nu,$$

that leads to *the analogy* $\varepsilon \sim h$ connecting Schrödinger's idea with the idea of auto-oscillations. So we come to the following heuristic statements:

− the quantum systems are characterized by inflow and dissipation not energy, but *action*;

− it is necessary to consider as energy of quantum auto-oscillation the average action inflow rate: $E = \overline{dS_0/dt}$;

− the action inflow for the period is always equal to the Planck constant $h = 2\pi\hbar = \Delta S_0 > 0$.

Here the last statement seems to be the most informative because it points to important feature of the quantum auto-oscillations. This feature, certainly, should have a general mathematical reason. This reason is easy for guessing.

We can suppose that *function* $\Phi = S_0/\hbar$ behaves like a polar angle on a plane, i.e. is *a multivalued* function, with periods, multiple $2\pi$. As a result we come to the following heuristic statement:

*there is a multivalued function $S_0$ having the dimensionality of action, defined on a phase space of a system, with periods, multiple Planck constant $h = 2\pi\hbar$, and defining the energy as time average*

$$E = \overline{dS_0/dt}.$$

Obviously, for additivity of such energy it is necessary to assume also the existence of similar property for function $S_0$.

Let's formulate one more hypothesis implicitly assumed till now.

**Hypothesis 4.** *The phase space of system should be symplectic. Attractors should be localized on Lagrange submanifolds.*

It is natural to assume also, that the dynamic system describing the quantum world should be *the close relative of the Hamilton systems.*

Let's set the problem to consider all assumptions stated above and to build the elementary dissipative system, «close to Hamilton one», with use of function of action $S_0$, and with an attractor of described «quantum type».

Let $I, \varphi$ be the variables of action-angle type for one-dimensional harmonic oscillator. We take a Hamiltonian

$$H = \omega_0 I - E_0 \varphi,$$

differing from usual by addition of the *multivalued* potential energy $U = -E_0\varphi$. To avoid unlimited increase of amplitude, we add in the corresponding Hamilton equation a dissipative term[13] $f = -\kappa I$. Then the dynamic equations

$$\begin{cases} dI/dt = -\partial H/\partial \varphi + f = E_0 - \kappa I \\ d\varphi/dt = \partial H/\partial I = \omega_0 \end{cases}$$

have the solution

$$\begin{cases} I(t) = E_0/\kappa + e^{-\kappa t}\left(I(0) - E_0/\kappa\right) \\ \varphi(t) = \varphi(0) + \omega_0 t \end{cases}$$

which shows, that oscillator converges to a limit cycle $C$, with the circular frequency $\omega = \omega_0$ and amplitude $I = E_0/\kappa$. In this case the energy inflow is provided by the potential energy $U = -E_0\varphi$, and a role of multivalued action function plays $S_0 = \hbar\varphi$. Therefore the «energy», by new definition, is a quantity

$$E = dS_0/dt = \hbar\omega_0.$$

Let's accept now, that $E_0 = \hbar\kappa$. Then the value of usual Hamiltonian $H_0 = \omega_0 I$ on a limit cycle $C$ will coincide with the «quantum» value of energy $\hbar\omega_0$:

$$E = H_0\big|_C = \hbar\omega_0,$$

---

[13] The reason for such choice a dissipative term will be clear later.

and, that is especially remarkable, the motion on a cycle will coincide with what is set by $H_0$. The possibility of such coincidence hints a hypothesis making the second definition of energy which we accept as a final one[14].

**Hypothesis 5.** *Dynamics on smooth «quantum attractors» $A$ is defined by the classic[15] Hamiltonian $H_0$. The energy spectrum is defined as the set of values of $H_0$ on connected components of $A$.*

The correctness of given formulation is guaranteed by Hypothesis 4, owing to a known fact of symplectic geometry:

*if the Hamilton field of some function $F$ touches of isotropic submanifold $L \subset M$, then this function is locally constant on $L$. If $L$ is a Lagrangian submanifold, the contrary is true.*

The example constructed above admits natural and extremely wide generalization. Really, the auto-oscillations in usual dissipative systems arise as a kind of dynamic balance between a pumping, i.e. the energy inflow from the outside, and its losses, dissipations. In the given example the pumping is provided by multivalued action function. In turn, a source of such functions is the multiconnectedness of a phase space[16], its general topological property. Therefore there is a tempting thought to use the *multivalued* Hamilton functions for construction of a special class of dynamic systems capable to auto-oscillations. For realization of this plan it would be desirable to find also the general mathematical mechanism playing a role similar to energy dissipation in usual auto-oscillatory systems.

There is a simple, and somewhat, a unique solution to this problem. Note that the phase flow of the oscillator belongs to the transformations of conformally-symplectic type, i.e. its influence on the symplectic form is reduced to its multiplication by decreasing exponential factor:

$$dI(t) \wedge d\varphi(t) = e^{-\kappa t} dI(0) \wedge d\varphi(0).$$

The solution consists in that this property of conformal contraction of the symplectic form *to postulate*. It should provide the desirable localization of q-states and attractors on Lagrange submanifolds. The corresponding dynamic systems will be called *conformally-dissipative*. It is that class of systems which is necessary to us.

Certainly, only the dynamically non-trivial systems of the given kind are interesting. It means that the steady regimes of such systems must differ from the set of stationary points. We will suppose that this property should be provided by the multivalued function $S_0$.

---

[14] Later we will find out, in what cases this definition coincides with the previous definition of energy in terms of frequency.

[15] The exact definition of a «classic Hamiltonian» is a problem for future.

[16] In the example above, it is a plane without the origin or the cylinder, if to admit the negative values of action *I*.

# Chapter 2. Elements of mathematics of CD-systems.

## §1. Definition and some general properties of CD-systems.

Let $V(x)$ be the vector field on symplectic manifold $(M,\Omega)$ generating *a semiflow* (one-parametrical semigroup of smooth mappings) $g_t : M \to M$, $t \geq 0$, conformally contracting the symplectic form $\Omega$:

$$(1.1) \qquad g_t^*\Omega = e^{-\kappa t}\Omega, \quad \kappa = const > 0.$$

Such dynamic system will be called the CD-system. In non-autonomous case (it is considered, that symplectic *form does not depend on time*) the role of semiflow plays the evolutionary system[17] $\{g_{t_2,t_1},\ t_2 \geq t_1\}$. In this case (1.1) is replaced on

$$(1.2) \qquad g_{t_2,t_1}^*\Omega = e^{-\kappa(t_2-t_1)}\Omega.$$

The infinitesimal variant of definition of CD-system is reduced to the formula:

$$(1.3) \qquad L_V\Omega = -\kappa\Omega.$$

Here $L_V$ denotes the differentiation along a vector field $V$. It is easy to see, that symplectic form $\Omega$ is always exact: $\Omega = d\alpha$, and the potential $\alpha$ is defined canonically:

$$\alpha = -\kappa^{-1}i_V\Omega.$$

Here $i_V$ is the operator of internal multiplication of the form by a vector $V$.

On the contrary, any 1-form $\alpha$ with non-degenerate external derivative $\Omega = d\alpha$ uniquely defines a vector field $V$ of some CD-system:

$$(1.4) \qquad i_V d\alpha = -\kappa\alpha.$$

Using the standard operator $\mathrm{J}$ transforming 1-forms $\omega$ to vector fields [1]:

$$X = \mathrm{J}\omega \Leftrightarrow i_X\Omega = -\omega,$$

the dynamic field of CD-system can be written down in an explicit form

$$(1.5) \qquad V = \kappa \cdot \mathrm{J}\alpha.$$

Comparing it with the definition of Hamilton fields $\mathrm{J}dH$, we obtain a rule: to obtain *CD-system field, it is necessary in expression for the general Hamilton fields to replace the partial derivatives $\partial H / \partial x^i$ to $\kappa\alpha_i$, where $\alpha_i$ are the coefficients of the form $\alpha = \sum \alpha_i dx^i$*.

---

[17] It means, that $g_{t_3,t_2} \circ g_{t_2,t_1} = g_{t_3,t_1}$ for any $t_3 \geq t_2 \geq t_1$.

We call the form $\alpha$ *the canonical defining form*, and will consider the formula (1.4) as the initial definition of CD-system, supposing also the dependence of the form $\alpha$ on time. But at first we will concentrate on the autonomous case.

Thus, the system is set by a pair $(M,\alpha)$, unlike of Hamilton systems which are set by a three $(M,\Omega,H)$. It should be noted, that the dynamic vector field $V$ is defined invariantly by a covector field $\alpha$. In particular, if $f$ is a diffeomorphism then a substitution $x = f(y)$ gives for $y$ the equations of CD-type with the defining form $f^*\alpha$.

The volume form is exact: $\Omega^n = d(\alpha \wedge \Omega^{n-1})$, so a phase space of CD-system is always *non-compact*. Besides, physically are interesting *multiconnected* $M$ with nonzero cohomology group: $H^1(M) \neq 0$. Natural examples of such manifolds give the complements of complex hypersurfaces in the Kähler manifolds, in particular, in the complex projective space. Other example is the cotangent bundles of multiconnected manifolds $M = T^*N$.

To system $(M,\alpha)$ corresponds the invariant distribution of hyperplanes $\ker \alpha$ defined everywhere out of the set $Nul(\alpha) = \{x \in M | \alpha_i(x) = 0\}$. Obviously, $Nul(\alpha)$ coincides with the set of stationary points. The distribution of hyperplanes are the symplectic compliment of dynamic vectors $V$: $\ker \alpha = V^\perp$. Analogue of this distribution in Hamiltonian mechanics is the integrable distribution $\ker dH$ corresponding to the foliation of level hypersurfaces of $H$.

Now we give the list of the basic formal consequences of the definition (1.4). They look like identities for the form $\alpha$, semiflow $g_t$, dynamic field $V$, any Hamilton field $\mathrm{J}\,dF$, and Poisson brackets. To shorten the notation, we will use $L^\kappa_X$ for the extended Lie derivative:

$$L^\kappa_X = L_X + \kappa,$$

and also $X.f$ for the Lie derivative $L_X f$ of a function $f$.

**Theorem.** *The algebraic consequences of definition* (1.4) *look like the following identities (formulas* (1.6) *and* (1.14) *relate to the autonomous case)*:

(1.6) $$g_t^* \alpha = e^{-\kappa t} \alpha,$$

(1.7) $$L^\kappa_V \alpha = 0,$$

(1.8) $$\alpha(V) = 0,$$

(1.9) $$[V, \mathrm{J}\,dF] = \mathrm{J}\,d(L^\kappa_V F),$$

(1.10) $$F - \alpha(\mathrm{J}\,dF) = \kappa^{-1} L^\kappa_V F, \quad \text{or} \quad V.F = -\kappa \alpha(\mathrm{J}\,dF),$$

(1.11) $$L_{JdF}\alpha = -\kappa^{-1}d\left(L_V^\kappa F\right),$$

(1.12) $$[L_V^\kappa, i_W] = i_{[V,W]},$$

(1.13) $$L_V^\kappa\{F_1, F_2\} = \{L_V^\kappa F_1, F_2\} + \{F_1, L_V^\kappa F_2\},$$

(1.14) $$g_t^*\{F_1, F_2\} = e^{\kappa t}\{g_t^* F_1, g_t^* F_2\},$$

(1.15) $$L_V^\kappa \alpha(W) = \alpha([V,W]).$$

Note, that any of equalities (1.7), (1.6) also can be used as a defining relation of CD-system because the form $\alpha' = \alpha + \kappa^{-1}d\alpha(V)$ satisfies (1.4). Then, equality (1.13) shows that the operator $L_V^\kappa$ is the differentiation of Poisson algebra. So, it is possible to use (1.13) as definition of conformally-Poisson systems.

**Examples.**

**1.** In the case $\alpha = \mathbf{p}d\mathbf{q}$ a vector field $V$ of CD-system is proportional to the Euler field:

$$V = -\kappa \mathbf{p} \cdot \partial/\partial \mathbf{p}.$$

**2.** For the oscillator of Chapter 1:

$$\alpha = (I - \hbar)d\varphi + \kappa^{-1}d(I\omega_0) = Id\varphi + \kappa^{-1}dH - dS_0.$$

We see that the action form $Id\varphi$, Hamilton function $H = I\omega_0$, and action function $S_0 = \hbar\varphi$ unite in the whole mathematical object.

**3.** If $x(t)$ is a phase trajectory of system $(M, \alpha)$ and $f: M \to M$ is a symplectomorfism:

$$f_t^* \Omega = \Omega,$$

then $y(t) = f^{-1}(x(t))$ is a phase trajectory of system $(M, \alpha + dF = f^*\alpha)$.

## §2. The correspondence with Hamilton equations.

At $\kappa = 0$ the formula (1.3) characterises Hamilton systems. Consider $\kappa$ as a parameter on which the form $\alpha$ and the vector field $V$ depends. Suppose, that the symplectic form $d\alpha(\kappa) = \Omega(\kappa)$ and the dynamic field $V(\kappa)$ is analytic at $\kappa = 0$. Then a defining relation (1.4) gives that $\alpha(\kappa)$ looks like

(2.1) $$\alpha(\kappa) = \alpha_0(\kappa) + \kappa^{-1}dH,$$

where the form $\alpha_0(\kappa)$ is analytic at zero, and $H$ (may be, $H$ is multivalued) does not depend on $\kappa$. Therefore the system field $V$ is a linear combination:

(2.2) $$V(\kappa) = \kappa \cdot J \alpha_0(\kappa) + J \, dH.$$

At the limit $\kappa = 0$ it is obtained the Hamiltonian system $(M, \Omega(0), H)$. In particular, the Hamilton's canonical equations are a limiting case of the equations of CD-system $(M, \alpha = \mathbf{p} d\mathbf{q} + \kappa^{-1} dH)$:

(2.3) $$\begin{cases} d\mathbf{p}/dt = -\partial H/\partial \mathbf{q} - \kappa \mathbf{p} \\ d\mathbf{q}/dt = \partial H/\partial \mathbf{p} \end{cases}$$

Equations (2.3) are, possibly, an elementary dissipative generalisation of the Hamilton equations. However the other form of these equations obviously reflecting the idea of «pumping» by means of *multivalued* Hamiltonian is more important here. Therefore we rewrite the full Hamiltonian in the form of difference $H - \kappa S_0$, selecting the singlevalued part $H$ and the multivalued action function $S_0$. More strictly $S_0$ is a local potential of the closed form $dS_0$ representing a nonzero cohomology class $H^1(M)$. However, sometimes it makes sense to consider singlevalued[18] action functions $S_0$.

So, we get a system

(2.4) $$\alpha = \mathbf{p} d\mathbf{q} + \kappa^{-1} dH(\mathbf{p}, \mathbf{q}, t) - dS_0(\mathbf{q}, t)$$

with equations of motion

(2.5) $$\begin{cases} d\mathbf{p}/dt = -\partial H/\partial \mathbf{q} - \kappa \cdot (\mathbf{p} - \partial S_0/\partial \mathbf{q}), \\ d\mathbf{q}/dt = \partial H/\partial \mathbf{p}. \end{cases}$$

## §3. The relation with Lagrange equations.

Let $L(\mathbf{q}, \dot{\mathbf{q}}, t)$ be a Lagrangian, obtained from $H$ by means of Legendre transformation: $L = \mathbf{p}\dot{\mathbf{q}} - H$, $\dot{\mathbf{q}} = \partial H/\partial \mathbf{p}$. Then the dynamic equation (2.5) are equivalent to *usual Lagrange equations* for any of two Lagrangians (note that $L''$ is multivalued in general):

(3.1) $$L' = e^{\kappa t}(L - dS_0/dt) = e^{\kappa t}(L - \dot{\mathbf{q}} \cdot \partial S_0/\partial \mathbf{q} - \partial S_0/\partial t),$$

(3.2) $$L'' = e^{\kappa t}(L + \kappa S_0).$$

In turn these equations are written as generalized Lagrange equations for $L$ obtained by addition of a «quantum force» $\mathbf{f}$ to right member of equations:

(3.3) $$\frac{d}{dt}\frac{\partial L}{\partial \dot{\mathbf{q}}} - \frac{\partial L}{\partial \mathbf{q}} = -\kappa \cdot \left(\frac{\partial L}{\partial \dot{\mathbf{q}}} - \frac{\partial S_0}{\partial \mathbf{q}}\right) = \mathbf{f}.$$

---

[18] The multivaluedness can disappear as a result of reduction of a system by the symmetry group.

We can consider formulas (3.1), (3.2) as a way of construction of CD-systems by the known *classic* Lagrangian $L$ and by the multivalued «nonclassic» summand $\kappa S_0$.

In this place we return to physics and consider two questions, concerning concepts of classic Lagrangian and initial attempt to determine energy through frequency (see Chapter 1).

### 3.1. What is a «classic Lagrangian»?

It is known, that in classic dynamics the Lagrangian is defined to within addition of total time derivative $df(\mathbf{q},t)/dt$. In our case it is already incorrect, and Lagrangians $L$, equivalent in classic sense, now lead *to different* systems. Obviously, this uncertainty is formally reduced to arbitrariness in the selection of function $S_0$.

But also within the limits of already fixed dynamic system there is a question what to consider as a classic Lagrangian[19], because the replacements are possible

(3.4) $$\tilde{L} = L - \kappa \cdot \sigma(\mathbf{q},t), \ \tilde{S}_0 = S_0 + \sigma(\mathbf{q},t),$$

still saving the form (3.3) of dynamic equations.

The main criterion which should eliminate the given arbitrariness is dictated by simple physical reasons (unfortunately, this criterion is unsufficiently accurate). Really, «a classical limit» in the sense of quantum theory is natural to connect with equality $\mathbf{f} \approx \mathbf{0}$, because then equations (3.3) turn to usual classic Lagrangian equations. But the same equality can be considered also as a quantization condition in quasiclassic approximation. Therefore it is natural to expect this equality only for the *steady* motion. Thus, the criterion of selection a «true» classic Lagrangian consists in the minimisation of a «quantum force» for the steady motions, by means of replacements (3.4):

(3.5) $$\tilde{\mathbf{f}} = \mathbf{f} + \kappa \cdot \partial \sigma / \partial \mathbf{q}.$$

From here, by the way, follows that unlike the modern quantum theory in future CD-theory (if it is possible), a classical Lagrangian (and Hamiltonian) is not obliged to be set a priori. In general situation it should be calculated, being based on the analysis of steady motions of a system with, probably, «not physical» primary Lagrangian.

### 3.2. Whether it is possible to determine energy through frequency?

Consider this question on an example of autonomous system (2.4), with the equations of motion (2.5). Assume, that the attractor of system belongs to the Lagrangian manifold

$$\mathbf{p} = \partial S_0 / \partial \mathbf{q},$$

and a Hamiltonian $H$ is classic. It means, that on every connected component of attractor the stationary Hamilton-Jacobi equation is satisfied (the total action looks like $S(\mathbf{q},t) = S_0(\mathbf{q}) - Et$):

---

[19] The similar remark also concerns to a Hamiltonian.

$$H(\partial S_0 / \partial \mathbf{q}, \mathbf{q}) = E,$$

and the constant $E$ has sense of energy, according to Hypothesis 5. It is required to find out, in what cases this energy can coincide with time average of a derivative $dS_0 / dt$. We have:

$$\frac{dS_0}{dt} = \frac{\partial S_0}{\partial \mathbf{q}} \cdot \dot{\mathbf{q}} = \mathbf{p} \cdot \frac{\partial H}{\partial \mathbf{p}}.$$

Therefore for this coincidence the equality is sufficient:

$$\mathbf{p} \cdot \frac{\partial H}{\partial \mathbf{p}} = H,$$

claiming that a classic Hamiltonian *is homogeneous first degree function of momentums*.

**Examples.**

**1.** Let $q$ be an arc coordinate on a circle, and let a Lagrangian be $L' = e^{\kappa t}\left(m\dot{q}^2/2 - h\dot{q}\right)$. This system converges to the limit cycle $m\dot{q} = h$. Here a function $S_0 = hq$ provides the *quantization* of a classic motion of Lagrangian system with $L = m\dot{q}^2/2$.

**2.** Let's consider a free relativistic particle with a fixed momentum $\mathbf{b}$:

$$H = c\sqrt{\mathbf{p}^2 + m^2 c^2}, \quad S(\mathbf{q}, t) = \mathbf{b}\mathbf{q} - ct\sqrt{\mathbf{b}^2 + m^2 c^2} = S_0(\mathbf{q}) - Et.$$

Then the value $dS_0 / dt$ on attractor $\mathbf{p} = \mathbf{b}$ coincides with the energy only at $m = 0$ when the Hamiltonian is a homogeneous function of momentum.

## §4. The compound systems.

As is known, a direct product of dynamic systems is described by the direct sum of vectorial fields $V = V_1 + V_2$ on a direct product of phase spaces $M = M_1 \times M_2$. To define a direct product of CD-systems, by analogy to the case of Hamiltonian systems it is natural to postulate the additivity of symplectic forms:

(4.1) $$\Omega = \Omega_1 + \Omega_2.$$

But then the sum $V = V_1 + V_2$ will generate CD-system if only constants $\kappa$ for various parts of system wil *coincide*. Finally, it leads to the important conclusion that application of systems of the given type for the description of physical subsystems means that this constant should be considered as *a fundamental dimensional constant*, like a speed of light and the Planck constant. Then a direct product of CD-systems will have a simple form:

(4.2) $$(M_1, \alpha_1) \times (M_2, \alpha_2) = (M_1 \times M_2, \alpha_1 + \alpha_2).$$

The compound systems *interacting* through Hamiltonian $U$ (probably, multivalued) have a similar kind. In this case the defining form also breaks up into the sum:

(4.3) $\qquad \alpha = \alpha_1 + \alpha_2 + \kappa^{-1}dU = (\alpha_1 + \kappa^{-1}d_1U) + (\alpha_2 + \kappa^{-1}d_2U) = \beta_1 + \beta_2.$

Here the decomposition of differential in the sum of partial differentials is meant

$$dU = d_1U + d_2U.$$

The equations of motion of such a system break up into a pair of equations for subsystems:

(4.4) $\qquad i_{V_1}d_1\beta_1 = -\kappa\,\beta_1, \quad i_{V_2}d_2\beta_2 = -\kappa\,\beta_2.$

The interaction does not allow to consider separate subsystems as true dynamic systems because for them, generally speaking, there is no the concept of a phase flow. Nevertheless, there is an important approach when influence of a «small» subsystem on «big» one can be neglected. Then to each dynamic trajectory of a «big» subsystem there corresponds a time-dependent phase flow in a «small» subsystem. That is the main physical sense of the concept of a non-autonomous dynamic system[20]. So, the relation (4.4) justifies the formal definition of *non-autonomous CD-systems* being made above.

**Example.**

Let's consider the compound system with the defining form:

$$\alpha = \mathbf{p}d\mathbf{q} + Id\varphi + \kappa^{-1}dH,$$

where $(\mathbf{p},\mathbf{q}) \in \mathbb{R}^{2n} = M_1$, $(I,e^{i\varphi}) \in \mathbb{R} \times S^1 = M_2$, $H = \mathbf{p}^2/2m + U(\mathbf{q}) + I\omega + \mathbf{fq}\cdot\cos\varphi$.

The equations of this system

(4.5) $\qquad \begin{cases} \dot{\mathbf{p}} = -\kappa\mathbf{p} - \partial U/\partial\mathbf{q} - \mathbf{f}\cdot\cos\varphi, & \dot{\mathbf{q}} = \mathbf{p}/m \\ \dot{I} = -\kappa I + \mathbf{fq}\cdot\sin\varphi, & \dot{\varphi} = \omega \end{cases}$

show, that the $\mathbb{R}^{2n}$-subsystem obeys the equation of the second order

(4.6) $\qquad \ddot{\mathbf{q}} + \kappa\dot{\mathbf{q}} + m^{-1}\cdot\partial U/\partial\mathbf{q} = -m^{-1}\cos(\omega t + \varphi_0)\cdot\mathbf{f},$

which coincides with the equation of forced oscillations of mechanic system with a friction. In case of quadratic potential energy the system $U(\mathbf{q})$ converges to single limit cycle with a frequency $\omega$.

---

[20] The non-autonomous systems can arise also in another way. For example, the equations in variations along the solution $p(t), q(t)$ to autonomous system $\alpha = pdq + \kappa^{-1}dH$ can be considered as the equations of motion of the system $\alpha_1 = p_1dq_1 + \kappa^{-1}dH_1(p_1,q_1,t)$, with a Hamiltonian $H_1 = \frac{1}{2}\left(H_{pp}p_1^2 + 2H_{pq}p_1q_1 + H_{qq}q_1^2\right)$.

## §5. Systems on Kähler manifolds.

Let $U$ be a global potential of the Kähler metrics on a complex manifold $M$. Then we can take the defining form of CD-system as

(5.1) $$\alpha = \operatorname{Im} \partial U,$$

where $\partial$ is a holomorphic exterior derivative. Such systems are the gradient ones. Then, we can add the differential of the Hamilton function to the defining form:

(5.2) $$\alpha = \operatorname{Im} \partial U + \kappa^{-1} dH.$$

Thus, the obtained system is set by two scalar functions on a complex manifold. In particular, the induced CD-system on a complex submanifold is set by restrictions of these functions. Note, that a potential and Hamiltonian are determined not uniquely because the transformations are possible

$$U' = U + \operatorname{Re} F, \quad H' = H + \varepsilon \cdot \operatorname{Im} F, \quad \varepsilon = \kappa/2,$$

with holomorphic functions $F$.

Often happens so that the phase space is constructed as *complement of a set of singularities* of local potential on the compact Kähler manifold. Then for the correctness of corresponding dynamic system it is necessary, that these singularities were *repelling*. For this purpose it is sufficient, that at approach to singularities the potential becomes positively infinite, and the Hamiltonian field remains regular:

(5.3) $$U \to +\infty, \quad \|\operatorname{J} dH\| < const.$$

Let's write the general dynamic equations of the Kähler systems. Let $U^{i\bar{k}}$ be a reciprocal matrix to the Hermite matrix of second derivatives $U_{i\bar{k}} = \partial^2 U / \partial x^i \partial \bar{x}^k$. Then the equations look like:

(5.4) $$dx^i / dt = \sum U^{i\bar{k}} \left( i \frac{\partial H}{\partial \bar{x}^k} - \varepsilon \frac{\partial U}{\partial \bar{x}^k} \right), \quad \varepsilon = \kappa/2,$$

or in non-coordinate form,

(5.5) $$dx / dt = \operatorname{J} dH - \varepsilon \operatorname{grad} U,$$

where the gradient is understood relatively to Riemannian metrics

(5.6) $$\|dx\|^2 = 2 \sum U_{i\bar{k}} dx^i d\bar{x}^k.$$

Note, that for pluriharmonic (or zero) Hamiltonian

$$H = \varepsilon \cdot \operatorname{Im} F, \quad \bar{\partial} F = 0,$$

the system (5.5) turns in pure gradient, with potential

$$U' = U + \operatorname{Re} F.$$

In the case of a single-valued potential $U'$ almost all trajectories of gradient systems are attracted by the stationary points of its local minima. Thus, the non-trivial attractor can arise only for multivalued potential.

If the system has a symmetry then usually critical points of potential form *a manifold*. To such situations is applicable

**Theorem.** *Let the set of critical points $\operatorname{Nul}(dU)$ of a potential of system* (5.2) *be a manifold. Then*

a) $\operatorname{Nul}(dU)$ *is symplecticly isotropic;*

b) *the law of evolution of a potential looks like:*

(5.7) $$\frac{dU}{dt} = \{H, U\} - \varepsilon \|\operatorname{grad} U\|^2, \quad \|\operatorname{grad} U\|^2 = 2 \sum U^{i\bar{k}} U_i U_{\bar{k}} ;$$

c) if $\{H, U\} = 0$, and a potential $U$ is single-valued, then $\operatorname{Nul}(dU)$ is an attractor, with Hamiltonian flow on it:

(5.8) $$dx/dt = \mathrm{J}\, dH, \quad x \in \operatorname{Nul}(dU).$$

**Example.**

The Kähler system with zero Hamilton function and potential $U = |z|^2 + 2\operatorname{Re}(c \log z)$, where $c \in \mathbb{C}$, is regular at $\operatorname{Re} c < 0$, and in the case $\operatorname{Im} c \neq 0$ converges to a limit cycle $|z| = \sqrt{-\operatorname{Re} c}$.

## §6. Description in the expanded space.

### 6.1. Defining form in the expanded space.

It is convenient to consider non-autonomous systems in the expanded phase space including additional time dimension:

(6.1) $$\bar{M}^{2n+1} = M^{2n} \times \mathbb{R}.$$

Let's agree about notations. The objects in the expanded space will be denoted by the barred symbols. For example, we will denote $\bar{V}$ the dynamic field on $\bar{M}$:

$$\bar{V} = V + \partial/\partial t.$$

Similarly, we will denote $\bar{g}_\tau$ the dynamic semigroup $\exp(\tau \bar{V})$ acting in the expanded space:

(6.2) $$\bar{g}_\tau(x,t) = (g_{t+\tau, t} x, \, t+\tau).$$

Sometimes, to emphasise that the differential is considered in usual phase space $M$ instead of $df$ we will write $d_M f$. It is convenient to use also a special denotation for the elongated operator of exterior derivative:

$$\bar{d} f = e^{-\kappa t} \cdot d(e^{\kappa t} f) = df + \kappa dt \wedge f.$$

In expanded space the role of defining form will play the 1-form

(6.3) $$\theta = \alpha_1 - H dt,$$

where $\alpha_1$ is some symplectic potential:

$$d_M \alpha_1 = \Omega.$$

As always, we will assume, that *symplectic form does not depend on time*. It is equivalent to the assumption of closedness of the form $\dot{\alpha}_1 = \partial \alpha_1 / \partial t$:

$$d_M \dot{\alpha}_1 = 0.$$

It is known [1], that form (6.3) generates a Hamiltonian system in $\bar{M}$. The corresponding vector field $\bar{W} = W + \partial / \partial t$ is uniquely determined by relations

(6.4) $$i_{\bar{W}} d\theta = 0, \quad dt(\bar{W}) = 1,$$

and for isochronous components of this field we get:

$$W = \mathrm{J}(dH + \dot{\alpha}_1).$$

It appears that the form $\theta$ as well generates a CD-system.

**Theorem 1.** a) *The vector field $\bar{V} = V + \partial / \partial t$ in expanded space, defined by the relations*

(6.5) $$i_{\bar{V}} \bar{d} \theta = 0, \quad dt(\bar{V}) = 1,$$

*is the field of CD-system $(M, \alpha)$ with canonical form $\alpha = \kappa^{-1}(dH + \dot{\alpha}_1) + \alpha_1$. This field is represented in the form of the sum of two fields:*

(6.6) $$\bar{V} = \bar{W} + \bar{R}.$$

*The first field is Hamiltonian one:* $\overline{W} = W + \partial/\partial t = \mathrm{J}(dH + \dot{\alpha}_1) + \partial/\partial t$. *The second is an isochronous dissipative field:* $\overline{R} = \kappa \cdot \mathrm{J}\alpha_1, \quad dt(\overline{R}) = 0$.

b) *Conversly, if the defining form* $\alpha$, $d\dot{\alpha} = 0$ *is known, then as the form* $\theta$ *it is possible to take the form with a zero Hamiltonian* $H$ :

(6.7) $$\theta = \alpha(t) - \int_{t_0}^{t} e^{-\kappa(t-\tau)} \dot{\alpha}(\tau) d\tau.$$

*In the autonomous case it is possible to accept* $\theta = \alpha$.

Such (non-autonomous, generally) system in the expanded space further will be denoted by $(\overline{M}, \theta)$.

**Corollary**. *The system* $(\overline{M}, \theta)$ *admits a relative integral invariant of* Poincaré-Cartan

(6.8) $$\oint e^{\kappa t}\theta = const.$$

Note a fundamental role of the 2-form $\overline{d}\theta$ determining a system. It is easy to see, that it has a rank $2n$, is $\overline{d}$-closed, expressed through the form $\alpha$ :

$$\overline{d}\theta = \Omega - \kappa\alpha \wedge dt = d_M \alpha - \kappa\alpha \wedge dt,$$

and is conformally transformed by dynamic semigroup:

$$\overline{g}_\tau^* \overline{d}\theta = e^{-\kappa\tau} \overline{d}\theta.$$

Unlike the form $\overline{d}\theta$ and canonical form $\alpha$, the form $\theta$ is defined non-uniquely because the gauge transformations are possible

(6.9) $$\theta' = \theta + \overline{d}\lambda,$$

which preserv a dynamic field $\overline{V}$. This gauge nonuniqueness will be used further as means for selection of the «physical» Hamiltonian (for it definition see §10).

Note also, that the transformations $\theta' = \theta + f(t)dt$ don't change both the fields $\overline{V}, \overline{W}$.

### 6.2. The features of standard gauge.

Here we will consider some features of «standard» gauge

(6.10) $$\theta(\overline{V}) = 0.$$

This gauge also equivalent to

(6.11) $$\theta(\overline{W}) = 0.$$

Note, that this gauge is automatically satisfied in the autonomous case if to accept $\theta = \alpha$.

About the degree of «rigidity» of standard gauge speaks

**Theorem 2.** *Residual gauge transformations* (6.9), *preserving a standard gauge* $\theta(\overline{V}) = 0$, *are generated by solutions to the equation* $L_{\overline{V}}^{\kappa} \lambda = 0$.

Other properties of standard gauge are described by following theorems.

**Theorem 3.** *Fields* $\overline{V}$ *and* $\overline{W}$ *commute in only case when a function* $\theta(\overline{V})$ *does not depend on a point* $x \in M$:

(6.12) $$\theta(\overline{V}) = f(t).$$

Note, that transition from gauge (6.12) to the standard one is trivial.

**Theorem 4.** *If denote* $\overline{g}_\tau = \exp(\tau \overline{V})$, $\overline{h}_\tau = \exp(\tau \overline{W})$, $\overline{r}_\tau = \exp(\tau \overline{R})$. *If* $\theta(\overline{V}) = 0$, *then:*

a) $\overline{g}_\tau = \overline{r}_\tau \circ \overline{h}_\tau = \overline{h}_\tau \circ \overline{r}_\tau$;

b) $\overline{g}_\tau^* \theta = \overline{r}_\tau^* \theta = e^{-\kappa\tau} \theta$, $\overline{h}_\tau^* \theta = \theta$;

c) $i_{\overline{V}} d\theta = -\kappa \theta$;

d) *the form* $\theta$ *admits the representation:*

(6.13) $$\theta = P_{t_0}^* \alpha_1(t_0),$$

*where* $P_{t_0} = P_{t_0} \circ \overline{h}_\tau$ *denote a projection mapping along the flow* $\overline{h}_\tau$ *on a hyperplane* $t = t_0$:

$$P_{t_0}(x,t) = \overline{h}_{t_0-t}(x,t) = (h_{t_0,t} x, t_0).$$

*In particular, the form* $\theta$ *can be expressed through* $2n$ *independent functions*[21], *being integrals of Hamilton field* $\overline{W}$.

Thus, in standard gauge the dynamics of system in the expanded space is reduced to «superposition» of usual Hamiltonian dynamics and isochronous dissipative dynamics. By the way, this fact allows to enter the concept of *time-dependent* attractor, as the dynamic trajectory of some «initial» attractor.

---

[21] I.e. the form has a class $2n$, see [4].

Because the relation c) in last theorem looks similar to definition of CD-system it will be useful to indicate the following facts.

**Theorem 5.** *The necessary and sufficient condition for the solvability of the equation*

(6.14) $$i_{\bar{X}} d\theta = -\kappa \theta$$

*is the equality* $\theta(\bar{V}) = 0$.

**Theorem 6.** *Let* $\theta(\bar{V}) = 0$. *Then the general solution to equation (6.14) looks like*

$$\bar{X} = \bar{V} + \mu \cdot \bar{W},$$

*where* $\bar{V} = V + \partial/\partial t$ *is a CD-system field and* $\mu = \mu(x,t)$ *is an arbitrary factor.*

Let's consider now a question of the *realizability* of standard gauge. It is easy to see that the problem is reduced to solution to the equation

(6.15) $$\bar{V}.u = -e^{\kappa t}\theta(\bar{V}) = f$$

for the function $u = u(x,t)$ connected with the gauge function by equality $u = e^{\kappa t}\lambda$. Hence, it is enough to show the solvability of equation (6.15). Unfortunately, it manages to be made only for the domain of the expanded space limited from above on time.

**Theorem 7.** *The solution to the equation*

(6.16) $$\bar{V}.u = f(x,t),$$

*in the domain* $t < t_0$ *of an expanded space* $\bar{M}$, *which satisfies to a «final» condition*

(6.17) $$u\big|_{t=t_0} = u_0,$$

*gives by the formula*

(6.18) $$u(x,t) = u_0\left(\bar{g}_{t_0-t}(x,t)\right) - \int_0^{t_0-t} f\left(\bar{g}_\tau(x,t)\right) d\tau.$$

**Examples.**

**1.** The autonomous Kähler system with potential $U$ and single-valued Hamiltonian $H$ is described in the expanded space by the form $\theta = \text{Im}\,\partial U - H dt$.

**2.** To the form

$$\theta = pdq - H(p,q,t)dt - dS(q,t)$$

there corresponds a canonical form

$$\alpha = pdq + \kappa^{-1} d_M H(p,q,t) - d_M S(q,t).$$

Write down an integral invariant:

$$\oint e^{\kappa t}\theta = \oint e^{\kappa t}(pdq - Hdt - dS) = \oint e^{\kappa t}(pdq + \kappa^{-1}dH - dS) = const.$$

For an isochronous contour of integration we will obtain

$$\oint \alpha = e^{-\kappa t} \cdot const,$$

or

(6.19) $$\oint pdq = \oint d_M S + e^{-\kappa t} \cdot const.$$

Taking account of the hypothesis about the periods of action function $S$, it should be clear that eventually the value of integral in the left part of (6.19) becomes whole multiple of the Planck constant. This fact is possible to understand as indication that the phase flow carries away a contour of integration on the attractor, satisfying the quantization conditions.

**3.** The standard gauge for $\theta = pdq - H(p,q,t)dt$ means, that Hamiltonian $H$ is a homogeneous function of the first degree on momentums.

## §7. Symmetries and integrals of motion.

Because to Lagrangian CD-systems the Noether theorem is applicable, we can be limited here by considering in Hamiltonian framewoks, i.e. in symplectic phase space.

At first we remind the notion of the symmetry of a dynamic system with the dynamic vector field $V$ (see [3]).

**Definition**. The diffeomorphism of a phase space in itself is called the *symmetry* of the vector field $V$ being set on it and the corresponding autonomous differential equation, if it maps a field in itself. In this case a field is said to be invariant relatively to diffeomorphism.

Let $G$ be a symmetry group. Then the group transformations commute with transformations of dynamic semigroup $g_t$:

(7.1) $$g_t \circ h = h \circ g_t, \quad h \in G.$$

The one-parametrical symmetry groups $h_\varepsilon = \exp(\varepsilon X)$ are generated by the vector fields $X$, commuting with a dynamic field $V$:

$$[V, X] = 0.$$

Such fields $X$ are called usually the *symmetry fields*. Also it will be convenient for us to adopt the following definition.

**Definition**. A function $F$ will be called *a quasi-integral of CD-system* if $e^{\kappa t}F$ is a constant of motion:

$$g_t^* F = e^{-\kappa t} F,$$

or, equivalently, if $L_V^\kappa F = 0$.

By means of identities §1 it is easily to derive

**Theorem 1.** a) *To any field of symmetry $X$ corresponds a quasi-integral $F = \alpha(X)$. Conversely, to any quasi-integral $F$ corresponds a Hamiltonian symmetry field $X = \mathrm{J}\, dF$.*

b) *Function F is a quasi-integral iff the equality*

$$F = \alpha(\mathrm{J}\, dF)$$

*is satisfied.*

c) *Hamiltonian symmetry fields are characterized by the property to preserve the defining form:*

$$L_X \alpha = 0.$$

d) *Any symmetry field $X$ uniquely breaks up into the sum of the Hamiltonian $(X_1)$ and non-Hamiltonian $(X_2)$ symmetry fields*

$$X = X_1 + X_2,$$

*and $X_1 = \mathrm{J}\, d\alpha(X)$, $\alpha(X_2) = 0$.*

e) *The symmetry fields $X$ form a Lie algebra. The quasi-integrals form a Lie algebra with respect to Poisson brackets.*

Let's formulate still some facts connected with quasi-integrals. First of all, to every quasi-integral there corresponds an invariant zero level set, if it is not empty.

Also, to finite set of functionally independent quasi-integrals $F_1, F_2, \ldots, F_k$ corresponds, except intersection of zeroes, other invariant manifolds. Namely, out of this intersection, through each point $x_0$ there passes an invariant manifold of codimension $k-1$, with the equation

$$\left(F_1(x) : F_2(x) : \ldots : F_k(x)\right) = \left(F_1(x_0) : F_2(x_0) : \ldots : F_k(x_0)\right) \in P^{k-1}.$$

**Theorem 2.** *The quasi-integral $F$ vanishes on any compact invariant subset $Q = g_t Q \subset M$.*

◀ Let $C = \sup_{x \in Q} |F(x)| < \infty$. Then

$$C = \sup_{g_t x \in Q} |F(g_t x)| = \sup_{x \in Q} |F(g_t x)| = e^{-\kappa t} \sup_{x \in Q} |F(x)| = e^{-\kappa t} C.$$

Hence $C = 0$. ▶

Further, there is a following analogue of reduction theorem of a Hamiltonian systems with symmetry.

**Theorem 3.** Let the defining form $\alpha$ be invariant with respect to some Lie group $G$, acting on $M$:

$$h^* \alpha = \alpha, \quad h \in G.$$

*Then the following statements are valid.*

a) *The action $G$ on $M$ is of Poisson type. The corresponding momentum mapping*

$$P: M \to Lie^*(G)$$

*is set by the formula*

(7.2) $$\langle P(x), \xi \rangle = \alpha(\tilde{\xi}), \quad \tilde{\xi}(x) = \frac{d}{dt}\bigg|_{t=0} e^{t\xi} x$$

*This momentum is a quasi-integral : $L_V^\kappa P = 0$.*

b) *The bundle projection (in the assumption that $M/G$ is a manifold)*

$$\pi: M \to M/G$$

is a Poisson morphism, and the projection of a dynamic field $\pi_* V$ determines on $M/G$ the confomally-Poisson system (see definition §1).

c) *Let $M_0 = P^{-1}(0)$ be the momentum level zero. Then $M_0$ is a $G$-space, and the phase flow on $M_0$ is $G$-invariant. The corresponding factor-system on $M_0/G$ (in the assumption that $M/G$ is a manifold) is a CD-system, and their defining form coincides with the projection of a form $\alpha$ (this projection is uniquely defined).*

As well as in Hamiltonian mechanic, the factor-system on $M_0/G$ will be called *the reduced system*.

**Examples.**

**1.** For a system with the form $\alpha = pdq + \kappa^{-1} d(pq)$ the Hamiltonian $H = pq$ is a quasi-integral. The phase flow looks like $\{p, q\} \to \{e^{-(\kappa+1)t} p, e^t q\}$. The set $H = 0$ breaks up into the union of two invariant coordinate lines. One line is attracting, another – is repelling.

**2.** For a system $\alpha = \mathbf{p}d\mathbf{q}+\kappa^{-1}d(\mathbf{p},\mathbf{W}(\mathbf{q}))$ the Euler field $\mathrm{E}=\mathbf{p}\cdot\partial/\partial\mathbf{p}$ is a non-Hamiltonian symmetry field. Hamiltonian component of this field is equal to quasi-integral

$$\alpha(\mathrm{E})=\kappa^{-1}\cdot(\mathbf{p},\mathbf{W}(\mathbf{q})).$$

**3.** Consider a magnetic monopole, fixed in a plane orthogonal to axis of a constant current, with a friction being proportional to speed. Such system is described by the form

$$\alpha = \mathbf{p}d\mathbf{q}+\kappa^{-1}dH - hd\varphi,$$

where $\varphi$ is a polar angle in $\mathbf{q}$-plane, and $H=\mathbf{p}^2/2m$. Passing to polar coordinates $r,\varphi,p_r,p_\varphi$:

$$q_1+iq_2 = re^{i\varphi},\ p_1+ip_2 = \left(p_r+ip_\varphi/r\right)e^{i\varphi},$$

we obtain

$$\alpha = p_r dr + \kappa^{-1}d\left(p_r^2+p_\varphi^2/r^2\right)/2m+\left(p_\varphi-h\right)d\varphi.$$

To the symmetry $\varphi\to\varphi+\varepsilon$ there corresponds the quasi-integral $Q=p_\varphi - h$. The reduced system is described by the form

(7.3) $$\tilde{\alpha} = p_r dr + \kappa^{-1}d\tilde{H},$$

with the reduced Hamiltonian $\tilde{H} = p_r^2/2m + h^2/2mr^2$.

The reduced Poisson structure is set by Poisson brackets of $U(1)$-invariant functions depending on $p_r, r, p_\varphi$:

(7.4) $$\{f,g\} = \partial f/\partial p_r \cdot \partial g/\partial r - \partial g/\partial p_r \cdot \partial f/\partial r,$$

The dynamic field of this confomally-Poisson system looks like

(7.5) $$\tilde{V} = \{\tilde{H},\cdot\} - \kappa p_r \cdot \partial/\partial p_r,$$

and $L_{\tilde{V}}^\kappa = L_{\tilde{V}} + \kappa$ is the differentiation of Lie algebra consisting of functions $f = f(p_r, r, p_\varphi)$.

## §8. The solutions of Lie's type.

It can happen, that some solutions $x(t)$ to system $(M,\alpha)$ coincide with solutions to the Hamiltonian system $(M,\Omega=d\alpha, H)$. Such solutions will be called the *Hamilton solutions*. The following theorem is convenient for search of Hamilton solutions to autonomous CD-systems with symmetry.

**Theorem 1.** *Let the defining form $\alpha$ be invariant with respect to the Lie group G, acting on M:*

(8.1) $$h^*\alpha = \alpha, \ h \in G,$$

*and let $x(t)$ be the solution to dynamic equations of this system. Let's make time-dependent change of variables (passage to the «moving coordinate system»):*

(8.2) $$x(t) = \exp(t\xi)y(t),$$

*where $\xi$ is some element of the Lie algebra $Lie(G)$. Then $y(t)$ satisfies the dynamic equations of an autonomous system $(M, \alpha_Q)$ with the defining form*

(8.3) $$\alpha_Q = \alpha - \kappa^{-1}dQ,$$

*where $Q = \alpha(\tilde{\xi})$ is a quasi-integral, corresponding to a fundamental field*

(8.4) $$\tilde{\xi}(y) = d/dt\big|_{t=0} \exp(t\xi)y = \mathbf{J}\, dQ.$$

*This equations looks like*:

(8.5) $$dy/dt = \kappa \mathbf{J}\alpha_Q = V(y) - \mathbf{J}\, dQ.$$

Further the solutions $x(t) = \exp(t\xi)x_0$ will be called a *Lie solution*. Obviously, in conditions of theorem 1 every Lie solution is at the same time a Hamilton solution, with a Hamiltonian $Q$. The reason of interest to Lie solutions is that they often are placed on attractors.

**Corollaries.**

a) *For $\alpha = \alpha_0 + \kappa^{-1}dH$ the passage to a moving coordinate system is equivalent to shift of a Hamiltonian: $H \to H - Q$.*

b) *The set of all points $x(t) = \exp(t\xi)x_0$ provided a fixed $\xi \in Lie(G)$, belonging to Lie solutions to system $(M, \alpha)$, coincides with the set $Nul(\alpha_Q)$ of stationary points of moving system $(M, \alpha_Q)$. With respect to a primary system $(M, \alpha)$ the set $Nul(\alpha_Q)$ is invariant, symplecticly isotropic (provided this set is a manifold), and also carrying the dynamics of a Hamiltonian system $(M, d\alpha, Q)$. The restriction of a quasi-integral $Q$ on $Nul(\alpha_Q)$ vanishes.*

*The points $x_0 \in Nul(\alpha_Q)$ is defined by the equation $\tilde{\xi}(x_0) = V(x_0)$ or, equivalently, by:*

(8.6) $$\alpha\big|_{T_{x_0}M} = \kappa^{-1}dQ\big|_{T_{x_0}M}.$$

*In a Kähler case when invariancy of potential and Hamiltonian with respect to a flow* $h_t = \exp(t\xi)$ *is supposed, the equation for search of points* $x_0 \in Nul(\alpha_Q)$ *looks like*

(8.7) $$\partial(H - Q - i\varepsilon U) = 0, \quad \varepsilon = \kappa/2.$$

**Theorem 2.** *Let the pair* $(\xi, x_0)$ *generates the Lie solution* $x(t) = \exp(t\xi)x_0$, *and let* $M_0$ *is a manifold of the zero level momentum mapping. Then*

a) *for any* $g \in G$ *the pair* $(\mathrm{Ad}(g)\xi, gx_0)$ *also generates the Lie solution, so the whole orbit* $Gx_0$ *consists of the Lie solutions;*

b) *if* $x_0 \in M_0$, *then dynamics on orbit* $Gx_0 \subseteq M_0$ *is set by a Hamiltonian. If* $x_0$ *is a critical point of a function* $H_P = H - \langle P, \xi \rangle$ *then the role of the corresponding Hamiltonian any* $G$-*invariant function* $H(x)$ *can play;*

c) *all* $G$-*orbits on* $M_0$ *is symplecticly isotropic;*

d) *manifold* $M_0$ *is symplecticly coisotropic.*

As we will see further, the properties of $G$-orbits of Lie solutions are characteristic for attractors, therefore these orbits are of interest in connection with the search of attractors.

Connection of Lie solutions with the reduced system is given by the following theorem.

**Theorem 3.** a) *Any stationary point of reduced CD-system is a projection of the corresponding Lie solution belonging to the set of the zero momentum (i.e. a relative equilibrium), and on the contrary.*

b) *That Lie solution* $x(t) = \exp(t\xi)x_0$ *belongs to the set of the zero momentum, it is sufficient, that the operator* $\mathrm{ad}(\xi)$ *had no nonzero real eigenvalues. For example, it is so for compact Lie groups.*

Let's compare a situation as concern of Lie solutions with the Hamiltonian case.

**Theorem 4.** *Let the potential* $\alpha$ *of the symplectic form* $\Omega = d\alpha$ *and a Hamiltonian* $H$ *are* $G$-*invariant*:

$$g^*\alpha = \alpha, \quad g^*H = H, \quad g \in G.$$

*Then the group G action is Poissonian, with the momentum mapping* $P(x)$:

$$\langle P(x), \xi \rangle = \alpha(\tilde{\xi}(x)), \quad \tilde{\xi}(x) = d/dt\big|_{t=0} \exp(t\xi)x.$$

*The mapping* $t \mapsto x(t) = \exp(t\xi)x_0$ *is a solution to Hamilton system* $\dot{x} = J\,dH$ *iff* $x_0$ *is a critical point of a shifted Hamiltonian*

$$H_P = H - \langle P, \xi \rangle$$

**Example.**

The Lie solutions can not belong to the level set of a momentum. To show this let's consider the Euler system, with the defining form $\alpha = pdq$ generating a phase flow

(8.8) $$g_t(p,q) = \left(e^{-\kappa t} p, q\right).$$

This system admits symmetries

$$A_\varepsilon(p,q) = \left(e^{-\varepsilon} p, e^\varepsilon q\right), \quad B_\varepsilon(p,q) = (p, q+\varepsilon),$$

with generators $A = pq$, $B = p$. The mapping $P:(p,q) \mapsto (A,B)$ coincides with the momentum mapping of Poisson action of the line affinities group

$$(p', q') = \left(a^{-1} p, aq + b\right).$$

Take $x_0 = (p_0, q_0) = (1, 0)$. Then we have the Lie solution $x(t) = g_t x_0 = \left(e^{-\kappa t}, 0\right) = A_{\kappa t} x_0$, and the momentum $P$ at this solution is not a constant.

## §9. The properties of invariant manifolds.

The basic interest for us is represented *by smooth attractors*[22], as the possible model of eigenstates of QM. However, are interesting as well the other, unstable invariant submanifolds. Really, there is also a major set of unstable objects in the quantum world.

It appears, invariant manifolds of CD-systems have the special properties. Here such three properties will be considered: isotropy, Hamiltonity, and quantizedness. The exact sense of these titles will be clear from the further text.

Further $N$ will denote a smooth *compact* invariant submanifold of autonomous CD-system:

$$g_t N = N \subset M, \quad \forall t > 0.$$

### 9.1. Isotropy.

**9.1.1. The autonomous case**. At first sight, it seems almost obvious, that the steady motion in system $(M, \alpha)$ in which there is a constant contraction of the symplectic form, should happen on isotropic submanifold $N \subset M$:

$$d\alpha\big|_{TN} = \Omega\big|_{TN} = 0.$$

---

[22] The exact definition of these objects (see, for example, [3]) play no role for us now. Still we will suppose the attractor to be a compact set.

It is really obvious if to use a *distality* condition which means, that the phase flow on $N$ cannot arbitrarily reduce the length of tangential vectors, in some Riemannian metric on $N$. Formally it means the existence of a constant $c > 0$ for which at all $t > 0$ and $\xi_x \in T_x N$ the inequality

$$\|g_{t*}\xi_x\| \geq c\|\xi_x\|$$

is valid. The proof follows from the fact, that a function

$$f(x;\xi) = |\alpha(x;\xi)| \cdot \|\xi\|^{-1}$$

satisfies the inequality

$$f(g_t x; g_{t*}\xi) \leq c^{-1} \cdot e^{-\kappa t} f(x;\xi),$$

Therefore, for $f_{\max} = \sup f(x;\xi)$ we get $f_{\max} \leq c^{-1} e^{-\kappa t} f_{\max}$, and $f_{\max} = 0$, i.e. $\alpha\big|_{TN} = 0$.

However without this condition to prove the property of isotropy, possibly, is not a simple exercise. The problem is reduced to the lemma proof:

**Lemma 1.** *Let the one-parameter group of transformations $g_t$ acts on a compact manifold, and such 1-form $\beta$ is set that $g_t^* \beta = e^{-\kappa t} \beta$. Then $\beta = 0$.*

Let me represent *the probable reasonings* in favour of this lemma.

We already know that CD-systems $(M^{2n}, \alpha)$ with a compact phase space don't exist. One of approaches to the proof of this fact consists in reviewing of the behaviour of total manifold volume

(9.1) $$v(M) = \int_M (d\alpha)^n,$$

under the influence of a phase flow $g_t$. If $M$ is compact, the inconsistency turns out: the volume should not vary because $g_t M = M$, and on the other hand, the formula (9.1) gives

$$v(g_t M) = e^{-n\kappa t} v(M).$$

In our case it is possible to operate in a similar way. If $\beta \neq 0$, we can consider the invariant integrable distribution $\ker d\beta$ defined in domain $D^{2p}$ where the form $d\beta$ has a constant rank $2p$. We choose the domain $D^{2p}$ of a maximum rank. This domain is invariant, together with its boundary, and the form $(d\beta)^p$ induces the volume form on any manifold $F^{2p}$ of dimension $2p$, which is transversal to distribution $\ker d\beta$, up to domain $D^{2p}$ boundary. The simple reasoning further follows: at small $t > 0$ the manifold $g_t F^{2p}$ is diffeomorphically projected on $F^{2p}$ along the stratums of distribution $\ker d\beta$ owing to what the volume should remain unchanged. On the other hand, this volume should be multiplied by an exponential $e^{-p\kappa t}$, i.e. we obtain an inconsistency.

Owing to the above there are reasons to assume, that the following theorem is valid.

**Theorem 1.** *Any compact invariant submanifolds* $N \subset M$ *of the autonomous system* $(M,\alpha)$ *is symplecticly isotropic, i.e.* $d\alpha|_{TN} = 0$. *Moreover* $\alpha|_{TN} = 0$.

**9.1.2. The nonautonomous case.** Here the systems, nonautonomous on a finite interval of time will be considered.

Let $(\bar{M}, \theta = \alpha_1 - Hdt)$ be such a system considered in expanded space. It will be assumed, that in the remote past $(t < t_0)$ the system moved on an «initial» isotropic smooth attractor $N \subset M$. Let $\bar{N} \subset \bar{M}$ be the trajectory of $N$ in the expanded space:

$$\bar{N} = \bigcup_{\tau \geq 0} \bar{g}_\tau(N, t_0) = \bigcup_{\tau \geq 0} \left( g_{t_0+\tau, t_0} N, t_0 + \tau \right).$$

If to consider $\bar{N}$ as nonautonomous generalisation of an attractor $N$, the analogue of isotropy property of $N$ is described by the following lemma.

**Lemma 2.** *The restriction of the form* $\bar{d}\theta = e^{-\kappa t} d(e^{\kappa t}\theta)$ *on the trajectory* $\bar{N} \subset \bar{M}$ *of isotropic manifold* $N \subset M$ *vanishes:* $\bar{d}\theta|_{T\bar{N}} = 0$.

◀ The tangent vectors $\bar{X}, \bar{Y}$ to manifold $\bar{N}$ in some point $(x,t)$ look like the sums

$$\bar{X} = X + \lambda \bar{V}, \ \bar{Y} = Y + \mu \bar{V}.$$

There $X, Y$ touches of manifold $\bar{N}$ isochronous cut $N(t)$, and $\bar{V}$ is a dynamic field on $\bar{M}$. Considering that $\bar{V}$ belongs to the form $\bar{d}\theta$ kernel we obtain

$$i_{\bar{Y}} i_{\bar{X}} \bar{d}\theta = i_Y i_X \bar{d}\theta.$$

Now the expression of form $\bar{d}\theta$ through the canonical form

(9.2) $$\bar{d}\theta = \Omega - \kappa\alpha \wedge dt$$

gives:

$$i_Y i_X \bar{d}\theta = i_Y i_X (\Omega - \kappa\alpha \wedge dt) = i_Y i_X \Omega = 0.$$

The last equality follows from isotropy of isochronous cut $\bar{N}$ which are passing through a point $(x,t)$. ▶

## 9.2. Hamiltonity.

We will call an invariant submanifold *a Hamiltonian one* if in its neighbourhood there is a Hamiltonian field $J\,dH$ coinciding on this submanifold with a dynamic field $V$. If this invariant submanifold is an attractor, then a Hamiltonian $H$ will be called *a physical one* (PH). Clearly, that this definition does not fix PH uniquely, or to within an additive constant. It would be interesting to discover the natural extra conditions eliminating this arbitrariness.

**9.2.1. The autonomous case.** To prove the Hamiltonity of attractors of autonomous CD-systems, the lemma is required.

**Lemma 3.** *Let the 1-form $\sigma$ vanishes on vectors, tangential to a connected smooth submanifold $N \subset M$: $\sigma|_{TN} = 0$. Then in a tubular neighbourhood of $N$ there is a function $H$ equal to zero on $N$ and with differential which coincides with the form $\sigma$ in any point $x \in N$:*

$$\sigma\big|_{T_N M} = dH\big|_{T_N M}.$$

◄ Let's introduce a Riemannian metric on $M$. We will identify a tubular neighbourhood of manifold $N$ with a neighbourhood of zero section of normal fiber bundle $\mathrm{N}(N)$. Let $x^i, y^j$ be local coordinates on $\mathrm{N}(N)$ ($x^i$ are coordinates in base $N$, $y^j$ are coordinates in fibers, the section $y^j = 0$ is identified with base). In these coordinates the form $\sigma$ looks like

$$\sigma = \sum a_i(x,y)dx^i + b_j(x,y)dy^j,$$

and $a_i(x,0) = 0$. Let $H = \sum b_j(x,y)y^j = i_W \sigma$. Here $W = \sum y^j \partial/\partial y^j$ is a local expression for a dilatation field $W$ of normal bundle. It is important that the function $H$ is invariantly defined in whole tubular neighbourhood of manifold $N$. To end the proof it is necessary to calculate a differential $dH$ at the point $(x, y = 0)$:

$$dH = \sum a_i(x,0)dx^i + b_j(x,0)dy^j = \sigma.$$

▶

It is already easy to prove the property of Hamiltonity.

**Theorem 2.** *Let the vector field $X$ touches of an isotropic submanifold $N \subset M$ of a symplectic phase space $(M, \Omega)$. Then in a tubular neighbourhood of this submanifold there is such a function $H$, that*:

(9.3) $$H\big|_N = 0, \quad X\big|_N = J\,dH\big|_N.$$

◄ Applying the Lemma 3 to the manifold $N$ and form $\sigma = -i_X \Omega$ gives a relation at points of manifold $N$:

$$i_X \Omega = -dH,$$

from which the theorem statement follows. ▶

For CD-systems it gives the following statement.

*There exists a singlevalued PH on a smooth isotropic attractor $N$ of autonomous CD-systems. Every PH take constant values on connected components of $N$.*

### 9.2.2. The nonautonomous case.

**Theorem 3.** *Let $\bar{N} \subset \bar{M}$ be the trajectory of isotropic manifold $N \subset M$. Then the dynamic field $\bar{V}$ on $\bar{N}$ coincides with some Hamilton field[23], defined in a tubular neighbourhood of $\bar{N}$:*

(9.4) $$\bar{V}\big|_{\bar{N}} = \left(\mathbf{J}\,dH + \partial/\partial t\right)\big|_{\bar{N}}.$$

*The equality (9.4) is equivalent to*

(9.5) $$\bar{d}\theta\big|_{T_{\bar{N}}\bar{M}} = \left(\Omega - dH \wedge dt\right)\big|_{T_{\bar{N}}\bar{M}}.$$

◀ We will use Poincare's relative lemma (see [2], p. 30): *the closed differential $k$-form on $\mathcal{M}$, equal to zero on $T\mathcal{N}$, in a tubular neighbourhood of $\mathcal{N} \subset \mathcal{M}$ is represented as a differential of the $(k-1)$-form equal to zero on $T_{\mathcal{N}}\mathcal{M}$.*

According to this lemma and by the previous theorem in a tubular neighbourhood of manifold $\bar{N} \subset \bar{M}$ there is such 1-form $\beta = \beta_1 + \beta_0 dt$, that

(9.6) $$d\left(e^{\kappa t}\theta\right) = d\left(e^{\kappa t}\beta\right),$$

(9.7) $$\beta\big|_{T_{\bar{N}}\bar{M}} = 0.$$

Writing more in detail the right and left members of equality (9.6), we obtain:

$$\begin{cases} \Omega = d_M \beta_1, \\ \kappa\alpha = \dot{\beta}_1 + \kappa\beta_1 - d_M \beta_0. \end{cases}$$

Taking account of (9.7), at points of $\bar{N}$ the equality becomes as

$$\kappa\alpha\big|_{T_{\bar{N}}\bar{M}} = \left(\dot{\beta}_1 - d_M \beta_0\right)\big|_{T_{\bar{N}}\bar{M}}.$$

Now it is necessary to check up the $d_M$-closedness of the form $\dot{\beta}_1 - d_M \beta_0$:

$$d_M \dot{\beta}_1 = \partial/\partial t\left(d_M \beta_1\right) = \partial\Omega/\partial t = 0.$$

---

[23] The Hamiltonian $H(x,t)$ here can be not single-valued.

Here our standard supposition about autonomy of the symplectic form have been essentially used. Hence, a dynamic field $V = \kappa \operatorname{J} \alpha$ on $\bar{N}$ is a Hamiltonian one:

$$V\big|_{\bar{N}} = \operatorname{J}\left(\dot{\beta}_1 - d_M \beta_0\right)\big|_{\bar{N}} = \operatorname{J} d_M H \big|_{\bar{N}},$$

that proves the theorem. ▶

As well as in an autonomous case a Hamiltonian $H(x,t)$ will be called *a physical Hamiltonian*. Further the physical Hamiltonian will be assumed to be *single-valued*.

### 9.3. The asymptotic representation of dynamics.

The Hamiltonity of attractors leads, apparently, to one important result. I will formulate it here as a hypothesis.

**9.3.1. The autonomous case.** Let $H$ be some everywhere defined PH, and let $h_t = \exp(t \operatorname{J} dH)$ be the corresponding Hamiltonian flow.

**Hypothesis.** *For almost all $x \in M$ the mapping is defined*

(9.8) $$x \mapsto r(x) = \lim_{t \to \infty} h_{-t} \circ g_t x \ .$$

*This mapping does not depend on choice of PH. For each smooth attractor $N \subset M$ the mapping $r$ defines a retraction on $N$ of its basin of attraction:*

$$r(r(x)) = r(x) \ .$$

**Corollaries.** 1) *There is the asymptotic representation of dynamics:*

(9.9) $$g_t \approx h_t \circ r, \quad t \to \infty.$$

2) *Relations are valid:*

(9.10) $$g_t \circ r = r \circ g_t = h_t \circ r.$$

3) *There is the «energy quantization» formula:*

(9.11) $$r^* H = \sum E_k I_k \ .$$

*Here $E_k$ is a value of PH on $k$-th connected component of an attractor, and $I_k$ is a characteristic function of this component of the basin of attraction.*

4) *The Poisson brackets of functions having the form $r^* F = F \circ r$ vanishes.*

Let's comment on some of these consequences. The formula (9.9) shows, that if our «observations» concern only to the steady motions the effectively dynamic system looks as Hamiltonian, and all dissipativity is reduced to discrete mapping $r$ (to a projection onto attractors). Obviously, this property calls associations with ideas about «quantum jumps», or with «a re-duction of a wave function».

Last consequence means, in particular, that in the basin of attraction of a $p$-dimensional smooth attractor $N^p \subset M^{2n}$ is defined *a foliation of* $(2n-p)$-dimensional *coisotropic* submanifolds:

$$r^*F_1 = c_1, \ r^*F_2 = c_2, \ ..., \ r^*F_p = c_p,$$

where $r^*F_1, \ r^*F_2, \ ..., \ r^*F_p$ is a maximal collection of functionally independent functions of the given sort. In the case of Lagrangian attractor it turns out a Lagrangian foliation. Such foliations being used in geometrical quantization are called the *polarizations*.

I will add that for *distal* attractors $N$ the mapping $r$ can be defined without using of PH. Really, let in some metric $\rho$ on $M$ the inequality $\rho(g_t x, g_t y) \geq C \cdot \rho(x,y), \ \forall x, y \in N$ holds. The equivalence relation can be introduced then, by taking the definition: $x' \sim y'$ iff

$$\lim_{t \to \infty} \rho(g_t x', g_t y') = 0.$$

Then the point $x'$ is equivalent no more than one point $y \in N$ and it is possible to define $r(x') = y$.

**9.3.2. The nonautonomous case.** The mapping $r$ in nonautonomous case can be defined by the former formula (9.8) if to pass to flows $\overline{g}_s, \overline{h}_s$ in the expanded space:

(9.12) $$\overline{r}(x,t) = \lim_{s \to \infty}(\overline{h}_{-s} \circ \overline{g}_s)(x,t).$$

The important property of mapping $\overline{r}$ is its *isochronism*: it is reduced to a time-dependent mapping $r_t$. It follows from that projections of flows $\overline{g}_s, \overline{h}_s$ on time axis are simple shifts:

$$\overline{h}_s(x,t) = \left(h_{t+s,t}x, \ t+s\right), \quad \overline{g}_s(x,t) = \left(g_{t+s,t}x, \ t+s\right).$$

Therefore the definition (9.12) is reduced to the following:

(9.13) $$r_t(x) = \lim_{s \to \infty} h_{t,t+s} \circ g_{t+s,t} x.$$

Note, that the basin of attraction of «instant attractors» to be defined by these mappings are functions of a time. As well as in the autonomous case, there is an asymptotic factorization:

$$g_{s,t} \approx h_{s,t} \circ r_t, \quad s \to \infty.$$

Apparently, there is some connection of mapping $\overline{r}$ with a dissipative semigroup $\overline{r}_\tau$ which was defined in §6. I will remind that in case of standard gauge $\theta(\overline{V}) = 0$ the dynamic semigroup factorized:

(9.14) $$\bar{g}_\tau = \bar{r}_\tau \circ \bar{h}_\tau = \bar{h}_\tau \circ \bar{r}_\tau.$$

Here $\bar{h}_\tau$ is a Hamilton flow on $\bar{M}$ with *a known* Hamiltonian. Comparing (9.14) with (9.12) we see, that a limit $\bar{f} = \lim_{\tau \to \infty} \bar{r}_\tau$ if it exists, is an isochronous retraction onto some invariant subset in the expanded space. It is natural to assume, that in some cases this subset can be identified with the trajectory $\bar{N} \subset \bar{M}$ of initial smooth attractor $N \subset M$ and mapping $\bar{f}$ can be identified with the mapping $\bar{r}$ considered above. In that case the flow $\bar{h}_\tau$ generator can be cosidered as PH.

### 9.4. Quantizedness.

Let me remind that according to the hypotheses of Chapter 1, we are interested only in CD-systems with non-trivial dynamics to be ensured by «multivalued action function». The quantizedness property concerns to such systems. As a matter of fact, this property *was postulated* in Chapter 1.

To formulate this property we assume the existence of some selected symplectic potential $\omega_1$ which is not coinciding with the defining form $\alpha$. Thanks to that there is a possibility of introduction of nontrivial Hamiltonian and action function. Such potential $\omega_1 = \mathbf{p}d\mathbf{q}$ always exists for the systems derived from Lagrangians. In the autonomous case

(9.15) $$\alpha = \mathbf{p}d\mathbf{q} + \kappa^{-1}dH(\mathbf{p},\mathbf{q}) - dS(\mathbf{q}),$$

and the quantizedness property means that integrals of 1-form $\omega_1$ over the cycles belonging to a smooth attractor $N$ are defined by periods $h_1,...,h_k$ of multivalued function $S(\mathbf{q})$:

(9.16) $$\oint \omega_1 = n_1 h_1 + ... + n_k h_k.$$

The equality (9.16) is an obvious consequence of isotropy of attractor $N$ and expression (9.15) for defining form.

The nonautonomous variant of the same system is described by the form in the expanded space

$$\theta = \mathbf{p}d\mathbf{q} - H(\mathbf{p},\mathbf{q},t)dt - dS(\mathbf{q},t).$$

Let's write out an integral invariant:

$$\oint e^{\kappa t}\theta = \oint e^{\kappa t}\left(\mathbf{p}d\mathbf{q} - Hdt - dS\right) = \oint e^{\kappa t}\left(\mathbf{p}d\mathbf{q} + \kappa^{-1}dH - dS\right) = C = const.$$

Here integration goes on a contour which lies on the attractor trajectory $\bar{N} \subset \bar{M}$.

For isochronous contours this invariant gives

$$\oint \alpha = Ce^{-\kappa t},$$

or

(9.17) $$\oint p\,dq = \oint d_M S + Ce^{-\kappa t}.$$

In view of that in the past the system is supposed to be autonomous, it is necessary to accept $C=0$, that again gives equality (9.16). The new is that now integration goes on the contours lying on arbitrary isochronous cuts of the attractor trajectory $\overline{N}$.

**Examples.**

**1.** The system $\left(M=\mathbb{R}^2,\ \alpha=e^x dy\right)$ has no attractors: $g_t(x,y)=(x-\kappa t, y)$.

**2.** The monopole system considered in an example 3 §7, also has no attractors, however there is the three-dimensional attracting manifolds corresponding to quasi-integral $Q=p_\varphi - h$. At the same time there is an attractor semblance in system. It means the fast transition to dynamic regime (say, to an instant $t_0$) with the approximate equality $p_r \approx h^2/\kappa m r^3$. In this regime the monopole moves on a circle which radius slowly grows:

$$r(t) \approx \left(r^4(t_0) + \frac{4h^2}{\kappa m^2}(t-t_0)\right)^{1/4}, \quad t > t_0.$$

**3.** An arbitrary dynamic system in $\mathbb{R}^n$:

(9.18) $$\dot{x}^i = W^i(x),\ i=1,2,\dots,n$$

extends to CD-system in $M=T^*\mathbb{R}^n$:

(9.19) $$\alpha = \sum_i \left(p_i dx^i + \kappa^{-1} d\left(p_i W^i\right)\right).$$

If the first derivatives of a vector field $W$ are bounded, for big enough a dissipative constant (for example, if $\kappa > \sup_{x,i} \sum_k \left|\partial W^k/\partial x^i\right|$) there is the attracting Lagrangian manifold $p_i=0$, and any attractor of system (9.18) becomes also the attractor of system (9.19).

**4.** To integrable Hamilton system there corresponds a CD-system

(9.20) $$\alpha = (\mathbf{I}-\mathbf{h})d\boldsymbol{\varphi} + \kappa^{-1} dH(\mathbf{I}),$$

where $\mathbf{h}=(h_1,\dots,h_n)$ is the collection of constants selecting an invariant torus $\mathbf{I}=\mathbf{h}$ of Hamilton system. Obviously, it will be a global attractor of system (9.20). The given system is integrable in the sense that a full collection of commuting quasi-integrals $Q_k = I_k - h_k$ generating a global attractor exists. Clearly also, that a Hamiltonian $H(\mathbf{I})$ is a physical one.

Simple computations show, that mapping $r$ here exists for any smooth Hamiltonian $H(\mathbf{I})$, and is given by the formula

$$r:(\mathbf{I},\boldsymbol{\varphi})\mapsto(\mathbf{h},\boldsymbol{\varphi}+\Delta(\mathbf{I})),$$

where the phase shifts are given by the formula

$$\Delta_i(\mathbf{I})=\kappa\sum_j(I_j-h_j)\cdot\int_0^\infty \tau e^{-\kappa\tau}\cdot\left.\frac{\partial^2 H}{\partial I_i \partial I_j}\right|_{\mathbf{I}=\mathbf{I}(\tau)}d\tau,\quad \mathbf{I}(\tau)=\mathbf{h}+e^{-\kappa\tau}(\mathbf{I}-\mathbf{h}).$$

**5.** Let's consider an oscillator

$$\alpha=(I-h(t))d\varphi+\kappa^{-1}\cdot d(I\omega_0),$$

being nonautonomous on a finite time interval, and $h(-\infty)=h_0$, $h(+\infty)=h_1$. As an initial attractor $N$ we take an attractor $I=h_0$ of this system at $t=-\infty$. It is easy to derive from the analysis of solutions, that the trajectory $\bar N$ of this attractor is described by the equation $I=f(t)$, where

$$f(t)=h(t)-e^{-\kappa t}\int_{-\infty}^t e^{\kappa\tau}\dot h(\tau)d\tau.$$

If the motion on $\bar N$ would be described by a single-valued Hamiltonian $H$, on a surface $I=f(t)$ we would have an inconsistent equality:

$$\int_0^{2\pi}\dot I d\varphi=-\int_0^{2\pi}\kappa(I-h(t))d\varphi=-2\pi\kappa(f(t)-h(t))=-\int_0^{2\pi}\frac{\partial H}{\partial\varphi}d\varphi=0.$$

Actually PH here *is multivalued*:

(9.21) $$H(I,\varphi,t)=I\omega_0+\kappa(f(t)-h(t))\cdot\varphi.$$

In the expanded space the same nonautonomous oscillator is described by the form

$$\theta=(I-f(t))\cdot(d\varphi-\omega_0 dt),$$

which satisfies a standard gauge $\theta(V)=0$. In this case the dissipative semigroup $\bar r_\tau$ acts by the formula: $\bar r_\tau(I,\varphi,t)=(I(t,\tau),\varphi,t)$, where $I(t,\tau)=f(t)+e^{-\kappa\tau}(I-f(t))$. The limit at $\tau\to\infty$ gives a retraction on $\bar N$:

$$\bar r(I,\varphi,t)=(f(t),\varphi,t).$$

**6.** At $t\to\infty$ the function $g_t^*F\approx F\circ h_t\circ r$ is constant on stratums of a retraction $r$ and it takes the same values, as on the attractor. For a limit cycle this function represents *a wave* rotating round it.

## §10. The invariant manifolds and the Hamilton-Jacobi equation.

As well as in the Hamilton mechanics in CD-case to a partial solution to the Hamilton-Jacobi equation there corresponds an invariant Lagrangian manifold. This solution allows to specify explicitly also the corresponding PH.

**Theorem 1**. *The submanifold* $\mathbf{p} = \partial S(\mathbf{q},t)/\partial \mathbf{q}$ *of a system* $(\overline{M}, \theta)$ *with the defining form*

(10.1) $$\theta = \mathbf{p}d\mathbf{q} - H(\mathbf{p},\mathbf{q},t)dt$$

*is invariant iff the action* $S$ *satisfies the generalised Hamilton-Jacobi equation:*

(10.2) $$H(\partial S/\partial \mathbf{q}, \mathbf{q}, t) + \partial S/\partial t + \kappa S(\mathbf{q},t) = f(t).$$

*The dynamics on this submanifold is set by a Hamiltonian* $H' = H + \kappa S$.

◄ For invariancy of the submanifold being set by the equations $\varphi_k(\mathbf{p},\mathbf{q},t) = 0$, $k = 1,2,...$ it is necessary and sufficient that relations were satisfied

(10.3) $$\overline{V}.\varphi_k\big|_{\overline{N}} = 0.$$

Here

$$\overline{V} = \frac{\partial}{\partial t} + \frac{\partial H}{\partial \mathbf{p}}\frac{\partial}{\partial \mathbf{q}} - \left(\kappa \mathbf{p} + \frac{\partial H}{\partial \mathbf{q}}\right)\frac{\partial}{\partial \mathbf{p}}$$

is a dynamic field in the expanded space. In this case

$$\varphi_k = \partial S/\partial q^k - p_k,$$

and we obtain

$$\overline{V}.\varphi_k\big|_{\overline{N}} = \left(\frac{\partial^2 S}{\partial q^k \partial t} + \frac{\partial H}{\partial p_i}\frac{\partial^2 S}{\partial q^k \partial q^i} + \kappa p_k + \frac{\partial H}{\partial q^k}\right)\bigg|_{\overline{N}} = \frac{\partial}{\partial q^k}\left(H(\partial S/\partial \mathbf{q}, \mathbf{q}, t) + \partial S/\partial t + \kappa S(\mathbf{q},t)\right) = 0,$$

that gives (10.2). The second statement of the theorem is obvious. ►

Similarly the following statement is derived.

**Theorem 2.** *The submanifold* $\overline{N}$ :

$$\mathbf{p} = \partial S^{(1)}(\mathbf{q},t)/\partial \mathbf{q} = ... = \partial S^{(m)}(\mathbf{q},t)/\partial \mathbf{q}$$

*is invariant iff the functions*

(10.4) $$f^{(i)}(\mathbf{q},t) = H(\partial S^{(i)}/\partial \mathbf{q}, \mathbf{q}, t) + \partial S^{(i)}/\partial t + \kappa S^{(i)}(\mathbf{q},t)$$

*satisfy the equations:*

$$\partial f^{(i)}(\mathbf{q},t)/\partial \mathbf{q}\big|_{\bar{N}} = 0, \quad i = 1, 2, \ldots, m.$$

*The dynamics on $\bar{N}$ is set by any of Hamiltonians $H' = H + \kappa S^{(i)}$.*

**Corollary.** *If all $S^{(i)}(\mathbf{q},t)$ are the various solutions to one and the same generalised Hamilton-Jacobi equation (10.2), then $\bar{N}$ is an invariant manifold.*

Let's give an equivalent statement of the theorem 1, possibly, more interesting in the physical plan.

**Theorem 3**. *The submanifold*

(10.5) $$\mathbf{p} = \partial S(\mathbf{q},t)/\partial \mathbf{q}$$

*of a system $(\bar{M},\theta)$ with the defining form*

(10.6) $$\theta = \mathbf{p}d\mathbf{q} - H(\mathbf{p},\mathbf{q},t)dt - dS(\mathbf{q},t)$$

*is invariant iff the action $S$ satisfies the Hamilton-Jacobi equation:*

(10.7) $$H(\partial S/\partial \mathbf{q}, \mathbf{q}, t) + \partial S/\partial t = f(t).$$

*The dynamics on this submanifold is set by a Hamiltonian $H(\mathbf{p},\mathbf{q},t)$.*

Thus, $H(\mathbf{p},\mathbf{q},t)$ is physical on manifold (10.5), if Hamilton-Jacobi equation holds.

**Remark**. The theorem 3 can be interpreted as follows. Suppouse, we have the Hamilton system described by the form $\theta_0 = \mathbf{p}d\mathbf{q} - H(\mathbf{p},\mathbf{q},t)dt$ and we want to construct the corresponding CD-system for which the Hamiltonian $H$ would be physical. For this purpose we can take the solution $S$ to the equation (10.7). Also we do a gauge transformation

$$\theta = \theta_0 - dS,$$

not changing Hamilton system, but changing the CD-system. Then, if it will appear, that manifold (10.5) contains an attractor, the problem will be solved.

Let's try now to generalise the theorem 3 on arbitrary systems $(\bar{M},\theta)$. Instead of Hamilton field $\bar{W} = \partial/\partial t + \mathrm{J}\,dH$ we will consider the Hamilton field defined by relations

$$i_{\bar{W}}d\theta = 0, \quad dt(\bar{W}) = 1.$$

It is required to find the invariant manifolds on which the dynamics is set by this field. Obviously, any such manifold belongs to the set, on which the field $\bar{W}$ coincides with a dynamic field:

(10.8) $$\bar{W} = \bar{V}.$$

Let's examine such set and when it is invariant. The answers to these problems gives

**Theorem 4**. *The set (10.8) coincides with the set of zeroes* $Nul(\theta \wedge dt)$. *For invariancy of this set it is necessary and sufficient that function* $F(x,t) = \theta(\overline{V})$ *satisfy equalities*

$$(10.9) \qquad \left.\frac{\partial F}{\partial x^i}\right|_{Nul(\theta \wedge dt)} = 0.$$

◀ Let's prove the first statement. We will use the identity

$$(10.10) \qquad \overline{d}\theta - d\theta = \kappa dt \wedge \theta.$$

Let in some point there is a coincidence (10.8). Then taking account of definition of fields $\overline{V}$ and $\overline{W}$, we obtain from (10.10):

$$0 = i_{\overline{V}}(dt \wedge \theta) = \theta - \theta(\overline{V})dt,$$

whence follows, that in this point $\theta \wedge dt = 0$. Inversely, from the same identity follows, that in points of a set $L = Nul(\theta \wedge dt)$ the form $\overline{d}\theta$ and $d\theta$ coincide. The vectors $\overline{V}$ and $\overline{W}$ therefore coincide.

For the proof of the second statement note, that $L$ is described by the equations $\theta_i = 0$, where $\theta_i$ is coefficients of form

$$\theta \wedge dt = \theta_i dx^i \wedge dt.$$

Therefore the condition of invariancy of manifold $L = Nul(\theta \wedge dt) \subset \overline{M}$ is the equality

$$\left.L_{\overline{V}}(\theta \wedge dt)\right|_{T_L \overline{M}} = 0.$$

Using the equality

$$L_{\overline{V}}(\theta \wedge dt) = -\kappa \theta \wedge dt + d\theta(\overline{V}) \wedge dt,$$

we get (10.9). ▶

It is easy to check up, that this theorem contains the previous one. Besides, it gives one more mode of deriving of invariant manifolds when instead of the Hamilton-Jacobi equation the standard gauge condition for the form $\theta$ is used.

**Corollary**. *Consider the system* $(\overline{M}, \theta)$, *satisfying a standard gauge condition* $\theta(\overline{V}) = 0$, *so* $Nul(\theta \wedge dt) = Nul(\theta)$. *Let the set* $Nul(\theta)$ *be a submanifold. Then* $Nul(\theta)$ *is an invariant submanifold, isotropic on isochronous cuts, the dynamics on which is set by a Hamiltonian field* $\overline{W}$.

**Remark**. The form $\theta - \overline{d}\lambda$ also satisfies a standard gauge condition, if $\lambda$ is a quasi-integral. Therefore, if the system admits quasi-integrals, it is possible to put a problem of search of such linear combinations

$$\lambda = c_1\lambda_1 + \ldots + c_m\lambda_m$$

of independed quasi-integrals $\lambda_i$, $\left(L_{\overline{V}}^{\kappa}\lambda_i = 0\right)$, for which the manifold $Nul(\theta - \overline{d}\lambda)$ is not empty. In some cases it can give the whole *spectrum* of invariant manifolds. In particular, within such scheme keeps a search of Lie solutions to autonomous systems.

**Examples.**

**1.** Consider an autonomous analogue of the generalised equation of Hamilton-Jacobi (10.2) for the system from an example §4, in case of one degree of freedom and with a square potential energy:

$$(\partial S / \partial q)^2 / 2m + kq^2 / 2 + \omega \cdot \partial S / \partial \varphi + fq\cos\varphi + \kappa S = const.$$

Let's search for the solutions to this equation which are looking like

$$S(q,\varphi) = aq^2 + \text{Re}\left(bqe^{i\varphi} + ce^{2i\varphi}\right),$$

where $a \in \mathbb{R}$, $b, c \in \mathbb{C}$. In the case $\kappa \geq \kappa_0 = 2(k/m)^{\frac{1}{2}}$ there are only *two solutions* of such sort. To them there correspond the attracting Lagrangian surfaces which are intersected on a limit cycle. These surfaces merge at $\kappa = \kappa_0$, and disappear at $\kappa < \kappa_0$. Unlike these surfaces, the limit cycle exists at all $\kappa > 0$.

**2.** In the expanded space the nonautonomous oscillator is described by the form

$$\theta = (I - f(t)) \cdot (d\varphi - \omega_0 dt),$$

for which $\theta(\overline{V}) = 0$. The set of zeroes $Nul(\theta)$: $I = f(t)$ is invariant. The Hamilton field $\overline{W}$ looks like:

$$\overline{W} = \partial / \partial t + \mathbf{J}\, d_M H, \quad H = \omega_0 I - \dot{f}(t) \cdot \varphi,$$

that agree with an example 5 §9, because $\ddot{f}(t) + \kappa \dot{f}(t) = \kappa h(t)$.

## §11. The statistical approach.

### 11.1. The basic assumtions.

Let's consider a system $\left(\bar{M}, \theta = \alpha_1 - H dt\right)$ in the expanded space, generally nonautonomous, but autonomous in «remote past». It will be convenient to us to introduce special terminology.

**Definitions.** 1) The connected compact invariant submanifold of autonomous CD-system is called *the quanton*[24] of this system.

2) The trajectory in the expanded space of the set of all initial quantons of nonautonomous system is called *the skeleton* of this system.

In spite of the absence of the strict proof, further we will suppose that all quantons are *isotropic*.

Let $\bar{N} \subset \bar{M}$ be a skeleton of the nonautonomous system. We will also assume further, that *the Hamiltonian $H$ is «physical»* for this skeleton. This property is ensured by the gauge transformation preserving a dynamic field $\bar{V}$:

$$\text{(11.1)} \qquad \theta' = \theta + \bar{d}\lambda .$$

Let $X_H = \mathrm{J}\, d_M H$ be a corresponding Hamiltonian field on $M$, and let in expanded space

$$\bar{X}_H = X_H + \partial/\partial t .$$

Thus, it is supposed, that on a skeleton $\bar{N}$ the dynamic field $\bar{V}$ coincides with a Hamiltonian field $\bar{X}_H$:

$$\text{(11.2)} \qquad \bar{V}\big|_{\bar{N}} = \bar{X}_H\big|_{\bar{N}} .$$

Accordingly we have also

$$\text{(11.3)} \qquad \bar{d}\theta\big|_{T_{\bar{N}}\bar{M}} = \left(\Omega - d_M H \wedge dt\right)\big|_{T_{\bar{N}}\bar{M}} .$$

### 11.2. The dynamics of ensemble.

Let's remind, that the statistical ensemble (at instant $t$) in a symplectic phase space $\left(M^{2n}, \Omega\right)$ is set by the generalised probability density $\rho_t(x) \geq 0$, allowing to calculate the averages of finite functions $f \in C_0^\infty(M)$:

$$\langle \rho_t ; f \rangle == \int_M f \cdot \rho_t \, \Omega^n .$$

---
[24] Possibly, it is useful also to include in this concept the isolated *points*.

The dynamics of ensemble is set by a condition of conservation of the probability, which is equivalent to the requirement that for any function $f \in C_0^\infty(M)$ and for any $\tau \geq 0$ the equality was satisfied

(11.4) $$\langle \rho_{t+\tau}; f \rangle = \langle \rho_t; g^*_{t+\tau,t} f \rangle.$$

Here the mappings $g_{t+\tau,t} : M \to M$ are taken from the action of a dynamic semigroup $\overline{g}_\tau$ in the expanded space $\overline{M} = M \times \mathbb{R}$:

$$\overline{g}_\tau(x,t) = (g_{t+\tau,t} x,\, t+\tau).$$

Using (11.4), it is easy to describe the dynamics of ensemble of nonautonomous CD-system with initial data $\rho_{t_0}(x)$ in terms of mappings $\{g_{t_2,t_1}\}$:

(11.5) $$\rho_{t_0+\tau}(x) = \begin{cases} e^{n\kappa\tau} \rho_{t_0}(g^{-1}_{t_0+\tau,t_0} x), & \text{if } x \in g_{t_0+\tau,t_0}(M), \\ 0, & \text{if } x \notin g_{t_0+\tau,t_0}(M). \end{cases}$$

In turn, (11.5) is equivalent to invariance of the generalized $2n$-form $\rho\Omega^n \sim \rho_t(x)\Omega^n$ in the expanded space[25]:

(11.6) $$\overline{g}^*_\tau(\rho\Omega^n) = \rho\Omega^n,$$

or to the Liouville equation

(11.7) $$L_{\overline{V}}(\rho\Omega^n) = 0.$$

Usually the Liouville equation (11.7) is written in the form of a continuity equation for a density $\rho_t$:

(11.8) $$\frac{\partial \rho_t}{\partial t} + \mathrm{div}(\rho_t V) = 0.$$

For CD-systems the last equation is:

(11.9) $$\partial \rho_t / \partial t + V.\rho_t - n\kappa\rho_t = 0.$$

The solution to Cauchy problem for the equation (11.9) is given by the formula (11.5).

Further us the description of dynamics of ensemble in terms of Hilbert space will interest. It will be assumed besides, that the ensemble *is localised on a skeleton*. Therefore we cannot in-troduce vectors $\psi_t$ of Hilbert space starting from an equality $\rho_t = |\psi_t|^2$. At first the selection of singular «background» $\rho$ of a density is necessary. After that we can consider the densities of the form $|\psi|^2 \cdot \rho$.

---

[25] If the time-dependent density is considered in $M$ we will denote it by $\rho_t$, if it is considered as object in the expanded space $\overline{M}$, we will write simply $\rho$. The similar remark concerns to notations of functions.

### 11.3. An invariant measure on a skeleton.

The identification a measure on a skeleton with a density $\rho$ here is meant. It can be constructed as a solution $\rho(x,t) = \rho_t(x)$ to Cauchy problem for the Liouville equation (11.9). As initial data $\rho_{t_0}$ we can take, for example, a time average of a measure on the initial cut $N$ of the skeleton, obtained from any Riemannian metric. Such distribution will be *invariant* distribution on $M$ for initial autonomous system[26].

In terms of the expanded space, the constructed solution $\rho(x,t)$ to Cauchy problem can be described as a generalized function of measure type on $\bar{M}$ which at $t \geq t_0$ is concentrated on a skeleton and for which the generalized $2n$-form $\rho\Omega^n$ is invariant:

$$\bar{g}_\tau^*(\rho\Omega^n) = \rho\Omega^n.$$

The following theorem shows, that there is a possibility to describe a dynamics of ensemble on a skeleton, starting only from an initial density $\rho_{t_0}$ with the Hamilton system $(M, \Omega, H)$. No other information about the initial CD-system for this purpose is required.

**Theorem 1.** *Let $\rho_t$ be a solution to Cauchy problem for the equation* (11.9) *with initial data $\rho_{t_0}$ of measure type on the initial cut of a skeleton. Then $\rho_t$ is also a solution to Cauchy problem, with the same initial data, for the Liouville equation of Hamilton system with a Hamiltonian $H$:*

(11.10) $$\partial \rho_t / \partial t + X_H . \rho_t = 0.$$

◄ Let $\bar{h}_\tau = \exp(\tau \bar{X}_H) : (x,t) \mapsto (\bar{h}_{t+\tau,t} x, \, t + \tau)$ be the one-parametrical group in expanded space generated by the Hamilton field $\bar{X}_H = \partial/\partial t + X_H$ of physical Hamiltonian $H$. Then, for any $f \in C_0^\infty(M)$ the function $\varphi(x) = f(g_{t+\tau,t} x) - f(h_{t+\tau,t} x)$ vanishes, if $x$ belongs to the distribution $\rho_t$ support, i.e. to the isochronous cut of skeleton $\bar{N}$. Since $\rho_t$ is of measure type, we get a relation

$$\langle \rho_t ; \varphi \rangle = 0,$$

---

[26] In Chapter 1 we consider q-states which turn out a time average of usual states and which, apparently, are physically more sensible, than arbitrary abstract ensembles. The invariance of initial distribution $\rho_{t_0}$ means, that this distribution can be expressed in terms of q-states.

or

(11.11) $$\langle \rho_t; g^*_{t+\tau,t} f \rangle = \langle \rho_t; h^*_{t+\tau,t} f \rangle.$$

Let now $\tilde{\rho}_t$ be a solution to Cauchy problem for the Liouville equation (11.10), with the same initial data $\rho_{t_0}$. Then, using (11.11), and considering the general principle expressed by the formula (11.4), we obtain:

$$\langle \rho_{t_0+\tau}; f \rangle = \langle \rho_{t_0}; g^*_{t_0+\tau,t_0} f \rangle = \langle \rho_{t_0}; h^*_{t_0+\tau,t_0} f \rangle = \langle \tilde{\rho}_{t_0+\tau}; f \rangle,$$

whence

$$\rho_{t_0+\tau} = \tilde{\rho}_{t_0+\tau},$$

that proves the theorem. ▶

### 11.4. Construction of Hilbert space.

Let's consider connection in a product fiber bundle $\bar{M} \times \mathbb{C} \to \bar{M}$, set by a covariant derivative:

(11.12) $$\nabla_{\bar{X}} = \bar{X} - \tfrac{i}{\hbar} \cdot \theta_1(\bar{X}).$$

We need that curvature of this connection would be proportional to the 2-form

$$\Omega - d_M H \wedge dt.$$

Let's assume for this purpose, that the symplectic form admits a potential which is not dependent on time:

$$\Omega = d_M \omega_1, \quad \partial \omega_1 / \partial t = 0.$$

Then we can take the $\theta_1$ which looks like

$$\theta_1 = \omega_1 - H dt - dF.$$

However, in view of possibility of the gauge transformation of connection, it is sufficient to take the case $F = 0$.

**Theorem 2.** *Let $\rho_t$ be a solution to the Cauchy problem for the Liouville equation* (11.10) *with the initial data concentrated on the initial cut of a skeleton $N = N_{t_0}$, and let $\psi(x,t) = \psi_t(x)$ be a solution to the Cauchy problem for the equation*

(11.13) $$\nabla_{\bar{X}_H} \psi = 0,$$

*with finite initial data*

(11.14) $$\psi_{t_0} = \psi\big|_{t=t_0} \in C_0^\infty(M).$$

*Then the function $|\psi_t|^2 \rho_t$ is a solution to the Cauchy problem for the equation* (11.10), *and an integral*

(11.15) $$\|\psi\|^2 = \int_M |\psi_t|^2 \cdot \rho_t \Omega^n$$

*does not depend on time.*

◄ From (11.13) follows, that function $|\psi_t|^2$ satisfies the Liouville equation

$$\bar{X}_H \cdot |\psi|^2 = 0,$$

or

(11.16) $$\partial |\psi_t|^2 / \partial t + X_H \cdot |\psi_t|^2 = 0,$$

and therefore as well a density $|\psi_t|^2 \rho_t$ satisfies the Liouville equation. Conservation of a total probability for this density leads to invariance of an integral (11.15). ►

**Remark.** The pair $(\rho_t, \psi_t)$ where $\rho_t$ is the background density localised on isochronous cuts of a skeleton, and $\psi_t$ is the complex function defined in a neighbourhood of the density $\rho_t$ support, would be possible to call *a coherent ensemble.* Within the limits of such terminology the equations (11.10), (11.13) are the equations of motion of coherent ensemble. To coherent ensemble there corresponds usual ensemble with a density $\rho'_t = |\psi_t|^2 \cdot \rho_t$, that should ensure the «square module» rule of quantum mechanics if to consider $\psi_t$ as an analogue of a wave function.

The integral (11.15) is possible to use for construction of the *time-dependent* Hilbert space $H_t$ of «wave functions» localised on the isochronous skeleton cut $N_t$. Obviously, the equation (11.13) defines a set of invertible isometric operators:

$$U_{t_2, t_1} : H_{t_1} \to H_{t_2},$$

that allows to identify these spaces.

In the autonomous case the system and a background density don't depend on time, $\rho_t = \rho_{t_0}$, so the density $\rho_{t_0}$ is invariant on $M$. This case is interesting because all $H_t$ are identified naturally, undependently on the dynamics. Then the operators

(11.17) $$U_{t_2 - t_1} = U_{t_2, t_1}$$

become usual unitary operators in the common Hilbert space $H$.

## 11.5. Comparison with geometrical quantization.

The basic problem of the theory of geometrical quantization (GQ), as is known, consists in *twice* to reduce the number of variables on which the wave function $\psi$ of prequantum Schrödinger equation depends[27], and to obtain thereby the true Schrödinger equation. The co-variant constancy of $\psi$ along the Lagrangian foliation (real *polarisation*) for this purpose is postulated. This condition selects separate stratums on which the wave function is localized. As a whole such stratums form the Bohr-Sommerfeld submanifold $BS \subset M$. Further two cases are possible:

- the Hamilton field $X_H$ does not touches $BS$;

- the Hamilton field $X_H$ touches $BS$.

The basic difficulty is represented by the first case. The second case is much easier, as then the prequantum Schrödinger equations is automatically reduced to the equation in space of functions depending on twice smaller number of arguments, that solves a quantization problem. Fortunately, we are interested in the second case.

The analogy of our constructions with the constructions of GQ consists in the following. We will begin that, as well as in GQ, we have a prequantum Schrödinger equation.

**Theorem 3.** *The equation* (11.13) *can be written in the form of a Schrödinger equation*

(11.18) $$i\hbar \frac{\partial \psi}{\partial t} = \hat{H}\psi,$$

*with a prequantum Hamilton operator*:

$$\hat{H} = H - \omega_1(X_H) - i\hbar \cdot X_H.$$

◀ Considering, that $\theta_1(\bar{X}_H) = \omega_1(X_H) - H$, we obtain

(11.19) $$\nabla_{\bar{X}_H} = \frac{\partial}{\partial t} + X_H - \frac{i}{\hbar} \cdot (\omega_1(X_H) - H) = \frac{\partial}{\partial t} + \frac{i}{\hbar} \cdot \hat{H},$$

where $\hat{H}$ is *the prequantization operator* (see, for example, [11]). ▶

Let's add to the theorem the statement, that in an autonomous case the operator $\hat{H}$ is correctly defined in a Hilbert space $H$ and generates there the unitary group (11.17):

(11.20) $$U_t = \exp(-it\hat{H}/\hbar).$$

---

[27] Roughly speaking, it is required to pass from the equation for $\psi(p,q,t)$ to the equation for $\psi(q,t)$.

The following problem concerns the analogue *of a polarisation*. In this relation we can point out to the foliation of *coisoptropic* submanifolds, existing in basin of attraction of every isotropic attractor (section 9.3.1 see). Further, the submanifold $BS$ can be compared with the union $N \subset M$ of all quantons. The isotropy of $N$, obviously, generalises the property of Lagrangianity of $BS$-manifolds. Besides, the manifold $N$, as well as $BS$, also obeys to some quantization conditions (section 9.4 see).

Further, similar to how in GQ the wave functions are concentrated on $BS$, our wave functions are localized on $N$. We will specify now the analogue of the property of wave function be a covariant constant along a polarisation.

**Theorem 4.** *In a tubular neighbourhood of the trajectory $\bar{N}_\alpha \subset \bar{M}$ of every quanton $N_\alpha$ there is a bundle $\bar{M} \times \mathbb{C} \to \bar{M}$ section $\psi_\alpha$, covariant constant over $\bar{N}_\alpha$, and having the unit module: $|\psi_\alpha| = 1$. The restriction $\psi_\alpha|_{\bar{N}_\alpha}$ is defined to within a constant phase factor.*

◄ The connection (11.12) is locally flat over a skeleton $\bar{N} \subset \bar{M}$:

$$d\theta_1|_{T\bar{N}} = 0.$$

Let's use the Poincare relative lemma. According to this lemma, there is such 1-form $\tilde{\theta}_1$ in a tubular neighbourhood of a skeleton $\bar{N}$ that

$$d\theta_1 = d\tilde{\theta}_1, \quad \tilde{\theta}_1\big|_{T_{\bar{N}}\bar{M}} = 0.$$

Let $\sigma = \theta_1 - \tilde{\theta}_1$. Then the form $\sigma$ is closed, and

(11.21) $$\theta_1\big|_{T_{\bar{N}}\bar{M}} = \sigma\big|_{T_{\bar{N}}\bar{M}}.$$

According to a quantization hypothesis (section 9.4 see.), taken in its elementary aspect, the periods of a form $\sigma$ on a skeleton (and in its tubular neighbourhood) are whole multiple of the Planck constant. Hence, in this neighbourhood there is a multivalued potential $S$ of this form:

$$\sigma = dS,$$

with periods multiple $h = 2\pi\hbar$. Also the section $\psi = e^{iS/\hbar}$ is correctly defined.

Let now $\bar{X}$ be the arbitrary vector field, tangential to a skeleton. Then we have

$$\nabla_{\bar{X}} \psi = \tfrac{i}{\hbar}\big(dS(\bar{X}) - \theta_1(\bar{X})\big) \cdot \psi.$$

Owing to equality (11.21) this expression vanishes over a skeleton. It also means, that $\psi = e^{iS/\hbar}$ is a section which is a covariant constant over a skeleton. We will note now, that a potential $S$ over a skeleton is defined up to the independent shift $S$ on the constant $c_\alpha$ over each connected component $\bar{N}_\alpha$. From here follows, that each a function

$$\psi^{(\alpha)} = \psi\big|_{\bar{N}_\alpha}$$

is defined to withing a phase factor. ▶

Note, that this theorem contains such statements:

- functions $\psi^{(\alpha)}$ satisfy the equation (11.18);

- functions $\psi^{(\alpha)}$ are covariant constant over quantons $N_\alpha$.

**Corollary**. *In the autonomous case each $\psi^{(\alpha)}$ is an eigenfunction for a prequantum Hamilton operator $\hat{H}$:*

(11.22) $$\hat{H}\psi^{(\alpha)} = E_\alpha \psi^{(\alpha)},$$

and $E_\alpha = H\big|_{\bar{N}_\alpha} = const$.

◀ In the autonomous case the field $\bar{X} = \partial/\partial t$ touches to $\bar{N}_\alpha$. Therefore, according to theorem 4,

$$0 = \nabla_{\bar{X}} \psi^{(\alpha)} = \bar{X}.\psi^{(\alpha)} - \frac{i}{\hbar}\theta_1(\bar{X})\psi^{(\alpha)} = \partial \psi^{(\alpha)}/\partial t + \frac{i}{\hbar} H \cdot \psi^{(\alpha)}.$$

Comparing this equality with (11.18), we obtain

$$\hat{H}\psi^{(\alpha)} = H \cdot \psi^{(\alpha)}.$$

To end the proof, we must consider the property of physical hamiltonian be a constant over the connected components of smooth isotropic attractors. ▶

Further, the covariant constant functions $\psi^{(\alpha)}$ have a finite norm owing to the compactness of quantons. Hence, there are reasons to assume that exactly these functions generate the Hilbert «quantization space» $H_Q$, containing in $H$.

So, we obtain a genuine Schrödinger equation, if:

a) the system is autonomous, and the background density $\rho_t$ is stationary: $\rho_t = \rho_{t_0}$;

b) the wave functions having a form $\psi_t = \sum_\alpha c_\alpha \psi_t^{(\alpha)} \in H_Q$ are considered.

Under these conditions $H_Q$ will be an invariant subspace in $H$, and restriction of prequantum dynamics on $H_Q$ gives *a stationary* Schrödinger equation. Besides, owing to that squares of modules of functions $\psi_t \in H_Q$ are constant over quantons, the ensembles $\rho'_t = |\psi_t|^2 \cdot \rho_t$ will be the *invariant* ensembles of the autonomous CD-system. Last property allows to formulate the following probable principle:

*The space $H_Q$ and background density $\rho_t$ of the autonomous system should be defined so that all invariant ensembles have the form*

$$(11.23) \qquad \rho'_t = |\psi_t|^2 \cdot \rho_t, \quad \psi_t \in H_Q.$$

However for the above definition $H_Q$ this principle, generally speaking, *is not valid*. It is possible, for example, a degenerate case when on some quanton $N_\alpha$ there is a *set* of invariant densities having the form $\lambda(x) \cdot \rho^{(\alpha)}$, where $\lambda(x) \geq 0$ are invariant functions, and $\rho^{(\alpha)}$ is some invariant background density[28].

In similar cases the modification of definition of space $H_Q$ is necessary. For example, in degenerate case it is possible to accept, that $H_Q$ is a direct sum of Hilbert spaces:

$$H_Q = \bigoplus_\alpha H_Q^{(\alpha)},$$

where every $H_Q^{(\alpha)}$ is generated by functions

$$\psi_\lambda^{(\alpha)} = \sqrt{\lambda(x)} \cdot \psi^{(\alpha)},$$

having a norm being defined by invariant density $\lambda \rho^{(\alpha)}$:

$$\left\| \psi_\lambda^{(\alpha)} \right\|^2 = \int_M \left| \psi^{(\alpha)} \right|^2 \lambda \rho^{(\alpha)} \Omega^n = \int_M \lambda \rho^{(\alpha)} \Omega^n.$$

For such definition the squares of modules of functions $\psi_t \in H_Q$ will be invariant functions, but these can be not constant over quantons.

But we will not go deep and consider here all conceivable possibilities. Further we will consider *non-degenerate* autonomous systems.

## 11.6. Connection with a quantum formalism.

The basic problem which us now interests is that the space $H$ and unitary group of Schrödinger equation as we saw, make sense only *in the autonomous case*. How, in that case, to explain the seeming (or genuine) *universality* of a quantum formalism?

For the purpose of situation clearing let's take an example of asymptotically autonomous system for which an external perturbation is introduced only on a finite time interval $(T_1, T_2)$.

---

[28] The case of existing of unique invariant normalized density on every quanton can be called *non-degenerate*.

Let's assume, that the initial condition of the system being set in the remote past, is described by the coherent ensemble $(\rho_{t_0}, \psi_{t_0})$, where function $\psi_{t_0} \in H$, and $\rho_{t_0}$ is an invariant background density. Then for the states $(\rho_t, \psi_t)$ of remote future the density $\rho_t$ will not coincide with initial stationary density:

$$\rho_t \neq \rho_{t_0},$$

and, moreover, will be non-stationary, in general. Hence, the dynamics after the «scattering» *goes beyond of Hilbert space H*.

This problem, however, can be solved easily. As far from a scattering interval $(T_1, T_2)$ we deal with the same autonomous system, the support of a final asymptotic density $\rho_t$ will be close to the support of an initial background density $\rho_{t_0}$ and hence there is an equality

$$\rho_t \approx \nu_t \cdot \rho_{t_0},$$

where $\nu_t$ is some function, which satisfy the Liouville equation

$$\partial \nu_t / \partial t + X_H . \nu_t = 0.$$

Note also, that this function does not depend in any way on $\psi_t \in H$. Hence, from the pair $(\rho_t, \psi_t)$ we can pass to the statistically equivalent pair[29]

(11.24) $$(\rho_t, \psi_t) \sim \left(\rho_{t_0}, \psi'_t = \sqrt{\nu_t} \cdot \psi_t\right),$$

where already $\psi'_t \in H$. And, as it is easy to see, the function $\psi'_t$ satisfies a prequantum Schrödinger equation, so the pair $(\rho_{t_0}, \psi'_t)$ satisfies the dynamic equations of coherent ensemble of autonomous system.

Now it is possible to compare the initial and final dynamics, being within the limits of *the same space*[30] *H*. Owing to the total probability conservation the functions $\psi_{t_1}$ and $\psi'_{t_2}$ at $t_1 \to -\infty$, $t_2 \to +\infty$ are connected by *a unitary operator*. Besides, the asymptotical behaviour of these functions assumes the existence of such functions $\psi_{in}, \psi_{out} \in H$ that

$$\psi_{t_1} \approx U_{t_1}\psi_{in}, \qquad \psi'_{t_2} \approx U_{t_2}\psi_{out},$$

where $U_t = \exp(-it\hat{H}/\hbar)$ is a unitary group of Schrödinger equation in $H$. Hence, in general situation it is natural to expect the existence *of a unitary scattering operator* $\hat{S}$ defined by the equality

$$\psi_{out} = \hat{S}\psi_{in}.$$

---

[29] Pairs $(\rho_t, \psi_t)$ and $(\rho'_t, \psi'_t)$ are supposed to be statistically equivalent, if $|\psi_t|^2 \cdot \rho_t = |\psi'_t|^2 \cdot \rho'_t$.

[30] Similarly, we can consider the case when asymptotic autonomous CDs at $t \to \pm\infty$ are different, and $H_{in} \neq H_{out}$. In that case we also can identify $H_{in}$ and $H_{out}$. Isomorphism of Hilbert spaces arises by itself, it has dynamic character.

On the other hand, the existence of scattering operator, together with a «square module» rule of probabilistic interpretation, are those beginnings which, apparently, are enough for a deduction of all quantum formalism. In this sense here the basic connection of CD-systems with the quantum theory has been established.

There is, however, one difficulty. In general case the «quantization» subspace $H_Q \subset H$ will not be an invariant subspace for the operator $\hat{S}$.

The possible way of solution this problem consists in reduction of initial autonomous system to turn the quantons into points. In this case $H = H_Q$ and to go out of $H_Q$ limits is impossible. In the models of Chapter 3 such reduction is really realizable.

**Example.**
For the system $\theta = (I - \hbar)d\varphi - H(I)dt$ the skeleton is set by the equation $I = \hbar$ and we can take $\theta_1 = Id\varphi - H(I)dt$. The operator of prequantization looks like

$$\hat{H} = H(I) - H'(I) \cdot (I + i\hbar \cdot \partial/\partial\varphi),$$

and

$$\psi = e^{iS/\hbar} = e^{i\varphi} \cdot e^{-itE/\hbar},$$

where $E = H(\hbar)$. Obviously, $\hat{H}\psi = E\psi$. The space $H_Q$ here is one-dimensional, and the space $H$ is infinite-dimensional.

## §12. The local field systems.

If we wish to build a dynamically nontrivial field CD-system it is necessary to begin, most likely, from a multivalued phase functional[31]

$$\Phi(f) = S(f)/\hbar.$$

---

[31] In order to avoid a confusion, instead of an «action functional» we will speak about a «phase functional».

More strictly, instead of this functional we can say about the closed not exact 1-form $\delta\Phi$ on a functional space of fields $f$. And, for the equations of motion to have the form of partial differential equations, the variation $\delta\Phi$ must have the form of an integral of a local field density. The remaining elements (the phase space, symplectic structure, Hamiltonian) are quite traditional. It is possible anyway to borrow them or to build just like their analogues in the known field Hamilton systems.

There is a description of fields for which such multivalued functionals are possible (see [8]). Such fields $f$ are interpreted geometrically as mappings

(12.1) $$f : N^q \to E$$

of fixed manifolds. In particular, the *bundle sections* can be such fields. Let the closed form $\omega_{q+1}$ represents a nonzero cohomology class $[\omega_{q+1}] \in H^{q+1}(E;\mathbb{R})$. Let's set in correspondence to a field $f$ and its small increment $\delta f$ a narrow strip $S_{\delta f} \subset E$ which in general case should be $(q+1)$-dimensional:

$$S_{\delta f} : (x, \tau) \mapsto f(x) + \tau \cdot \delta f(x), \quad x \in N^q, \ 0 \leq \tau \leq 1.$$

Then the form $\delta\Phi(f, \delta f)$ is defined as a principal linear part of a flux of the form $\omega_{q+1}$ through this strip:

(12.2) $$\int_{S_{\delta f}} \omega_{q+1} = \delta\Phi(f, \delta f) + O\left(\|\delta f\|^2\right).$$

If the form $\omega_{q+1}$ is written in local coordinates $(y^i)$ in the form of sum (over all indices)

$$\omega_{q+1} = \sum \omega_{i_0 \ldots i_q}(y) dy^{i_0} \wedge \ldots \wedge dy^{i_q},$$

with antisymmetric coefficients $\omega_{i_0 \ldots i_q}$ then the functional $\Phi$ is defined by the formula:

(12.3) $$\delta\Phi(f; \delta f) = (q+1) \cdot \sum \int_{N^q} \omega_{i_0 \ldots i_q}(f(x)) \cdot \delta f^{i_0}(x) \cdot df^{i_1}(x) \wedge \ldots \wedge df^{i_q}(x).$$

The main hypothesis of the given work is that the united field theory should be a local field CD-system with a local multivalued functional.

Obviously, any ideas limiting arbitrariness in a choice of a field $f$ are important. Now we know two general restrictions of this kind. The first is the requirement to the periods of the functional $\Phi$. The second is the condition *of asymptotical additivity* (which was already mentioned in §4 of Chapter 1) of this functional. For the field theory it means the additivity over the «physical space». It can be valid, if the physical space is identified with a base $N^q$, or with a base quotient set. For example, we can take the manifold $N^q$ from the theory of Kaluza-Klein type.

It should be noted also that multivalued functionals naturally arise for gauge fields. For example, in space of modules of non-Abelian field in $E^3$ the Chern-Simons functional is known, with a variation

$$\delta\Phi = -\tfrac{1}{4\pi}\int \varepsilon^{ijk}\,\mathrm{Tr}\!\left(F_{ij}\cdot\delta A_k\right)d^3\mathbf{x}.$$

## §13. A remark about the singularities.

There is one more important condition for to limit arbitrariness in a choice of CD-systems. It is necessary to ensure *a regularity* of the constructed dynamic system. It means, that the system should not be «blowing up», passing to infinity for finite time. Besides, the singularities *in a phase space* should be *repelling*.

The existence of the repellent singularities, however, can seem a serious drawback. Really, this fact means the presence of singularities in the past , because time reversion converts the repellent singularities into the attracting ones.

How to concern to such singularities? Note, that mathematically the system is quite correct, and problems arise at attempt to invert the dynamics of irreversible system, in general.

From the physical point of view any concrete system actually is a part of wider system, and the singularity can be a price for the simplified description. For example, in the remote past our system generally could not exist as object, and have been separated from the big system at the moment of its birth. Then, obviously, not every extrapolation of the states of our system in the past will have a physical sense. Also the recipe of liquidation of singularities from here implies: we must expand the system.

Sometimes, possibly, will be enough simple prolongation (with preservation of dimension) a system phase space. The example gives the oscillator model. We can describe the states either by an action-angle variables $I>0, \varphi$, or by a holomorphic coordinate $z=\sqrt{I}e^{i\varphi}$. In the second case we have the Kahler system with a repellent singularity in a point $z=0$. To liquidate this singularity we can return to the action-angle variables $I, \varphi$ and then admit arbitrary values of the action.

## §14. Dynamic symmetries.

The definition of CD-system symmetry accepted in §7 is insufficient because it don't include the usual constants of motion. However it is natural also to suppose their connection with symmetry.

**Definition.** *A dynamic symmetry* is a Lie algebra $Lie(\tilde{G})$ of vector fields in a phase space, containing a dynamic field of CD-system $(M,\alpha)$. It means, that there are linearly independent fields $X_0 = V, X_1, \ldots, X_m$ for which

(14.1) $$[X_i, X_j] = \sum_k C_{ij}^k X_k.$$

We don't speak here about the symmetry *group* because the dynamic semigroup $g_t = \exp(tV)$ is not obliged to be a group and consequently the Lie algebra $Lie(\tilde{G})$ generally cannot be integrated to the genuine action of the Lie group $\tilde{G}$ in $M$.

Let's select a part of relations (14.1) corresponding to the values of indexes $i = 0$, $j > 0$:

(14.2) $$[V, X_j] = \sum_{k>0} A_j^k X_k + B_j V.$$

Applying the form $\alpha$ to (14.2) and using identities of §1, we get for the functions $F_j = \alpha(X_j)$ a linear equations

(14.3) $$L_V^\kappa F_j = \sum_{k>0} A_j^k F_k.$$

From here the generalisation of theorem about the quasi-integrals follows.

**Theorem 1.** *In the case of dynamic symmetry there is a relation*

(14.4) $$F(g_t x) = e^{(\mathbf{A} - \kappa)t} F(x),$$

where $F$ is a vector-function having the components $F_j = \alpha(X_j)$, and $\mathbf{A}$ is the matrix $(A_j^k)$.

Obviously, for usual symmetry a matrix $\mathbf{A}$ is zero, and (14.4) means, that all functions $F_j$ are quasi-integrals.

Let's assume now, that all fields $X_i$, except a field $V = X_0$, are the *Hamiltonian ones*:

(14.5) $$X_i = \mathbf{J} dH_i, \quad i > 0.$$

Then this fields generate a Lie subalgebra $Lie(G) \subset Lie(\tilde{G})$. It follows from (14.1), if to consider, that a commutator of the Hamilton fields is the Hamilton field. The commutator $[V, X_j]$ also is the Hamilton field, therefore[32] $B_j = 0$. Besides, using the identities of §1, we derive the equalities

---

[32] Note, that the equality $B_j = 0$ excludes a symmetry of relativistic type: the group $\tilde{G}$ cannot be the Poincare group, de Sitter group, etc.

(14.6) $$L_V^\kappa H_j = \sum A_j^k H_k + b_j, \quad j, k > 0,$$

(14.7) $$\{H_i, H_j\} = \sum C_{ij}^k H_k + b_{ij}, \quad i, j, k > 0,$$

with some constants $b_i, b_{ij}$. If it is impossible to liquidate these constants by overdetermination of Hamiltonians $H_i \to H_i + c_i$, we will introduce the Hamiltonian $H_0 \equiv 1$ linearly independent with the others Hamiltonians $H_i$, $i > 0$. The expanded collection of Hamiltonians sets the momentum mapping

(14.8) $$x \mapsto \{H_0(x), H_1(x), \ldots, H_m(x)\} = P(x) \in Lie^*(G_1)$$

corresponding to Poisson action of some new Lie group $G_1$. This group is a one-dimensional central expansion of group $G$, in the sense that $G \approx G_1 / Z_1$, where $Z_1$ is a one-dimensional subgroup from the center.

As a result of such expansion of a group the constants $b_{ij}$ in equality (14.7) disappear, but in equality (14.6) constants $b_j$ remain, in general. As a result, without loss of generality, we can accept the following

**Definition**. By *Hamiltonian dynamic symmetry* of CD-system we will call a connected Lie group $G$, if the action of this group is Poissonian, and if the momentum mapping $P: M \to Lie^*(G)$ corresponding to this action is transformed linearly by a phase flow:

(14.9) $$L_V^\kappa P = \mathbf{A} P + b,$$

where $b \in Lie^*(G)$ is a vector constant.

**Remark**. Here is meant, that the operator $\mathbf{A}$ and vector $b$ are *uniquely* defined by equality (14.9). We will assume for this purpose the following conditions. If among the Hamiltonians $H_\xi = \langle P, \xi \rangle$ there are nonzero constants, should be $b = 0$. Besides, the Lie algebra $Lie(G)$ homomorphism $P: \xi \mapsto H_\xi = \langle P, \xi \rangle$ in a Lie algebra of functions should have *the zero* kernel: $\ker P = 0$. If it not so the given drawback is eliminated by replacement of $Lie(G)$ by the factor-algebra $Lie(G') = Lie(G) / \ker P$.

**Theorem 2.** *Let the action of connected Lie group $G$ in a phase space of CD-system with the defining form $\alpha$ is Poissonian. Then the following three conditions are equivalent:*

a) $G$ *is a Hamilton dynamic symmetry, i.e. there is an equality*

(14.10) $$L_V^\kappa P = \mathbf{A} P + b.$$

b) *The relations are satisfied*

(14.11) $$L_\xi \alpha = -\kappa^{-1} d\langle \mathbf{A}P, \xi \rangle, \quad \xi \in Lie(G).$$

c) *The law of transformation of the defining form looks like*

(14.12) $$g^*\alpha = \alpha + \kappa^{-1} d\langle P, \psi(g) \rangle, \quad g \in G,$$

where $P$ is a momentum, and $\psi(g) \in Lie(G)$ is some function.

*Besides, in the case of dynamic symmetry, the following statements are valid.*

d) *There is an equality*

(14.13) $$P(g_t x) = e^{(\mathbf{A}-\kappa)t} P(x) + \frac{e^{(\mathbf{A}-\kappa)t} - 1}{(\mathbf{A}-\kappa)} b.$$

e) *The adjoint operator* $\mathbf{A}^*$ *is the algebra* $Lie(G)$ *differentiation. It is defined by the function* $\psi(g)$

(14.14) $$\xi \mapsto \mathbf{A}^*\xi = -\left.\frac{d}{dt}\right|_{t=0} \psi(e^{t\xi}).$$

f) *The constant* $b$ *satisfies a condition*

(14.15) $$\langle b, [\xi, \eta] \rangle \equiv 0.$$

g) *As a* $\psi(g)$ *it is possible to take the following function*

(14.16) $$\psi(g) = \mathrm{Ad}(g^{-1})\psi_0 - \psi_0 \in Lie(G),$$

where $\psi_0$ is some element of Lie algebra $Lie(G)$ or of their one-dimensional expansion $Lie(G_1)$, for which $\mathbf{A}^* = \mathrm{ad}(\psi_0)$. In the second case $Lie(G)$ is an ideal in $Lie(G_1)$, and an operator of the adjoint representation $\mathrm{Ad}(g^{-1})$ in (14.16) acts in $Lie(G_1)$.

◀ The equivalence a) <=> b) directly follows from the identity:

(14.17) $$L_\xi \alpha = -\kappa^{-1} d\langle L_V^\kappa P, \xi \rangle.$$

Let's prove the implication c) => b). For the proof we take $g = e^{t\xi}$ and differentiate (14.12) with respect to $t$ at $t = 0$. At the same time we will get the formula (14.14).

Let's prove the implication b) => c). It is easy to derive the differential equation for the form $\alpha(t) = (e^{t\xi})^* \alpha$:

$$\text{(14.18)} \qquad \frac{d}{dt}\alpha(t) = -\kappa^{-1} d\left\langle P, \mathrm{Ad}\left(e^{-t\xi}\right) D\xi \right\rangle,$$

where $D = \mathbf{A}^* \in \mathrm{End}(Lie(G))$. Integrating, we get:

$$\text{(14.19)} \qquad \alpha(t) - \alpha = -\kappa^{-1} d\left\langle P, \int_0^t \mathrm{Ad}\left(e^{-s\xi}\right) D\xi \, ds \right\rangle,$$

that gives c) for exponential elements of a group $g = e^{t\xi}$. So, the same is true for various finite products of such elements, i.e. for any elements of a group, owing to its connectivity.

The statement d) is obvious. For the proof of e) I will remind that $L_V^\kappa$ is a differentiation of Lie-Poisson algebra. Applying an operator $L_V^\kappa$ to the Poisson brackets

$$\text{(14.20)} \qquad \{H_\xi, H_\eta\} = H_{[\xi,\eta]}, \ H_\zeta = \langle P, \zeta \rangle, \ \zeta \in Lie(G),$$

and considering a), we get a relation

$$\text{(14.21)} \qquad \langle P, D[\xi,\eta] - [D\xi,\eta] - [\xi, D\eta] \rangle \equiv -\langle b, [\xi,\eta] \rangle.$$

Further, let $b = 0$. Consider, that the kernel of momentum mapping as the Lie algebra $Lie(G)$ homomorphism into the Lie algebra of functions is null. Therefore the operator $D$ is a differentiation of Lie algebra:

$$\text{(14.22)} \qquad D[\xi,\eta] = [D\xi,\eta] + [\xi, D\eta].$$

Now, let $b \neq 0$. Owing to the agreement, that in this case the Hamiltonians $H_\zeta = \langle P, \zeta \rangle$ are not nonzero constants, both members of equality (14.21) should be zero. From here follows (14.22), and also the statement f). The formula (14.14) was already proved above.

Let's prove the statement g). If $D$ is interior differentiation then

$$\text{(14.23)} \qquad D = \mathrm{ad}(\psi_0), \ \psi_0 \in Lie(G).$$

If $D$ is no interior differentiation it is possible to construct one-dimensional expansion $Lie(G_1)$ of a Lie algebra, adding to $Lie(G)$ a new basis vector $\psi_0$, with new commutators

$$\text{(14.24)} \qquad [\psi_0, \xi] = D\xi, \ \xi \in Lie(G).$$

In both cases we have

$$\text{(14.25)} \qquad D = \mathrm{ad}(\psi_0).$$

So, we can calculate an integral in (14.19):

$$\int_0^t \mathrm{Ad}\left(e^{-s\xi}\right) D\xi ds = -\int_0^t e^{-s\cdot \mathrm{ad}(\xi)} \mathrm{ad}(\xi)\psi_0 ds = \left(\mathrm{Ad}(e^{-t\xi})-1\right)\psi_0.$$

Hence, the formula (14.16) is proved for all exponential elements of group, that is why and for all finite products of such elements. Owing to connectivity of group $G$ the formula is true for all elements of group. ▶

**Corollary.** *A dynamic symmetry is reduced to usual symmetry if the differentiation $D=\mathbf{A}^*$ of the Lie algebra $Lie(G)$ is interior.*

Really, in this case $\mathbf{A} = \mathrm{ad}^*(\psi_0)$, where $\psi_0 \in Lie(G)$, and the defining form is represented in the form of the sum

(14.26) $$\alpha = \alpha_{inv} + \kappa^{-1} d\langle P, \psi_0\rangle,$$

where $\alpha_{inv}$ is a $G$-invariant form. Passing in a «moving frame», we get a system with the defining form $\alpha_{inv}$.

The Hamiltonian dynamic symmetry is the special case of the defined above symmetry of general type because the Hamilton fields $X_i = \mathbf{J}\, dH_i$ together with the dynamic field $V$ satisfy relations (14.2). But the theorem 1 is in this case applicable, and a problem arises about the connection of functions $F_j = \alpha(X_j)$ with the momentum $P$.

**Theorem 3.** In the case of Hamilton dynamic symmetry the vector function $F = \alpha(\mathbf{J}\, dP)$ expresses through the momentum by the formula

(14.27) $$F = -\kappa^{-1}\left((\mathbf{A}-\kappa)P + b\right).$$

**Examples.**
1. We take the system describing a falling of a rain drop[33] (force of resistance is considered proportional to velocity):

(14.28) $$\alpha = pdq + \kappa^{-1} dH, \quad H(p,q) = p^2/2 + q.$$

Also we will consider the Hamilton fields $X_1, X_2$ with Hamiltonians $H_1 = p$, $H_2 = q$. Let's note the following facts. First, the dynamic field $V$ of a system is included in the Lie algebra:

$$[V, X_1] = 0, \quad [V, X_2] = X_1 + \kappa X_2, \quad [X_1, X_2] = 0.$$

Second, the linear differential equations with constant coefficients are satisfied

---

[33] The same system describes a free particle with the fixed momentum.

$$(d/dt+\kappa)H_1 = L_v^\kappa H_1 = -1, \qquad (d/dt+\kappa)H_2 = L_v^\kappa H_2 = H_1 + \kappa H_2.$$

They give an explicit dependece on time *of two* functions $H_1, H_2$ while *only one* quasi-integral $F_1 = \alpha(X_1) = H_1 + \kappa^{-1}$ is derived from a usual symmetry field $X_1$.

**2.** Let's show, that the dynamic symmetry of the «rain drop» system (the previous example see) is not reduced to usual symmetry, in sense of a corollary to Theorem 2. In this case the group $G$ and the Lie algebra $Lie(G)$ elements are identified with triangular matrixes

(14.29) $$g = \begin{pmatrix} 1 & a & c \\ 0 & 1 & b \\ 0 & 0 & 1 \end{pmatrix}, \quad \xi = \begin{pmatrix} 0 & \dot a & \dot c \\ 0 & 0 & \dot b \\ 0 & 0 & 0 \end{pmatrix} \leftrightarrow (\dot a, \dot b, \dot c),$$

and the Poisson action of group on $M$ looks like

(14.30) $$g(p,q) = (p-b, q+a).$$

The corresponding generators are equal

(14.31) $$H_\xi = \langle P, \xi \rangle = p\dot a + q\dot b + \dot c = H_1 \dot a + H_2 \dot b + H_0 \dot c,$$

that gives the expression for differentiation $D = \mathbf{A}^*$ of a Lie algebra

(14.32) $$D : (\dot a, \dot b, \dot c) \mapsto (\dot b, \kappa \dot b, \kappa \dot c - \dot a).$$

Considering the expression for a commutator

(14.33) $$[\xi_1, \xi_2] \leftrightarrow (0, 0, \dot a_1 \dot b_2 - \dot b_1 \dot a_2),$$

we see, that this differentiation is not interior.

## §15. Remarks about the constants of motion.

Note the properties *of usual* constants of motion of CD-systems, which distinguish them from the constants of motion of Hamilton systems:

a) any constant of motion $F$ should satisfy almost everywhere a relation

(15.1) $$F = r^* F,$$

where $r$ is a mapping from section 9.3;

b) the Poisson brackets of all constants of motion vanishes.

Let $H$ be a physical Hamiltonian of CD-systems and let $F$ be a constant of motion of the corresponding *Hamilton system*:

$$\{H, F\} = 0.$$

What is the general behaviour of such a constant of motion over the quanton $N$? In this relation it is possible to tell the following.

Clearly, that $F$ is a constant of motion on $N$. Therefore or hypersurfaces $F = const$ stratify $N$ on invariant manifolds, or the function $F$ is constant on $N$. In the first case the Hamilton field $J\,dF$ does not touch to $N$. In the second case the tangency is not obligatory, but if the manifold $N$ is Lagrangian, the field $J\,dF$ always touches to $N$.

**Example.**

Let's consider the system

$$\alpha = (I_1 - h)d\varphi_1 + (I_2 - h)d\varphi_2 + \kappa^{-1}d(\omega_0 I_2),$$

in which the role of PH plays the $H = \omega_0 I_2$. Here there are two constant of motion of the Hamilton system with the Hamiltonian $H$:

$$F_1 = \cos\varphi_1, \ F_2 = I_2,$$

and $F_1$ is as well constant of motion the CD-system. The constant of motion $F_1$ takes *various* values on the Lagrangian attractor $N$: $I_1 = I_2 = h$, and the Hamilton field $J\,dF_1 = \sin\varphi_1 \cdot \partial/\partial I_1$ does not touch the attractor $N$. The constant of motion $F_2$ takes on $N$ a constant value $h$, and it is because the field $J\,dF_2 = \partial/\partial\varphi_2$ touches $N$.

### §16. A remark about the relativistic symmetry.

As it was noted above, the symmetry of Poincare type cannot be the exact dynamic symmetry of CD-system. Nevertheless, a lack of this symmetry at fundamental dynamic level does not exclude a possibility of it occurrence at the empiric level when the steady motions are analyzing basically. We will consider a simple example.

Let's take the CD-system $(\overline{M}, \theta)$ describing a free relativistic particle with a fixed momentum $p = b$:

(16.1) $\quad \theta = pdq - Hdt - dS, \ H(p) = \sqrt{p^2 + m^2}, \ S(q,t) = bq - Et, \ E = H(b),$

setting the light speed $c = 1$. Using the possibilities of an expanded space, we will make the Lorentz transformation of variables and parameters:

$$\begin{cases} t = \alpha t' + \beta q', & q = \alpha q' + \beta t' \\ H = \alpha H' + \beta p', & p = \alpha p' + \beta H' \\ E = \alpha E' + \beta b', & b = \alpha b' + \beta E' \end{cases}$$

The coefficients $\alpha, \beta$ in these formulas are expressed through the velocity $v$ of a moving frame of reference:

$$\alpha = \left(1-v^2\right)^{-\frac{1}{2}}, \quad \beta = \alpha \cdot v.$$

If to express the form $\theta$ through the shaded quantities, we will obtain the form $\theta'$ which will look like $\theta$. Let $\overline{V}'$ be a dynamic field on $\overline{M}$, expressed through the shaded quantities. As we know, it is defined by the relations

(16.2) $\quad i_{\overline{V}'} d\left(e^{\kappa t}\theta'\right) = 0, \quad dt(\overline{V}') = 1$.

Because the variable $t'$ should play a role of new dynamic parameter, we accept *a new normalisation* of a dynamic field: $dt'(\overline{V}') = 1$. The integral lines in the expanded space, representing the dynamics of a system will not vary.

As a result we get the defining relations for a CD-system:

(16.3) $\quad i_{\overline{V}'} d\left(e^{\kappa \alpha t'}\left\{e^{\kappa \beta q'}\theta'\right\}\right) = 0, \quad dt'(\overline{V}') = 1$,

with a *new dissipative constant* $\tilde{\kappa} = \alpha\kappa$ and with the form $\tilde{\theta} = e^{\kappa\beta q}\theta$.

**Lemma**. *The equations of motion of CD-system*

$$\left(\overline{M}, \theta = e^{\kappa\beta q}\left[pdq - H(p,q)dt - dS(q,t)\right]\right),$$

*with a dissipative constant $\tilde{\kappa} = \alpha\kappa$ look like*:

(16.4) $\quad \begin{cases} \dot{p} = -\partial H / \partial q - \alpha\kappa \cdot (p - \partial S / \partial q) + \kappa\beta \cdot (H + \partial S / \partial t), \\ \dot{q} = \partial H / \partial p. \end{cases}$

Because the action function $S$ in our example satisfies the Hamilton-Jacobi equation, and regardung this Lemma, we see that the motion over the attractor of system (16.4) is described by the usual Hamilton equations. Thus, the steady dynamics of system (16.1) *is Lorentz-invariant*, and the dynamics of relaxation *is non-invariant*.

# Chapter 3. Models.

## §1. Deducing and discussing the CS-model.

### 1.1. Generalization of the oscillator model.

Let's generalize the primary oscillator model from Chapter 1. The obvious shortcoming is that here is only one steady state, with energy $E = \hbar\omega_0$, instead of *a spectrum* of such states

(1.1)
$$E_n = n\hbar\omega_0.$$

Let's introduce a holomorphic coordinate $z = \sqrt{I}e^{i\varphi}$. Then the defining form will become

(1.2)
$$\alpha = \operatorname{Im}(\bar{z}dz) + \kappa^{-1} \cdot dH - dS_0,$$

with the Hamiltonian $H = \omega_0 \cdot |z|^2$ and action function $S_0 = \hbar \arg(z)$, and the point $z=0$ is excluded from a phase plane $M = \mathbb{C}$. Being guided by heuristic analogy of function $\Phi = S_0/\hbar$ with a phase of wave function, it is natural to assume, that a wave function in this case is equal to $F(z) = z$. It is confirmed by the following. If $F(z)$ is equal to a wave function of the $n$-th excited state of quantum oscillator in the Fock-Bargmann representation, $F(z) \sim z^n$, the energy of auto-oscillation will be equal to $E = n\hbar\omega_0$, because $\arg(z^n) = n \cdot \arg(z)$. To generalize we take an arbitrary entire function[34] $F(z)$, and then take $S_0 = \hbar \arg(F(z))$. Then the resulting system will be of Kahler type, with a potential

(1.3)
$$U(z) = |z|^2 - \hbar \cdot \log|F(z)|^2.$$

There is also the many-dimensional generalization, which confirms the analogy $F$ with a wave function. Let the phase space be $\mathbb{C}^N$, and $z = (z^1, z^2, ..., z^N) \in \mathbb{C}^N$. Let the potential $U$ looks like (further for simplicity $\hbar = 1$):

(1.4)
$$U(z) = \|z\|^2 - \log|F(z)|^2,$$

where $F(z)$ is an entire function. Let's take the square Hamiltonian $H(z) = \sum A_{i\bar{k}} z^i \bar{z}^k$. Then if to introduce the Hamilton operator $\mathbf{H} = \sum A_{i\bar{k}} z^i \cdot \partial/\partial z^k$ the condition $\{H, U\} = 0$ is noted in the form of the requirement that $F$ was an eigenfunction: $\mathbf{H}F = \lambda F$, $\lambda \in \mathbb{R}$. Besides, to points of a minimum of potential $U$ there correspond the points of maxima of density

(1.5)
$$\rho(z) = e^{-U} = |F(z)|^2 e^{-\|z\|^2}.$$

Nevertheless, being guided by a principle, that a wave function should be an element not dynamic, but *the statistical* theory, we will consider this analogy only as some favorable indication. However, this analogy should be explained at future.

Let's return to an initial one-dimensional case. The analogy demands that the function $F$ could vary. Hence, the function $F$ should appear in a role *of an independent dynamic variable*, and the direct product $M \times H$, where $H$ is a Hilbert space of Fock-Bargmann entire functions, should be a phase space:

(1.6)
$$\|F\|^2 = \pi^{-1} \int |F(x)|^2 \exp(-\bar{x}x) d^2x.$$

---
[34] The poles of $F$, unlike the zeroes, are *attracting* singularities, and should be excluded.

As the potential of such expanded system we take the elementary generalization of potential (1.3) coming to mind:

(1.7) $$U = \|F\|^2 + |z|^2 - \log|F(z)|^2.$$

The necessity of such choice of potential is dictated also by aspiration to reach of some symmetry between the «field» $F$ and «particle» $z$. Then, a Hamiltonian of the system we will leave the same, depending only on $z$. It corresponds to the idea that $F$ describes a vacuum with zero energy, and is confirmed by the analysis of solutions to the dynamic equations. As it will be shown in §3, the resulted system converges to the attractors coinciding at $\kappa \to 0$ with auto-oscillations described above, and at $\kappa \to 0$ the function $F$ becomes proportional to $z^n$.

The further generalizations are that. At first, we substitute *the arbitrary* Hamiltonian depending on $z$ for $H$. Secondly, let's pay attention to known identity

(1.8) $$F(z) = (F, f(\bar{z})) = \pi^{-1} \int F(x) \exp(\bar{x}z - \bar{x}x) d^2x,$$

where $f(\bar{z}; x) = \exp(x\bar{z})$ are *the coherent states*, antiholomorphicly depending on a point $z \in \mathbb{C}$. Regarding that $|z|^2 = \log\|f(\bar{z})\|^2 + const$, we can write potential (1.7) in terms of coherent states, regardless to a choice of concrete realization of the Hilbert space $H$:

(1.9) $$U = \|F\|^2 - \log \rho(z, F).$$

Here we have introduced a density in «particle» phase space:

(1.10) $$\rho(z, F) = \frac{|(F, f(\bar{z}))|^2}{\|f(\bar{z})\|^2}.$$

Further $f(\bar{z})$ will denote the coherent states of general type, antiholomorphicly depending on a point $z \in M$, defining the symplectic and Kahler structures on $M$, or, on the contrary, being constructed by such structure. For example, they can be obtained from the representations of the Lie groups [10].

As a result we obtain the «model of coherent states» or *the CS-model*. This is the Kahler CD-system with phase space $M \times H$, set by the potential (1.9), with Hamiltonian $H(z)$, and with equations of motion

(1.11) $$\begin{cases} \dfrac{dz}{dt} = \mathrm{J}\, dH + \varepsilon \operatorname{grad}_z \log \rho(z, F), & \varepsilon = \kappa/2 \\ \dfrac{dF}{dt} = -\varepsilon F + \varepsilon \dfrac{f(\bar{z})}{(f(\bar{z}), F)} \end{cases}$$

### 1.2. The analogy with Schrödinger equation.
Let's consider the case of maximal likeness of CS-model vacuum with a wave function. It will be so, if

a) the Hamiltonian looks like a covariant symbol of the self-adjoint operator:

(1.12) $$H(z) = \frac{(\mathbf{A}f(\overline{z}), f(\overline{z}))}{\|f(\overline{z})\|^2};$$

b) the corresponding unitary dynamics preserves the variety of coherent states:

(1.13) $$\exp(-it\mathbf{A})f(\overline{z}_0) = f(\overline{z}(t)).$$

At these suppositions it is easy to show, that $z(t)$ in (1.13) is a solution to the classic Hamilton equations: $z(t) = \exp(t\,\mathrm{J}\,dH)z_0$.

Let's note, that the density $\rho(z, F)$ is invariant with respect to transformations

$$\{z, F\} \mapsto \{\exp(t\,\mathrm{J}\,dH)z,\ \exp(-it\mathbf{A})F\}.$$

Therefore these transformations are included into the group of symmetries of the system. The corresponding quasi-integral $Q$ is given by the formula:

$$Q(z, F) = H(z) - (\mathbf{A}F, F).$$

From here one curious conclusion follows. Using this symmetry we can pass to a «moving frame»

(1.14) $$z = \exp(t\,\mathrm{J}\,dH)z',\quad F = \exp(-it\mathbf{A})F'.$$

According to Chapter 2, the outcome will be such as though we replace the Hamiltonian $H(z)$ by a new Hamiltonian $H(z) - Q(z, F) = (\mathbf{A}F, F)$. Accordingly, the equations of motion (1.11) will be replaced by the new:

(1.15) $$\begin{cases} \dfrac{dz}{dt} = \varepsilon\,\mathrm{grad}\log\rho(z, F), \\ \dfrac{dF}{dt} = -i\mathbf{A}F - \varepsilon F + \varepsilon\dfrac{f(\overline{z})}{(f(\overline{z}), F)}. \end{cases}$$

It is evident, that the modification of vacuum dynamics at passage from the equations (1.11) to the equations (1.15) is similar to one that happens in QM at passage from Heisenberg representation (where a wave function does not vary) to Schrödinger representation.

It is convenient to use the same terminology and for the given case. Thus, equations (1.15) describe dynamics of system in Schrödinger representation which basic peculiarity is that at the $\varepsilon \to 0$ the equation for a field $F$ turns to the Schrödinger equation.

However at $\varepsilon > 0$, thanks to the general law, the system (1.15) converges to one of the steady states (i.e. to an attractor), so $F$, eventually, is approaching to one of *eigenvectors* of the Hamilton operator $\mathbf{A}$:

$$(1.16) \qquad F(t) \approx e^{-i\omega_n t} F_n, \quad \mathbf{A} F_n = \omega_n F_n.$$

As regards a «particle» $z$, from the first equation we see, that it aspires to appear on a density $\rho(z, F)$ crest in phase space. Owing to the density invariance mentioned above, the same is valid also in initial «physical» Heisenberg representation.

The established analogy can be used in the heuristic purposes, for designing of CD-systems corresponding to quantum systems.

### 1.3. Multiparticle generalization.

The construction of CS-model by the classic Hamilton system can be considered as some analogue of the usual procedure of quantization. However, despite the general character of such «quantization», its result can be unacceptable from the physical point of view.

For example, considering $M$ as a phase space of one particle, in the multiparticle case we must use the direct product $M \times ... \times M$. Further, the analogy with QM demands to consider $F$ as a vector of tensor product $H \otimes ... \otimes H$ of Hilbert spaces. On the other hand, a field interpretation of the theory assumes that coherent states should be understood as the approximate description of solitons[35]. Therefore for the description of multiparticle system it should be sufficient the collection of coherent states $\{f_k(\overline{z}_k)\}$ and a vacuum field $F$ in *physical space*.

With keeping in mind of such restriction, for the Kahler potential and Hamiltonian of the multiparticle system it is possible to offer the following variant[36]:

$$(1.17) \qquad U = \|F\|^2 - \sum_k \log \frac{|(F, f_k(\overline{z}_k))|^2}{\|f_k(\overline{z}_k)\|^2}, \quad H = \sum_k H_k(z_k).$$

Let's note some properties of the model (1.17).

a) The distant particles (when the wave packets $f_k(\overline{z}_k)$ are not overlapping) move as independent ones.

b) Let's assume, that particles can be of different sorts $A, B, C, ...$ in the sense that vectors $f$ of particles of one sort coincide, and vectors of particles of different sorts are orthogonal: $(f_A, f_B) = ... = 0$. Then the vacuum field breaks up into the direct sum of the fields $F_A$, $F_B$,...

---

[35] May be, these solitons will arise thanks to interaction between *the fields of matter and vacuum*.

[36] For the sake of simplicity the coupling Hamiltonian of particles is omitted.

belonging to mutually orthogonal subspaces $H_A, H_B, \ldots$ of total Hilbert space $H$, and each sort of particles interacts only with the corresponding component of vacuum field. Then, as it is easy to see, the dynamic system *breaks up into a direct product of independent systems*.

c) The basic shortcoming of the given model is *the nonlocality* of dynamic equations, because of the nonlocality of action function $S_0 = \hbar \sum_k \arg(F, f_k(\overline{z}_k))$. However, as it was noted above, in a field theory we can use the local multivalued functionals leading to partial differential equations.

*Example.*

To descript the system of free particles, moving with the speed of light $c = 1$, we can take the Hamiltonians $H_k(\mathbf{z}_k) = |\mathbf{p}_k|$, and the Schrödinger coherent states $f_k(\overline{z}_k) \in L_2(\mathbb{R}^3)$:

(1.18) $\qquad f_k(\overline{\mathbf{z}}_k; \mathbf{x}) = \exp\left(-(\mathbf{x} - \overline{\mathbf{z}}_k)^2 / 2a^2\right), \quad \mathbf{z}_k = \mathbf{q}_k - ia^2 \mathbf{p}_k \in \mathbb{C}^3, \mathbf{x} \in \mathbb{R}^3.$

This choice leads to the potentials of 1-particle spaces:

$$V_k(\mathbf{z}_k) = \log\|f_k(\overline{\mathbf{z}}_k)\|^2 = -\frac{(\mathbf{z}_k - \overline{\mathbf{z}}_k)^2}{4a^2} + \frac{3}{2}\log(\pi a^2) = a^2 \cdot \mathbf{p}_k^2 + const,$$

and to standard $N$-particle symplectic form

$$\Omega = d\sum_k \operatorname{Im} \partial V_k = \sum_k d\mathbf{p}_k \wedge d\mathbf{q}_k.$$

The one-particle case of this model is considered in section 3.4.

### 1.4. Comparison with the de Broglie-Bohm mechanics.

The last model (1.17) gives a description in quantum area and contains the objects of two sorts – the particles and a field, reminding the known model of a wave-pilot, or a mechanics de Broglie-Bohm. Therefore it is interesting to make clear connection of our model (and generally CD-systems) with a mechanics of de Broglie-Bohm.

Remind, that in this mechanics the state of dynamic system is set by the coordinates $\mathbf{q}$ of a particle and by the wave function $\psi(\mathbf{x})$. The dynamical equations are the Schrödinger equation for the wave function $\psi$ and the de Broglie guidance equation

(1.19) $\qquad\qquad\qquad m\dot{\mathbf{q}} = \mathbf{p} = \partial S / \partial \mathbf{q},$

where $S(\mathbf{q}, t)/\hbar = \arg \psi(\mathbf{q}, t)$ is a phase of the wave function at the point of disposition of a particle.

Let's consider a semiclassic asymptotics of this system. According to polar representation of a wave function $\psi = \sqrt{\rho} e^{iS/\hbar}$, at $\hbar \to 0$ the system breaks up into three subsystems. The states of these three subsystems are described, accordingly, by the particle coordinates $\mathbf{q}$, action field $S(\mathbf{x})$, and density field $\rho(\mathbf{x})$.

The first subsystem $(\mathbf{q})$ obeys to the equation (1.19). The basic, driving system is the second $(S)$, for which a wave equation role plays the Hamilton-Jacobi equation

(1.20) $$\partial S/\partial t + H(\partial S/\partial \mathbf{x}, \mathbf{x}, t) = 0.$$

The dynamics of the third subsystem $(\rho)$ is set by the continuity equation

(1.21) $$\frac{\partial \rho}{\partial t} + \mathrm{div}\left(\frac{\rho}{m}\frac{\partial S}{\partial \mathbf{x}}\right) = 0.$$

Now for us the basic interest is represented by the closed subsystem $(\mathbf{q})+(S)$ described by the equations (1.19), (1.20). The matter is that the similar equations arise for the CD-system

(1.22) $$\theta = \mathbf{p}d\mathbf{q} - H(\mathbf{p},\mathbf{q},t)dt - dS(\mathbf{q},t).$$

According the Chapter 2, the Hamilton-Jacobi equation (1.20) ensures for this system the existence of invariant Lagrangian submanifold

(1.23) $$\mathbf{p} = \partial S/\partial \mathbf{q}.$$

The motion over this submanifold is described by a Hamiltonian $H$. It is obvious from the equations of motion

(1.24) $$\begin{cases} d\mathbf{p}/dt = -\partial H/\partial \mathbf{q} - \kappa \cdot (\mathbf{p} - \partial S/\partial \mathbf{q}), \\ d\mathbf{q}/dt = \partial H/\partial \mathbf{p}. \end{cases}$$

Hence, the dynamics of a subsystem $(\mathbf{q})+(S)$ can be described by combined equations (1.24), (1.20). Thus, something common at system (1.22) with a mechanics de Broglie-Bohm *is available*. To strengthen this effect, let's construct the CD-system, for which the equations (1.24), (1.20) arise as a quasi-classic approximation.

To do this, it will be useful to consider an example of a free particle:

$$H = \mathbf{p}^2/2m, \quad S = \mathbf{b}\mathbf{q} - \mathbf{b}^2 t/2m.$$

In this example the equation (1.23) describes *an attractor*[37] of the system (1.22), which consists of trajectories of a Hamiltonian system, with the Hamilton function $H$. For these trajectories a particle momentum has a preset value $\mathbf{p} = \mathbf{b}$. However to embrace *all* classic motions, it is necessary to consider the *various* solutions $S(\mathbf{q},t)$ to the Hamilton-Jacobi equation, and, even-tually, the various CD-systems. This, at first sight, a strange situation, speaks about incompleteness of our system. As well as in an example with oscillator, the system should be expanded *by vacuum degree of freedoms*, in the form of some field $F(\mathbf{x})$. For this purpose, by analogy with the CS-model, we take the system

(1.25) $$\theta' = \mathbf{p}d\mathbf{q} - H(\mathbf{p},\mathbf{q},t)dt - dS + \mathrm{Im}(\delta F, F).$$

---

[37] Note, that this «attractor» is noncompact, and breaks the assumption § 9 of Chapter 2.

Here the action $S$ means a vacuum phase at the point, where a particle is disposed

(1.26) $$S(\mathbf{q}, F) = \hbar \arg F(\mathbf{q}).$$

We see also that the action $S$ play a role similar to the role of particle-vacuum coupling Hamiltonian.

Then equations (1.24) become a part of wider set of equations including in addition the equations for a vacuum component

(1.27) $$\frac{dF}{dt} = -\varepsilon F + \frac{\hbar \varepsilon f(\mathbf{q})}{F(\mathbf{q})}, \quad \varepsilon = \kappa/2.$$

Here the vectors $f(\mathbf{q})$ are defined from equality

$$F(\mathbf{q}) = (F, f(\mathbf{q})).$$

At the heart of such formal construction two hypotheses lie:

a) a particle is described by $\delta$-shaped field $f_\mathbf{q}(\mathbf{x}) \approx c(\mathbf{q}) \cdot \delta(\mathbf{x} - \mathbf{q})$, or by general «coherent state» $f_{\mathbf{p},\mathbf{q}}(\mathbf{x}) \approx c(\mathbf{p},\mathbf{q}) \cdot \delta(\mathbf{x} - \mathbf{q})$;

b) beside a particle the vacuum field looks like the wave $F(\mathbf{x}) \approx a(\mathbf{x}) \exp(iS(\mathbf{x})/\hbar)$.

If to compare the obtained system with the system of de Broglie-Bohm, then, at first sight, a vacuum field $F(\mathbf{x},t)$ of our system should correspond to the de Broglie-Bohm wave function $\psi(\mathbf{x})$. However this variant is wrong, as the equation (1.27) has no anything common with a Schrödinger equation.

To bypass this obstacle, let's be restricted by the frameworks of CS-model. Let the particle be described by $\delta$-shaped coherent state:

$$f(\overline{z}; \mathbf{x}) = f_{\mathbf{p},\mathbf{q}}(\mathbf{x}).$$

As an example it is possible to take the wave packets (1.18), adopting their breadth $a \to 0$.

Let's demand also that the Hamiltonian looked like a covariant symbol (1.12) of the Hamilton operator $\mathbf{A}$ with a simple discrete spectrum. Then (section 2.4 see.) there is a whole series of the «spectral» solutions to our system presumably disposed on an attractor for which $F(\mathbf{x},t)$ are close *to operator* **A** *eigenfunctions*:

$$F(\mathbf{x},t) \approx F_n(\mathbf{x}).$$

The closeness to eigenfunctions is guaranteed by a supposed smallness of the dissipative constant.

So, though the equation for a field $F$ is not a Schrödinger equations, the dynamics *of system as a whole* is that, that in steady states this field is a solution to the stationary Schrödinger equation. Hence, the steady motions of our system differ from corresponding solutions in the de Broglie-Bohm mechanic only by a phase factor:

$$\psi(\mathbf{x},t) \sim e^{-i\omega_n t} F_n(\mathbf{x}).$$

Full coincidence cannot be, at least for the reason, that the de Broglie-Bohm system does not describe a relaxation to the quantum steady-states.

## §2. The analysis of general CS-model.

### 2.1. Further generalization and the equations of motion.

The Kahler potential of an oscillator from section 1.1 can be rewritten without explicit use of coherent states $f(\overline{z})$. To derive the corresponding generalization it is enough to note, that expression

$$(2.1) \qquad \frac{|(F,f(\overline{z}))|^2}{\|f(\overline{z})\|^2} = \frac{|F(z)|^2}{\|f(\overline{z})\|^2} = \langle F,F \rangle_z$$

is possible to consider as Hermitian square of holomorphic section $F: z \mapsto F(z) = (F, f(\overline{z}))$ of some fiber bundle. Therefore the oscillator potential will be noted in the form

$$(2.2) \qquad U(z,F) = \|F\|^2 - \log \langle F,F \rangle_z.$$

Also it is possible to start with *holomorphic Hermitian line bundle* $L \to M$, assuming, that $L$ admits a rich enough set of holomorphic sections $F \in \Gamma(L)$.

For example, if $M = \mathbb{C}P^n$, the holomorphic sections of $L$ are identified with the homogeneous polynomials $F(s)$, $s \in \mathbb{C}^{n+1}$ of a fixed degree $m$, and Hermitian metric in fibers is set by the formula $\langle F,F \rangle_{p(s)} = |F(s)|^2 \cdot \|s\|^{-2m}$. Here $p(s) \in \mathbb{C}P^n$ is a projection to projective space. At $n=1$ such system *the spin $m/2$* describes. It will be considered in section 3.5.

Apparently from this example, the sections are convenient to set by homogeneous functions on the corresponding principal $\mathbb{C}^*$-bundle $p: E \to M$, and homogeneity degree $m$ de-pends on the line bundle $L$. In such realization the Hermitian square of section is calculated by means of a homogeneous function $B(\lambda s, \overline{\lambda s}) = |\lambda|^{2m} \cdot B(s, \overline{s}) > 0$:

(2.3) $$\langle F, F \rangle_z = \frac{|F(s)|^2}{B(s,\overline{s})}, \quad z = p(s).$$

Note, that the Kahler and symplectic structures on $M$ are completely defined by the Hermitian metric on $L$. Namely, at the fixed section $F$, the expression (2.2) is necessary to consider as a local potential of the Kahler metric on $M$. This metric, obviously, does not depend on a section. Besides, the Hilbert space $H \subseteq \Gamma(L)$ of holomorphic sections is defined by the expression for quadrate of norm of the section

(2.4) $$\|F\|^2 \sim \int_M \langle F, F \rangle_z \Omega^n.$$

Further, we will assume, that the Lie group $G$ acts on the fiber bundle $E$ by holomorphic automorphisms, preserving the Hermitian metric on $L$. Last requirement is equivalent to a condition of $G$-invariance of the function $B(s,\overline{s})$. Owing to these suppositions, the unitary representation of group acts in $H$:

(2.5) $$\bigl(T(g)F\bigr)(s) = F(g^{-1}s).$$

It is convenient to define the «coherent states» $E(\overline{s}) \in H$ by the relation

(2.6) $$F(s) = \bigl(F, E(\overline{s})\bigr),$$

as everywhere defined homogeneous antiholomorphic vector functions $E(\overline{s}) \in H$ on a principal bundle $E$. They are connected with the $f(\overline{z})$ by formula

(2.7) $$f(\overline{z}) = E\bigl(\overline{s(z)}\bigr),$$

where $s: M \to E$ is any local holomorphic section, and $z = p(s)$. Further we will use $f(\overline{z})$ or $E(\overline{s})$, depending on convenience. It should be noted, that the various vector functions $f(\overline{z})$ obtained by this way, can differ only by a scalar antiholomorphic multiplier. Also, these functions $f(\overline{z})$ *could be not globally defined*. For example, in the case $M = \mathbb{CP}^1$, the functions $f(\overline{z})$, which is defined everywhere on $M$, don't exist.

Directly from definition of coherent states their important property follows[38]:

(2.8) $$E\bigl(\overline{gs}\bigr) = T(g)E(\overline{s}).$$

---

[38] It can be invalid for the functions $f(\overline{z})$.

Further we will consider a special case, apparently, closer to physics, when (at a due normalization) the equality is satisfied[39]

(2.9) $$B(s,\overline{s}) = \|E(\overline{s})\|^2,$$

meaning, that coherent states $E(\overline{s})$ set the antiholomorphic isometric embedding of $M$ in a projective Hilbert space $P(H)$.

Now we can formulate the «one-particle» CD-system connected with given coherent states[40]. In the Kahler case which here is considered, the defining form looks like:

(2.10) $$\alpha = \operatorname{Im}\partial U + \kappa^{-1}dH,$$

and a phase space is the complement of set of incident pairs

$$N = \{(z,F): F(s) = 0,\ z = p(s)\} \subset M \times H.$$

Let's consider, that the Hamiltonian $H$ does not depend on a vacuum component $F$. To write out the equations of motion, it is convenient to use that our system locally is composite. We can take the direct product potential

(2.11) $$U' = \|F\|^2 + V(z),$$

where $V(z)$ is the local potential of Kahler metric on $M$,

(2.12) $$V(z) = \log \|f(\overline{z})\|^2,$$

and Hamiltonian

$$H' = H - \kappa S_0,$$

containing an action function

(2.13) $$S_0 = \arg(F, f(\overline{z})).$$

Then from formulas of §5 Chapters 2 we obtain the equations of motion:

---

[39] This equality is characteristic for Berezin's quantization of homogeneous Kahler manifolds, see [9].

[40] The adjective «one-partial» should not be taken here literally, because of absence of restrictions on physical interpretation of manifold $M$.

(2.14)
$$\begin{cases} V_{i\bar{k}} \dfrac{dz^i}{dt} = iH_{\bar{k}} - \varepsilon V_{\bar{k}} + \varepsilon \dfrac{(f_{\bar{k}},F)}{(f,F)} = iH_{\bar{k}} + \varepsilon \dfrac{\partial}{\partial \bar{z}^k} \log \dfrac{|(F,f(\bar{z}))|^2}{\|f(\bar{z})\|^2} \\ \dfrac{dF}{dt} = -\varepsilon F + \varepsilon \dfrac{f}{(f,F)} \end{cases}$$

Here $\varepsilon = \kappa/2$ $V_{\bar{k}} = \partial V / \partial \bar{z}^k$, $V_{i\bar{k}} = \partial^2 V / \partial z^i \partial \bar{z}^k$ etc.

Let's note the important fact. According to the equations (2.14) the singularity $(F,f) = 0$ is repelling because near to it the inequality holds

$$\frac{d}{dt}|(F,f)|^2 > 0.$$

## 2.2. The symmetries and quasi-integrals.

The potential (2.2) and corresponding Kahler structure on $M \times H$ are invariant with respect to transformations $\{z,F\} \mapsto \{gz, T(g)F\}$. As a consequence, the action of group $G$ is Poissonian, and the momentum mapping $P: M \times H \to Lie^*(G)$ is set by the formula:

(2.15)
$$\langle P(z,F); \xi \rangle = i(\tilde{\xi}) \operatorname{Im} \partial U = A(z,\xi) - (\mathbf{A}(\xi)F, F).$$

Here $i(\tilde{\xi})$ is an interior multiplication by a vector $\tilde{\xi}$ of an infinitesimal group action on space $M \times H$, and $\xi \in Lie(G)$ is an element of a Lie algebra. The function $A(z,\xi)$ in this formula is *a covariant symbol* of the self-adjoint operator

$$\mathbf{A}(\xi) = iT_*(\xi) = i \left. \frac{d}{dt} \right|_{t=0} T(e^{t\xi})$$

with respect to system of coherent states $E(\bar{s})$:

(2.16)
$$A(z,\xi) = \frac{(\mathbf{A}(\xi)E, E)}{\|E\|^2}.$$

Note the infinitesimal analogue of relation (2.8)

(2.17)
$$\tilde{\xi}.E = -i \cdot \mathbf{A}(\xi)E.$$

Considering, that the fundamental vector field $\tilde{\xi}(s)$ on fiber bundle $E$ is divided into holomorphic and antiholomorphic components:

(2.18)
$$\tilde{\xi} = \xi' + \xi'' = \xi^a(s) \frac{\partial}{\partial s^a} + \overline{\xi^a(s) \frac{\partial}{\partial s^a}},$$

we can rewrite (2.17) in the form

(2.19) $$\xi''.E = -i \cdot \mathbf{A}(\xi)E.$$

Note also the property of invariance of a symbol

(2.20) $$A(z,\xi) = A(gz, \mathrm{Ad}(g)\xi), \quad g \in G.$$

It can be calculated by the formulas:

(2.21) $$A(z,\xi) = -i\frac{d}{dt}\bigg|_{t=0} \log B(e^{t\xi}s, \overline{s}) = -i\xi'.\log B(s,\overline{s}).$$

Owing to a decomposability of the Kahler metric on $M \times H$, from the equality (2.15) follows that the group actions on $M$ is Poissonian, and the corresponding momentum mapping $P(z): M \to Lie^*(G)$ is defined by the same symbol $A(z,\xi)$:

(2.22) $$\langle P(z); \xi \rangle = A(z,\xi).$$

The obvious symmetry of given system is the group $U(1)$ of phase transformations

(2.23) $$\{z, F\} \to \{z, e^{it}F\}.$$

To theses transformations correspond the quasi-integral

(2.24) $$Q_1 = \|F\|^2 - 1.$$

Other symmetry $G_H \subseteq G$ is represented by transformations preserving both the potential $U$ and Hamiltonian $H$. To each one-parameter group of transformations from $G_H$

(2.25) $$\{z, F\} \to \{e^{t\xi}z, T(e^{t\xi})F\}, \quad \xi \in Lie(G_H),$$

there corresponds the quasi-integral $Q_2$, being a component of momentum (2.15):

(2.26) $$Q_2 = A(z,\xi) - (\mathbf{A}(\xi)F, F) = \langle P(z,F); \xi \rangle.$$

It is interesting, that the equality $Q_1 = 0$ is similar to a condition of normalization of wave function in quantum mechanics, and the equality $Q_2 = 0$ states about coincidence of corresponding momentums of «particle» $z$ and «vacuum field» $F$.

### 2.3. The equations for search of Lie solutions.

According to a definition, the Lie solutions to the dynamic equations (2.14), with the symmetry group $U(1) \times G_H$, look like

(2.27) $$z(t) = e^{t\xi}z_0, \quad F(t) = e^{i\omega t}T(e^{t\xi})F_0,$$

where $\xi \in L(G_H)$, and $\omega$ is a real parameter.

According to §8 Chapter 2, the equations for the Lie solutions for a subgroup $G_H$, are reduced to one relation at the point $\{z_0, F_0\}$:

(2.28) $$\kappa\alpha = d(\omega Q_1 + Q_2).$$

In more details, we will obtain the equalities (a lower index 0 at $z$ and $F$ will be omitted further):

(2.29) $$\left(\mathbf{A}(\xi) - \omega + i\varepsilon\right)F = \frac{i\varepsilon E}{(E,F)},$$

(2.30) $$\left.\frac{\partial}{\partial z^k}\right|_{F=const}\left(H(z) - A(z,\xi) + i\varepsilon \cdot \log\frac{|(E,F)|^2}{\|E\|^2}\right) = 0.$$

It is easy to see, that necessary and sufficient condition of resolvability of the equation (2.29) with respect to the vector $F$ is the equality

(2.31) $$\mathrm{Re}(\mathbf{R}E, E) = 0,$$

where $\mathbf{R} = \mathbf{R}_{\omega - i\varepsilon} = \left(\mathbf{A}(\xi) - \omega + i\varepsilon\right)^{-1}$ is the resolvent of operator $\mathbf{A}(\xi)$. If this condition is true, the general solution to (2.29) looks like

(2.32) $$F = \mu \cdot \mathbf{R}E,$$

where $\mu \in \mathbb{C}$ is the arbitrary factor ensuring a normalization $\|F\|^2 = 1$. Hence, the vector $F$ is excluded from the equations (2.29), (2.30). As a result of such elimination we obtain the

**Theorem 1.** a) *The search of solutions* (2.27) *to dynamic equations* (2.14) *is reduced to a solution to finite $G_H$-invariant set of equations with respect to unknowns* $z \in M$, $\xi \in Lie(G_H)$, $\omega \in \mathbb{R}$:

(2.33) $$\mathrm{Re}\, R(z, \xi, \omega - i\varepsilon) = 0,$$

(2.34) $$\frac{\partial}{\partial z^k}\{H(z) - A(z,\xi) + i\varepsilon \cdot \log R(z,\xi,\omega - i\varepsilon)\} = 0,$$

*where* $R(z,\xi,\omega - i\varepsilon) = (\mathbf{R}E, E)/\|E\|^2$ *is the resolvent* $\mathbf{R}$ *symbol.*

b) *If $\{z, \xi, \omega\}$ is the solution to the equations (2.33), (2.34), then $\{gz, \mathrm{Ad}(g)\xi, \omega\}$ is also the solution, at any $g \in G_H$.*

c) *From the equalities (2.34) follows*

(2.35) $$\mathrm{Re}\frac{1}{R(z,\xi,\omega - i\varepsilon)} = A(z,\xi) - \omega.$$

*Therefore a consequence of equations (2.33), (2.34) the equality is*

(2.36) $$A(z,\xi) = \omega.$$

*The system (2.33), (2.34) is equivalent to system (2.36), (2.34).*

d) *Elimination of frequency $\omega$ from the system (2.33), (2.34) with the help of (2.36) gives an $G_H$-invariant finite set of equations on the $G_H$-space $M \times L(G_H)$, which is reduced to a system on the set of orbits $\{M \times Lie(G_H)\}/G_H$.*

e) *From the equations (2.33), (2.34) follows, that on a corresponding solution to the dynamic equations the momentums (2.24), (2.26) vanish.*

Let's present now some properties of resolvent symbol $R(z,\xi,\omega-i\varepsilon)$.

**Theorem 2.** *The function $R(z,\xi,\omega-i\varepsilon)$:*

a) *is $G$-invariant:*

(2.37) $$R(z,\xi,\omega-i\varepsilon) = R\big(gz, \mathrm{Ad}(g)\xi, \omega-i\varepsilon\big);$$

b) *satisfies the differential equation (it is used for a deriving of (2.35)):*

(2.38) $$-i \cdot \xi'.R + \big(A(z,\xi) - \omega + i\varepsilon\big) \cdot R = 1,$$

*(here $\xi' = \xi^k(z)\partial/\partial z^k$ is a holomorphic component of the fundamental symmetry field $\tilde{\xi}(z)$ on $M$);*

c) *satisfies the inequalities*

$$-\varepsilon^{-1} \leq \mathrm{Im}\, R < 0;$$

d) *admits the integral representation:*

(2.39) $$R(z,\xi,\omega-i\varepsilon) = -i\int_0^\infty e^{-(i\omega+\varepsilon)t} \frac{B\big(e^{t\xi}s,\overline{s}\big)}{B(s,\overline{s})} dt.$$

## 2.4. Connection with a problem on eigenvalues.

The sense of equality (2.36) clears up in that specific case which here will be considered. Namely, we will assume in addition, that the Hamiltonian has a special form (for example, this is the property of oscillator model):

(2.40) $$H(z) = A(z,\xi_0),$$

where $\xi_0 \in Lie(G)$ is some fixed element of a Lie algebra. To such a Hamiltonian naturally corresponds *the Hamilton operator* $\mathbf{H} = \mathbf{A}(\xi_0)$.

It appears, in this case the equations (2.33), (2.34) have the solutions for which $\xi = \xi_0$. As a consequence, the (2.36) turns into *the Planck formula*

(2.41) $$H(z) = \omega.$$

**Theorem 3.** *Let $\omega_n$ be the isolated simple eigenvalue of an operator $\mathbf{A}(\xi_0)$:*

(2.42) $$\mathbf{A}(\xi_0)F_n = \omega_n F_n, \quad \|F_n\| = 1,$$

*and let $z_n$ be a critical point of the density*

(2.43) $$\rho_n(z) = |(F_n, f(\bar{z}))|^2 / \|f(\bar{z})\|^2 = \langle F_n, F_n \rangle_z.$$

*Let this critical point be nondegenerate in quotient space $M/G_{\xi_0}$, where $G_{\xi_0}$ is the isotropy group of element $\xi_0$:*

$$G_{\xi_0} = \{g \mid \mathrm{Ad}(g)\xi_0 = \xi_0\}.$$

*Then, at small enough $\varepsilon = \kappa/2 > 0$ there is a solution to the equations (2.33), (2.34) for which $\xi = \xi_0$, $\omega \approx \omega_n$, $z \approx z_n$. On such «spectral» solution the Hamiltonian (2.40) approximately coincides with a «quantum average»:*

(2.44) $$H(z) \approx (\mathbf{A}(\xi_0)F_n, F_n).$$

*The corresponding Lie solution to the dynamic equations looks like*

(2.45) $$z(t) \approx e^{it\xi_0} z = \exp(t\,\mathrm{J}\,dH)z, \quad F(t) \approx e^{i\alpha} F_n, \quad \alpha = const,$$

*i.e. the Hamiltonian $H(z) = A(z, \xi_0)$ is close to a physical Hamiltonian for this series of solutions*.

◄ The equations (2.33), (2.34) in this case are reduced to the following ones:

(2.46) $$\mathrm{Re}\, R(z, \xi_0, \omega - i\varepsilon) = 0,$$

(2.47) $$\frac{\partial}{\partial z^k} R(z, \xi_0, \omega - i\varepsilon) = 0.$$

Using the properties of a resolvent, at $\omega \approx \omega_n$ and at small $\varepsilon$, we get:

(2.48) $$R(z, \xi_0, \omega - i\varepsilon) = \frac{\rho_n(z)}{\omega_n - \omega + i\varepsilon} + \varphi(z, \omega, \varepsilon),$$

where $\rho_n(z) = |(F_n, f(\overline{z}))|^2 / \|f(\overline{z})\|^2$ is a density function on $M$ corresponding to an eigenvector $F_n$, $\varphi(z,\omega,\varepsilon)$ is some function bounded at $\varepsilon \to 0$, $\omega \to \omega_n$. Regarding, that at endpoints of an interval $(\omega_n - \varepsilon, \omega_n + \varepsilon)$ the function $\psi(\omega) = \operatorname{Re}(\omega_n - \omega + i\varepsilon)^{-1}$ accepts values $\psi = \pm 1/2\varepsilon$ we come to conclusion, that at small enough $\varepsilon$ there should be a solution $\omega = \omega(z)$ to the equation (2.46) in this interval. Further, the representation (2.48) supposes the $z$-differentiation, therefore the equations (2.47) are noted in the form

$$(2.49) \qquad \frac{\partial \rho_n(z)}{\partial z^k} = -(\omega_n - \omega + i\varepsilon)\frac{\partial}{\partial z^k}\varphi(z,\omega,\varepsilon).$$

We are interested in the resolvability of these equations with respect to $z \in M$, $z \approx z_n$, at $\omega = \omega(z)$. As a right member in (2.49) has an order $O(\varepsilon)$, the existence of a nondegenerate critical point $z_n$ of the density $\rho_n$ will be sufficient to derive the solution by a theory of perturbations. However, the $G_H$-invariance of a density forms an obstacle. In this case $G_H$ coincides with an isotropy group $G_{\xi_0}$ of an element $\xi_0 \in Lie(G)$. But functions $\varphi(z,\omega,\varepsilon)$ and $\omega(z)$ have the same symmetry. Therefore the equation (2.49) (where after $z$-differentiation should be set $\omega = \omega(z)$) admits a projection on the quotient space $M/G_{\xi_0}$, where the image $z'_n$ of $z_n$ becomes a nondegenerate critical point of a density, owing to the theorem condition. Therefore in quotient space there is a solution $z' \in M/G_{\xi_0}$ close to $z'_n$. To solution $z'$ there corresponds the whole $G_{\xi_0}$-orbit of solutions on $M$, and as a solution $z \in M$, close to $z_n$, it is possible to take the point of an orbit nearest to $z_n$. The remaining statements of theorem are checked up directly. ▶

**Remarks**. 1) At the «steady-state», described by the theorem, the vacuum field $F$ is not strictly static, and slightly fluctuates round an average value, in agreement with (2.45).

2) At small $\varepsilon$ the equation (2.46) has, in general, the «superfluous» roots[41] $\omega = \omega(z)$ not close to eigenvalues $\omega_n$. To these roots can correspond the Lie solutions which are not «spectral». If to consider, that a problem on eigenvalues (2.42) is similar to corresponding problem in quantum mechanics, then there is a problem on reasons of non-observability of such «not spectral» solutions. Their *instability* can be an obvious reason. Probably also, that some «not spectral» solutions actually are observed, because in the quantum world there is the whole set of worthy *unstable* objects[42].

---

[41] For example, between the roots $\omega_1(z), \omega_2(z)$ being close to eigenvalues $\omega_1, \omega_2$ must be one more root, because of the same signs of derivatives $\frac{\partial}{\partial \omega}\operatorname{Re} R(z,\xi_0,\omega-i\varepsilon)$ in this roots at $\varepsilon \to 0$.

[42] For example, stable electron has unstable analogues in the form of mu-meson and tau-lepton.

3) In the case of simple discrete spectrum of operator $\mathbf{A}(\xi_0)$ the deviation from eigenvalues $\omega_n$ has a second order $O(\varepsilon^2)$:

(2.50) $$\omega - \omega_n \approx \varepsilon^2 \cdot \sum_{j \neq n} \frac{\rho_j \rho_n^{-1}}{\omega_j - \omega_n}.$$

## §3. Special cases of CS-model.

### 3.1. The matrix model.

**1.** Let's consider the Kahler system having a phase space $M$ (with complex coordinates $\psi_1, \psi_2, \ldots, \psi_N, \chi_1, \chi_2, \ldots, \chi_N$) as the complement in $\mathbb{C}^{2N}$ the zero set of the bilinear form

(3.1) $$\chi \psi = \chi_1 \psi_1 + \chi_2 \psi_2 + \ldots + \chi_N \psi_N.$$

Here $\psi$ is a column, and $\chi$ is a line of length $N$. The system is set by the potential

$$U = \chi \chi^* + \psi^* \psi - \log |\chi \psi|^2,$$

and Hamiltonian

$$H = \psi^* \mathbf{A} \psi.$$

For example, it is possible to take as $\mathbf{A}$ a Hamilton operator of finite-dimensional quantum system. Equivalently, it is possible to take the potential $U' = \chi \chi^* + \psi^* \psi$ and multivalued Hamiltonian

(3.2) $$H' = \psi^* \mathbf{A} \psi + \kappa \cdot \arg(\chi \psi) = H - \kappa \cdot S_0,$$

including the action function $S_0 = -\arg(\chi \psi)$. The equations of motion have the following form

(3.3) $$\begin{cases} d\psi / dt = i\mathbf{A}\psi - \varepsilon \psi + \varepsilon (\chi / \chi \psi)^* \\ d\chi / dt = -\varepsilon \chi + \varepsilon (\psi / \chi \psi)^* \end{cases}$$

where $\varepsilon = \kappa/2$. At $\varepsilon \to 0$ they are reduced to the first equation becoming the Schrödinger equation, and admit the *spectrum* of auto-oscillations. Namely, to every normalized eigenvector $\mathbf{x}$,

$$\mathbf{A}\mathbf{x} = \omega_m \cdot \mathbf{x},$$

corresponds the solution

(3.4) $$\psi(t) = e^{i\omega_m t} \mathbf{x}, \quad \chi(t) = e^{i\alpha} \mathbf{x}^*,$$

where $e^{i\alpha}$ is an arbitrary constant phase factor.

**2.** Let's consider the symmetry of system. To phase transformations

$$\psi \to e^{i\alpha}\psi, \quad \chi \to e^{i\beta}\chi$$

correspond the quasi-integrals

$$Q_1 = \psi^*\psi - 1, \quad Q_2 = \chi\chi^* - 1.$$

The reduction by this symmetry gives the reduced system defined on a subset $\chi\psi \neq 0$ of product of projective spaces $P(H) \times P(H^*)$, with the reduced potential and Hamiltonian ($\chi, \psi$ are considered as homogeneous coordinates):

$$\tilde{U} = \log\frac{\psi^*\psi \cdot \chi\chi^*}{|\chi\psi|^2}, \quad \tilde{H} = \frac{\psi^* A \psi}{\psi^*\psi}.$$

After the given reduction the solutions (3.4) turn into attracting equilibrium states (for a simple spectrum of $A$, these are isolated points). They can be described as intersection of critical manifolds of reduced potential and Hamiltonian[43].

The complete symmetry group $G$, besides the specified above phase transformations, includes also the unitary transformations generated by various self-adjoint matrixes $C$, commuting with the Hamiltonian $A$:

(3.5) $$\psi \to \exp(isC)\psi, \quad \chi \to \chi\exp(-isC), \quad s \in \mathbb{R}.$$

To these transformations there corresponds a set of quasi-integrals

$$Q_C = \psi^* C \psi - \chi C \chi^*.$$

**3.** Let the matrix $A$ has a simple spectrum. Then it is possible to consider $A$, $C$ as diagonal matrixes

(3.6) $$A = \text{diag}(\omega_1, ..., \omega_N), \quad C = \text{diag}(c_1, ..., c_N).$$

Let's describe *the reduced* system, which corresponds to action (3.5) of the group $\mathbb{T}^N$. As coordinates on the corresponding quotient set $\tilde{M}$ of zero momentum is possible to take the coordinates of action-angle type:

$$\chi_k \psi_k = I_k e^{i\varphi_k}, \quad I_k \geq 0, \quad k = 1, 2, ..., N.$$

In these coordinates the defining form of reduced system on $\tilde{M}$ will become:

---

[43] Such intersection is an invariant set.

(3.7) $$\tilde{\alpha} = \sum_k \left( I_k d\varphi_k + \kappa^{-1}\omega_k dI_k \right) - d\arg\sum_k I_k e^{i\varphi_k}.$$

**4.** The reduced system (3.7) has the residual phase symmetry

(3.8) $$\varphi_k \to \varphi_k + t, \quad k = 1, 2, \ldots, N,$$

with the quasi-integral generator

$$Q = I_1 + I_2 + \ldots + I_N - 1.$$

Let's find all Lie solutions $\left( I_k(t) = I_k,\ \varphi_k(t) = \varphi_k + \lambda\kappa t \right)$ corresponding to this symmetry. The problem is reduced to searching of nonempty subsets $Nul(\tilde{\alpha} - \lambda dQ)$. Denote, for brevity,

$$F = \sum I_k e^{i\varphi_k}.$$

Then

$$d\arg F = \operatorname{Im}\sum F^{-1} e^{i\varphi_k} \left( dI_k + iI_k d\varphi_k \right),$$

and, equating to zero the coefficients of the form $\tilde{\alpha} - \lambda dQ$, we obtain the equations for unknown quantities $\lambda,\ I_k,\ \varphi_k,\ k = 1, 2, \ldots, N$:

$$\begin{cases} I_k \left( 1 - \operatorname{Re}(F^{-1} e^{i\varphi_k}) \right) = 0, \\ \operatorname{Im}(F^{-1} e^{i\varphi_k}) = \kappa^{-1}\omega_k - \lambda. \end{cases}$$

Note, that for each nonzero $I_k$ the equality is satisfied

$$1 + \left( \kappa^{-1}\omega_k - \lambda \right)^2 = |F|^{-2}.$$

In view of that all frequencies $\omega_k$ is different, such case is possible if only no more than two nonzero values $I_k$ exist. From here we obtain two series of solutions on $\tilde{M}$.

1) Solutions with unique nonzero value $I_a$ (all remaining $I_k = 0$):

$$\lambda = \kappa^{-1}\omega_a, \quad I_a = 1.$$

2) Solutions with two nonzero values $I_a, I_b$:

$$\lambda = \frac{\omega_a + \omega_b}{2\kappa}, \quad I_a = I_b = \frac{1}{2}, \quad e^{i(\varphi_b - \varphi_a)} = \frac{1 + \frac{i}{2\kappa}(\omega_b - \omega_a)}{1 - \frac{i}{2\kappa}(\omega_b - \omega_a)}.$$

Solutions to the first series correspond to the solutions (3.4) specified above and describe the steady limit cycles. Solutions to the second series are unstable.

To all these solutions correspond Lie solutions to *primary* system on $M$, having the form

(3.9) $$\psi(t) = \exp(it\mathbf{C})\mathbf{x}e^{i\omega t}, \quad \chi(t) = \mathbf{y}\exp(-it\mathbf{C}).$$

They also consist of two series.

1) Solutions (3.4) with $\mathbf{C} = \mathbf{0}$. These solutions fill *two-dimensional tori*.

2) Solutions, having vectors $\mathbf{x}$, $\mathbf{y}$ with two nonzero components with numbers $m = a, b$. Besides, the unique nonzero diagonal elements of matrix $\mathbf{C}$ are the elements $c_a, c_b$. These solutions fill *three-dimensional tori*.

For the solutions to second series we have

$$|x_a|^2 = |x_b|^2 = |y_a|^2 = |y_b|^2 = 1/2,$$
$$\omega = (\omega_a + \omega_b)/2, \quad c_a = (\omega_a - \omega_b)/4 = -c_b,$$

and phases of components of vectors $\mathbf{x}$, $\mathbf{y}$ should satisfy the relation

$$y_b x_b = y_a x_a \frac{1 + \frac{i}{2\kappa}(\omega_b - \omega_a)}{1 - \frac{i}{2\kappa}(\omega_b - \omega_a)}.$$

In other respects these phases can be arbitrary.

**5.** Let's check up, that solutions to the first series are steady, and the solutions of the second series are unsteady. For this we fulfil reduction by residual symmetry (3.8). Note the expression for defining form (3.7), restricted on manifold of zero momentum:

(3.10) $$\tilde{\alpha}|_{Q=0} = \sum_{k \neq a}\left(I_k d\phi_k + \kappa^{-1}\mu_k dI_k\right) - d\arg\Phi.$$

Here $1 \leq a \leq N$ is an arbitrarily selected index number,

$$\phi_k = \varphi_k - \varphi_a, \quad \mu_k = \omega_k - \omega_a, \quad \Phi = 1 + \sum_{k \neq a} I_k(e^{i\phi_k} - 1).$$

It is possible to consider a collection $\left\{(I_k, \phi_k) \mid k \neq a, I_k \geq 0, \sum_{k \neq a} I_k \leq 1\right\}$ as local coordinates of action-angle type on «completely reduced» space $\tilde{M}$, and the form (3.10) as the defining form of reduced system. Write out the equations of motion:

(3.11) $$\begin{cases} \dfrac{dI_k}{dt} = -\kappa I_k \cdot \left(1 - \mathrm{Re}\left(\Phi^{-1}e^{i\phi_k}\right)\right), \\ \dfrac{d\phi_k}{dt} = \mu_k - \kappa\,\mathrm{Im}\left(\Phi^{-1}(e^{i\phi_k} - 1)\right). \end{cases}$$

Let's consider a solution from the first series on $\tilde{M}$, for which $I_a = 1$. The solutions close to this solution, are characterized by a smallness of all action variables

$$I_k \approx 0.$$

Therefore the projections on $\tilde{\tilde{M}}$ of these close solutions are described by the equations

$$\begin{cases} \dfrac{dI_k}{dt} \approx -\kappa I_k \cdot (1 - \cos\phi_k), \\ \dfrac{d\phi_k}{dt} \approx \mu_k - \kappa \sin\phi_k. \end{cases}$$

It is clear, that at $|\mu_k| > \kappa$ there is an asymptotic stability:

$$\lim_{t\to\infty} I_k(t) = 0.$$

But if the difference of eigenvalues $\mu_k = \omega_k - \omega_a$ is small, there is an interesting phenomenon *of a relaxation delay*. Really, in this case the phase $\phi_k$ tends to the small constant values $\phi_k \approx \mu_k / \kappa$, and for action variables the exponential law of decrease is obtained

$$I_k \sim \exp\left(-\dfrac{\mu_k^2}{2\kappa} \cdot t\right),$$

with extremely slow damping.

Now we will consider the solutions close to the solutions to the second series on $\tilde{\tilde{M}}$, for which $I_a \approx I_b \approx \dfrac{1}{2}$, and others $I_k \approx 0$. Projections on $\tilde{\tilde{M}}$ of these solutions are described by a following collection of equations.

At first, it is the equations for variables $\left\{\tilde{I}_b = I_b - \dfrac{1}{2},\ \phi_b\right\}$:

(3.12)
$$\begin{cases} \dfrac{d\tilde{I}_b}{dt} = \kappa \cdot \tan^2 \dfrac{\phi_b}{2} \cdot \tilde{I}_b + \dfrac{\kappa}{\cos^2 \dfrac{\phi_b}{2}} \cdot \sum_{k \neq a,b} \sin^2 \dfrac{\phi_k}{2} \cdot I_k + O(|I|^2), \\ \dfrac{d\phi_b}{dt} = \mu_b - 2\kappa \cdot \tan \dfrac{\phi_b}{2} + O(|I|). \end{cases}$$

Secondly, it is the equations for the remaining variables $\{(I_k, \phi_k)|\ k \neq a, b\}$:

(3.13)
$$\begin{cases} \dfrac{dI_k}{dt} = (-\kappa + \kappa \dfrac{\mu_b}{2} \cos\phi_k + \dfrac{\mu_b}{2} \sin\phi_k) \cdot I_k + O(|I|^2), \\ \dfrac{d\phi_k}{dt} = (\mu_k - \dfrac{\mu_b}{2}) + \dfrac{\mu_b}{2} \cos\phi_k - \kappa \sin\phi_k + O(|I|). \end{cases}$$

Here the denotation is introduced

$$|I| = |\tilde{I}_b| + \sum_{k \neq a,b} |I_k|.$$

The equation for $\phi_b$ shows, that asymptotically $\tan\dfrac{\phi_b}{2} \approx \dfrac{\mu_b}{2\kappa} + O(|I|)$. Hence, asymptotically

$$\frac{d\tilde{I}_b}{dt} > \frac{\mu_b^2}{4\kappa} \cdot \tilde{I}_b + O(|I|^2),$$

whence follows the instability of solution, because $\mu_b = \omega_b - \omega_a \neq 0$, under the supposition.

And, the less dissipative constant $\kappa$, the more instability is.

**6.** I will make some more remarks concerning the given model.

1. The dynamics of solutions (3.4) to the first series is defined by the physical Hamiltonian $H = \psi^* \mathbf{A} \psi$. Accordingly, for these solutions the energy is defined

$$E = \psi^* \mathbf{A} \psi = \omega_m.$$

2. The solutions to the first series are steady irrespective of signs of eigenvalues $\omega_m$.

3. If to draw an analogy of $\mathbf{A}$ with the Hamilton operator of Schrödinger (or Dirac) equation for electron of hydrogen atom, then existence of solutions to the second series becomes the unpleasant fact as anything similar is not discovered in experiment. The supposition about *a smallness of the dissipative constant* ensuring acceptable instability of these solutions would be unique solution to this problem.

4. The quantum force (i.e. a right member of dynamic equations, vanishing at $\kappa = 0$) vanishes only on the solutions to the first series. Besides, the solutions to the second series differ from solutions to the first series by nontrivial dynamics of a vacuum component: $d\chi/dt \neq 0$.

5. The analogy with Schrödinger equation suggests also to investigate the behaviour of system under the influence of short-term *exterior perturbation* $\mathbf{A} \to \mathbf{A} + \mathbf{V}(t)$, for the purpose of modeling of dynamics of quantum transition. The problem however is that Schrödinger equation includes an exterior field, which is *statistically averaged* and has a little in common with the real field describing, presumably, photons or other particles. Differently, the self-oscillatory model of quantum transitions can work only when we know not averaged, but actual perturbation $\mathbf{V}(t)$. And if the Schrödinger equation predicts the transitions for *whatever weak exterior field*, within the limits of self-oscillatory model such is impossible, as the consequence of stability of auto-oscillations.

### 3.2. The model of fermion type.
**1.** The impossibility of the dynamic quantum theory is argued sometimes by the reference to quantum statistics which show essential difference from the classic statistic in the very concept of composite system. Therefore it is useful to have models which would make clear, what, as a matter of fact, is required from the dynamic system describing bosons or fermions.

And it would be desirable not to appeal to quantum formalism, and to use those bases of quantum statistics which were known before developing this formalism. For bosons such basis is reduced to the Planck formula of quantization of oscillator energy $E = n\hbar\omega$, and, accordingly, for energy of system of independent oscillators:

(3.14) $$E = n_1\hbar\omega_1 + n_2\hbar\omega_2 + ... + n_N\hbar\omega_N.$$

In the fermion case it is enough to check up the Pauli Exclusion Principle. Also it is possible to express this principle in terms of the energy spectrum of system of noninteracting fermions:

(3.15) $$E = \hbar\omega_{i_1} + \hbar\omega_{i_2} + ... + \hbar\omega_{i_k}, \quad 1 \leq i_1 < i_2 < ... < i_k \leq N.$$

Here we will consider the elementary model of fermion type[44]. It is obtained as immediate generalization considered above matrix model. We will describe the states of a system by the matrixes $\psi, \chi$ of sizes $N \times k$ and $k \times N$ respectively, such, that $\det \chi\psi \neq 0$, and $k \leq N$ will be the analogue of a number of particles. The defining form is set by following formulas (further $\hbar = 1$):

(3.16) $$\alpha = \alpha_0 + \kappa^{-1}dH - dS,$$

where

(3.17) $$\alpha_0 = \operatorname{Im}\operatorname{Tr}\left(\psi^*d\psi + \chi^*d\chi\right), \quad S = -\arg\det(\chi\psi), \quad H = \operatorname{Tr}\left(\psi^*\mathbf{A}\psi\right),$$

and $\mathbf{A}$ is a constant Hermitian matrix, with eigenvalues $\omega_1 < \omega_2 < ... < \omega_N$. Note, that this system is of Kahler type, with potential and Hamiltonian

(3.18) $$U = \operatorname{Tr}\left(\psi^*\psi + \chi\chi^*\right) - \log|\det \chi\psi|^2, \quad H = \operatorname{Tr}\left(\psi^*\mathbf{A}\psi\right).$$

Write out the dynamic equations:

(3.19) $$\begin{cases} d\psi/dt = i\mathbf{A}\psi - \varepsilon\psi + \varepsilon\chi^*\left(\psi^*\chi^*\right)^{-1}, \\ d\chi/dt = -\varepsilon\chi + \varepsilon\left(\psi^*\chi^*\right)^{-1}\psi^*. \end{cases}$$

The regularity of this system is provided by inequality:

(3.20) $$\tfrac{d}{dt}\log|\det \chi\psi|^2 \geq 0,$$

which is satisfied near to singularity. Let's represent the law of change of action function:

---

[44] To obtain the boson variant, we must substitute a permanent for a determinant in the expression for the action function (see below).

(3.21) $$dS/dt = \operatorname{Re}\operatorname{Tr}(\mathbf{A}\Pi), \quad \Pi = \psi(\chi\psi)^{-1}\chi = \Pi^2.$$

The symmetry of system is described by transformations

(3.22) $$\psi \mapsto \exp(is\mathbf{C})\cdot\psi\cdot\mathbf{u}_k^{-1}, \quad \chi \mapsto \mathbf{v}_k\cdot\chi\cdot\exp(-is\mathbf{C}), \quad s \in \mathbb{R},$$

with arbitrary unitary matrixes $\mathbf{u}_k, \mathbf{v}_k$ and with traceless Hermitian matrixes $\mathbf{C}$, commuting with marix $\mathbf{A}$. The corresponding quasi-integrals are equal to

(3.23) $$\mathbf{Q}_1 = \psi^*\psi - \mathbf{1}_k, \quad \mathbf{Q}_2 = \chi\chi^* - \mathbf{1}_k, \quad Q_C = \operatorname{Tr}(\psi^*\mathbf{C}\psi - \chi\mathbf{C}\chi^*),$$

where $\mathbf{1}_k$ is a unit matrix of $k \times k$.

Let's search for the elementary solutions to dynamic equations in the form of exponential[45]:

(3.24) $$\psi(t) = \psi\cdot\exp(-it\omega), \quad \chi(t) = \exp(it\rho)\cdot\chi$$

where $\omega, \rho$ are the Hermitian matrixes of $k \times k$ that reduces (3.19) to algebraic relations. These relations are easily resolved if at first to get a relations

(3.25) $$[\omega, \mathbf{Q}_1] = 2i\varepsilon\mathbf{Q}_1, \quad [\rho, \mathbf{Q}_2] = 2i\varepsilon\mathbf{Q}_2.$$

Really, from (3.25) the vanishing of quasi-integrals follows:

(3.26) $$\mathbf{Q}_1 = \mathbf{Q}_2 = 0.$$

Using it, we obtain sequentially

(3.27) $$\rho = 0, \quad \chi = \mathbf{v}\psi^*, \quad \mathbf{A}\psi = \psi\omega,$$

where $\mathbf{v}$ is some unitary matrix. Thus, also $Q_C = 0$. As a result, all solutions (3.24) can be discribed as follows. Let $\mathbf{e}_1, \mathbf{e}_2, ..., \mathbf{e}_N$ be a fixed orthonormalized basis of eigenvectors-columns of the matrix $\mathbf{A}$, $\mathbf{A}\mathbf{e}_i = \omega_i\mathbf{e}_i$. Then

(3.28) $$\psi = (\mathbf{e}_{i_1}\mathbf{e}_{i_2}...\mathbf{e}_{i_N})\cdot\mathbf{u}^*, \quad \omega = \mathbf{u}\cdot\operatorname{diag}(\omega_{i_1}, \omega_{i_2}, ..., \omega_{i_N})\cdot\mathbf{u}^*, \quad \chi = \mathbf{v}\psi^*, \quad \rho = 0,$$

where $i_1 < i_2 < ... < i_N$, and $\mathbf{u}, \mathbf{v}$ are arbitrary unitary marixes of $k \times k$. Therefore the manifold $Q$ of all solutions consists of $C_N^k$ connected components $Q(i_1, i_2, ..., i_k)$, isomorphic to a direct product of unitary groups $U_k \times U_k$. This manifold is symplecticly isotropic, because all solutions (3.28) lie on Lagrangian manifold $Nul(dU)$ of critical points of a potential,

(3.29) $$Q \subset Nul(dU) = \{\psi, \chi \mid \psi^*\psi = \mathbf{1}_k, \chi = \mathbf{u}_k \psi^*, \mathbf{u}_k \mathbf{u}_k^* = \mathbf{1}_k\}.$$

---

[45] Obviously, they are the Lie solutions corresponding to the subgroup $U_k \times U_k$ of a complete symmetry group.

For the same reason the dynamics on $Q$ is set by Hamiltonian $H$. Besides, the action function restriction on $Q$ varies according to the hypothesis about energy, as

(3.30) $$dS/dt\big|_{Nul(dU)} = \operatorname{Re}\operatorname{Tr}(\mathbf{A}\mathbf{\Pi})\big|_{Nul(dU)} = \operatorname{Tr}(\mathbf{A}\psi\psi^*)\big|_{Nul(dU)} = H\big|_{Nul(dU)},$$

and

(3.31) $$H\big|_{Q(i_1, i_2, \ldots, i_k)} = \omega_{i_1} + \omega_{i_2} + \ldots + \omega_{i_k}.$$

It should be noted, that the restriction of matrix $\mathbf{P} = \psi\psi^*$ on $Nul(dU)$ is an orthogonal projection operator of rank $k$, and the values (3.31) are strictly all extreme values of restriction $H\big|_{Nul(dU)} = \operatorname{Tr}(\mathbf{AP})$.

Thus, the property (3.15) is observed, and in this sense we really have the system of fermion type. It would be desirable to prove also, that except of manifolds $Q(i_1, i_2, \ldots, i_k)$ other attractors don't exist.

**2.** We will look now, how close the obtained finite-dimensional model can be approached to a field model in usual three-dimensional space $\mathbb{R}^3$. The field theory makes the new essential demand – a possibility of an asymptotic decomposability of system into a direct product that should correspond to passage to the usual statistics of distinguishable particles.

The elementary field version of the given model, which allows to hope for realization of this condition, consists in considering the columns of matrixes $\psi, \chi^*$ as some fields in $\mathbb{R}^3$, that corresponds to a limit $N \to \infty$. In «a distinguishable limit» each such field should look like a wave packet describing a separate particle. And the system decomposition on two mutually remote groups of fermions could be connected with decomposition of these matrixes into submatrixes of identical sizes:

(3.32) $$\psi = (\psi_1, \psi_2), \quad \chi^* = (\chi_1^*, \chi_2^*).$$

This implies also the orthogonality relations $\chi_2 \psi_1 = \chi_1 \psi_2 = 0$. As it is easy to see, at such suppositions the dynamic equation (3.19) really break up into two pairs the same equations for these groups, as it is required.

However on this way there is a difficulty connected with symmetry transformations from the subgroup $SU_k \times SU_k$. These transformations can «mix» the columns of matrixes $\psi, \chi^*$ describing mutually remote wave packets, generating their superpositions. But such superposi-tions contradict the idea of field treatment of particles. Exception is made only *by the permutations* preserving localization, but interchanging the position of packages.

Possibly, it can be explained by primitiveness of given model. The true Hamiltonian can be slightly non-invariant. Or, more realistic model can arise by means of *a holomorphic enclosure* in the given model, and this enclosure can break the specified symmetry.

Anyhow, in the given model the system of $k$ fermions, after a relaxation to a «steady-state», is described by orthonormalized system of fields

(3.33) $$\left(\psi_1(x,t), \psi_2(x,t), \ldots, \psi_k(x,t)\right) = \boldsymbol{\psi}(t),$$

according to the equality $\mathbf{Q}_1 = \boldsymbol{\psi}^* \boldsymbol{\psi} - \mathbf{1}_k = 0$. The full description of system, of course, should include as well $\chi$-fields, but if us the states close to steady interest only, the $\psi$-fields will be enough as shows (3.28).

It is interesting, that the similar classic field model of fermions was offered in 1983 in work [14]. A. Ranada considered the Lagrangian system, including $k$ Dirac fields $\psi_i(x,t)$, with preserving scalar products

(3.34) $$N_{ij} = \int_{\mathbb{R}^3} \psi_i^+ \psi_j d^3\mathbf{r}.$$

Also he had *postulated*, that real dynamics of a system flows on *the subset* of a phase space restricted by equations $N_{ij} = \delta_{ij} C$, where a constant $C$ is of dimensionality of action. I will quote, how this model explains the process of shaping of atoms:

*As an example, let us consider two fields $\psi_1$ and $\psi_2$ of charge $e$, in the Coulomb potential of a helium nucleus. Let $\psi_1$ be bound in the state $1S$ spin up and let $\psi_2$ be a wave packet far away from the nucleus at time $t = 0$. As $N_{12} = 0$, the state $1S$ spin up is not accessible to $\psi_2$ if $\psi_1$ remains in it. If $\psi_2$ becomes bound, it must go to a state orthogonal to $\psi_1$. Consequently the process of formation of complex atoms by classic evolution follows the prescription of the Pauli principle, and two fields can never occupy the same state.*

Obviously, the relations (3.34) actually coincide with our relations $\boldsymbol{\psi}^* \boldsymbol{\psi} = \mathbf{I}_k$ following from the general approach. It would be possible to try, starting from Ranada model, to expand this model by the «vacuum» $\chi$-fields and by the phase function $\Phi = \arg \det(\chi\psi)$, and to obtain the field CD-system, the dynamics of which, after a relaxation, would be close to dynamics of Ranada model. However such model is difficult for accepting seriously, in view of its general aesthetic unattractiveness and non-locality.

The genuine field model of fermions, certainly, should be mathematically more refined. Therefore let's formulate the general conclusions which are extracted from our model.

a) First of all, possibly, it is necessary to connect *fermion* behaviour of system with the property of phase function $\Phi = S/\hbar$ to vary on $\pi$ under the permutation of arguments corresponding to particles: $\Phi \to \Phi + \pi \pmod{2\pi}$, or with the property of antisymmetry of function $\exp iS/\hbar$ (that should generalize the antisymmetry of $\det(\chi\psi)$ as function of matrix $\psi$ columns).

b) The property of asymptotic decomposability follows from the asymptotic additivity of defining form. It is reduced to the additivity of 1-form $\alpha_0$, and to the asymptotic additivity of the Hamiltonian $H$ and phase function (or functional) $\Phi$.

Hence, the fermion character of theory basically should follows from the corresponding properties of phase function. A perspective classic field theory for realization of this scenario is *the theory of monopoles* [7]. In this book the low-energy dynamics of system of $k$ «static monopoles» is considered. The phase space of the system looks like $M = T^*N$, and, it appears, that states in configuration space $N$ can be identified with rational functions $s(z)$:

$$(3.35) \qquad s(z) = p_{k-1}(z)/q_k(z) = \frac{a_1}{z-b_1} + \ldots + \frac{a_k}{z-b_k},$$

where all $b_i$ are different, as is supposed. In this theory the important role play *the eliminant* $\Delta(p_k, q_{k-1})$ of a numerator and denominator: its zero *are excluded* from the space $N$. Last property suggests, that the phase function $\Phi$ in the theory of monopoles needs to be built by means of an eliminant. Namely, it can be defined as the quantity proportional *to argument* of eliminant:

$$(3.36) \qquad \Phi = \frac{1}{2} \cdot \arg \Delta(p_k, q_{k-1}) = \frac{1}{2} \cdot \sum_i \arg a_i + \sum_{i<j} \arg(b_i - b_j).$$

Obviously, such definition ensures both the asymptotic additivity and necessary property of antisymmetry. Besides, the summands $\frac{1}{2}\arg a_i$ generate non-standard (multiple $\pi$) periods of phase functional. Possibly, these periods should be connected with *the spinor* character of the wave functions of fermions. Standard periods, multiple $2\pi$, happen from the summands $\arg(b_i - b_j)$. Obviously, these summands can be interpreted as multivalued Hamiltonians of paired *interaction* of particles. The general supposition from here follows, that in a basis of fermi-statistics the special, nonclassic interaction described by multivalued Hamiltonians lies. The definitive clearing of this problem demands the elaboration of the corresponding field theory. The first problem which should be become clear when elaborating such a model, is whether the variation $\delta\Phi$ has a local form in the primary field variables.

### 3.3. The harmonic oscillator model.

Here the Kahler system with potential (for simplicity here $\hbar = 1$)

$$(3.37) \qquad U = \|F\|^2 + |z|^2 - \log|F(z)|^2,$$

and with Hamiltonian $H = \omega_0 |z|^2$ is considered. The function $F$ is understood as an element of the Fock-Bargmann Hilbert space $H$ of holomorphic entire functions:

$$(3.38) \qquad \|F\|^2 = \pi^{-1} \int |F(x)|^2 e^{-\bar{x}x} d^2 x.$$

Let $\varepsilon = \kappa/2$. Then equations of motion of such system can be noted in the form

$$(3.39) \qquad \begin{cases} \dfrac{dz}{dt} = (i\omega_0 - \varepsilon) z + \varepsilon \overline{\dfrac{\partial}{\partial z} \log F(t,z)} \\ \dfrac{\partial F(t,x)}{\partial t} = -\varepsilon F(t,x) + \varepsilon \dfrac{f(\bar{z};x)}{F(t,z)} \end{cases}$$

Here $f(\bar{z})$ is a coherent state

$$(3.40) \qquad f(\bar{z};x) = \exp(x\bar{z}).$$

Let's search for the Lie solutions ($\omega, \xi \in \mathbb{R}$):

$$(3.41) \qquad z(t) = e^{it\xi} z, \quad F(t,x) = e^{i\omega t} F\left(e^{-it\xi} x\right),$$

as we expect, that auto-oscillations interesting us are the solutions of such type.

The basic equations for such solutions are derived in section 2.3. To formulate these equations let's specify the corresponding general objects.

The potential of Kahler metric of the very oscillator here is $V(z) = |z|^2$. The scalar product of coherent states is equal to

$$(3.42) \qquad (f(\bar{z}_1), f(\bar{z}_2)) = \exp(\bar{z}_1 z_2).$$

In this case we have the action of symmetry group $U(1) \times U(1)$:

$$(3.43) \qquad (z, F) \mapsto \left(e^{i\xi} z, e^{i\omega} T(e^{i\xi}) F\right),$$

where $T$ is the unitary representation of subgroup $U(1)$:

$$(3.44) \qquad \left(T(e^{i\xi}) F\right)(x) = F(e^{-i\xi} x).$$

The infinitesimal operator of this representation looks like

$$(3.45) \qquad \mathbf{A}(\xi) = \xi x \cdot \partial/\partial x,$$

and its covariant symbol is equal to

$$(3.46) \qquad A(z, \xi) = \xi |z|^2.$$

Write out the differential equation for a symbol of resolvent of operator $\mathbf{A}(\xi)$:

(3.47) $$\left(\xi|z|^2 + \xi z \partial/\partial z - \omega + i\varepsilon\right) R(z,\xi,\omega - i\varepsilon) = 1.$$

It allows to express a derivative $\partial R/\partial z$ through a function $R$ that will be used below.

Let's accept as unknowns the dimensionless quantities

(3.48) $$p = |z|^2, \quad q = \varepsilon/\xi,$$

considering, that frequency $\omega$ will be defined from the formula $A(z,\xi) = \omega$. Let's introduce also *the dimensionless parameter of model*

(3.49) $$\mu = \omega_0/\varepsilon.$$

Let $g(p,q)$ be the function $i\xi \cdot R(z,\xi,\omega - i\varepsilon)$ expressed through the $p$, $q$, after of substitution (the dependences $g$ on known parameters is not written out)

(3.50) $$\omega = |z|^2 \cdot \xi = p\xi.$$

The function $g(p,q)$ admits various representations. The first – in the form of a series:

(3.51) $$g(p,q) = \sum_{n=0}^{\infty} \frac{p^n e^{-p}}{n!} \frac{1}{q + i(p-n)},$$

the second – in the form of an integral

(3.52) $$g(p,q) = \frac{1}{1 - e^{-2\pi(q+ip)}} \int_0^{2\pi} \exp\{p(e^{it} - 1 - it) - qt\} dt.$$

For a search of $p, q$ the general relations of Chapter 2 give one complex equation

(3.53) $$g(p,q) = \frac{q}{q^2 + \mu pq - p}.$$

Note that for a qualitative computer analysis it is convenient to reduce the equation (3.53) to the system of two real equations:

(3.54) $$\begin{cases} q^2 + \mu pq - p - \dfrac{q}{\operatorname{Re} g(p,q)} = 0, \\ \operatorname{Im} g(p,q) = 0. \end{cases}$$

First of all, it is interesting, whether equation (3.53) has the solutions corresponding to the known stationary states of quantum oscillator. To obtain such solutions, we must suppose, that $\mu = \omega_0/\varepsilon$ is *the big parameter*, $\mu \to \infty$. From the physical point of view this means that the dissipative constant $\kappa$ should be *small* with respect to characteristic atomic frequencies. In that case from the representation (3.51) follows, that equation (3.53) supposes the asymptotic solutions

(3.55) $$p = n + O(\mu^{-2}), \quad q = \mu^{-1} + O(\mu^{-3}),$$

where *n* is arbitrary natural number. The orders of small deviations are explained by that functions $p(1/\mu)$, $q(1/\mu)$ are even and odd respectively. Besides, for any value $\mu$ there is a solution describing *the ground state* of oscillator. It corresponds to the value $p = 0$, but equation (3.53) does not concern to it, as (3.53) has been obtained with use of (3.47) in the supposition $p \neq 0$.

In the same conditions there is also a second infinite series of solutions for which, at increasing, *p* approuches to half-integer values. A computer shows, at least for small quantum numbers $n$, that auto-oscillations of these two series differ by stability of solutions. Solutions to the first series are steady, and the second are unsteady.

Further, it is easy to check up, that the vector $F$ of Hilbert space for $n$-th solution to the first series converges at $\mu \to \infty$ to the eigenvector of $n$-th steady state of quantum oscillator.

Generally solutions exist for all $n \geq N(\mu)$, where $N(\mu)$ is the number of first (after a ground state) raised level of oscillator. It is clear from the qualitative analysis of solutions to equations (3.54) in a plane $(p, q)$ and from the study of asymptotics of big quantum numbers, i.e. the solutions with great values $p$, at arbitrary fixed value of parameter $\mu$. By the way, in view of the formula $p = |z|^2 / \hbar$ such asymptotics should describe a «classic limit» $\hbar \to 0$.

The solution to last problem is reduced to the search by a saddle point method of asymptotics of function $g(p,q)$:

(3.56) $$g(p,q) = \sqrt{\frac{\pi}{2p}} \left\{ 1 + \frac{1}{p}\left(q^2 - 1/6 + iq/3\right) \right\} \cdot \text{cth}\left(\pi(q+ip)\right) - \frac{q}{p} + O(p^{-3/2}).$$

After that, using the formula

(3.57) $$\text{cth}\left(\pi(q+ip)\right) = \frac{\text{sh}(2\pi q) - i\sin(2\pi p)}{\text{ch}(2\pi q) - \cos(2\pi p)},$$

we obtain two series of asymptotic solutions to equation (3.53):

(3.58) $$p = n + O(n^{-1}), \quad q = \mu^{-1} + \sqrt{\frac{2}{\pi n}\frac{\text{th}(\pi/\mu)}{\mu^2}} + O(n^{-1}),$$

(3.59) $$p = n + \frac{1}{2} + O(n^{-1}), \quad q = \mu^{-1} + \sqrt{\frac{2}{\pi n}\frac{\text{cth}(\pi/\mu)}{\mu^2}} + O(n^{-1}).$$

To the solutions (3.58) there correspond solutions to the equations of motion

$$z(t) \approx e^{i\omega_0 t} z_0, \quad F(t,x) \approx e^{i\alpha} x^n / \sqrt{n!}.$$

From here follows, that *the physical Hamiltonian* (in the sense of §9 Chapters 2), in the region of the big quantum numbers, is close to a classic Hamiltonian $H = \omega_0 |z|^2$, and does not depend on $F$. Hence, we have obtained a usual classic limit (the Hamilton dynamics of a $z$-subsystem with Hamiltonian $H = \omega_0 |z|^2$). The solutions described by the formula (3.59), as unstable, don't concern to a quantum spectrum.

Let's estimate now how the number $N(\mu)$ grows at $\mu = \omega_0 / \varepsilon \to 0$. Let's use asymptotics (3.56) for a deriving of the approximate equation *of $n$-th* oval $\operatorname{Im} g(p,q) = 0$, for which $p \approx n$:

(3.60) $$\frac{q \operatorname{sh}(2\pi q)}{3n} \approx \sin(2\pi p).$$

Replacing, according to (3.58), $q$ by $\mu^{-1}$, we obtain the inequality

(3.61) $$n \geq \frac{1}{3\mu} \operatorname{sh}(2\pi / \mu) \approx N(\mu).$$

We see that this growth is of exponential type. Actually it means the existence of boundary from below for the oscillator frequency $\omega_0$ behind which a ground state is possible only. For clearing up of behaviour of an oscillator at such small values of frequencies, we will write out the dynamic equations in the second order of the perturbations theory:

(3.62) $$\begin{cases} \dfrac{dz}{dt} = (i/2 - \mu/8)\omega_0 z + O(\mu^2), \\ \dfrac{d\varphi}{dt} = \omega_0 |z|^2 / 2 + O(\mu^3), \end{cases}$$

considering $\mu = \omega_0 / \varepsilon$ as a small parameter. Here the phase $\varphi(t)$ enters into the expression for a vacuum component of solution

$$F(t) = e^{i\varphi(t)} S(z(t)) / \|S(z(t))\|,$$

and

$$S(z) = f(\bar{z}) + i\mu/2 \cdot \left(\bar{z} \cdot \partial / \partial \bar{z} - |z|^2\right) f(\bar{z}) + O(\mu^2).$$

We see, that such low-frequency oscillator behaves as a damping classic oscillator, with frequency $\omega_0 / 2$ and a damping factor $\omega_0^2 / 8\varepsilon$. Such behaviour differs from that, which the usual quantum model gives. For example, in application to the oscillators of electromagnetic field[46] it means, that the vacuum is opaque for the radio-waves, with wave length exceeding some critical value. Therefore the constant $\kappa$ should be extremely small. But too small it cannot be, as it will reduce to the growth of the relaxation time to quantum-steady states.

---

[46] For that, of course, we have not a solid ground now.

### 3.4. The model of a free massless particle.

Let's show, that the idea of special interaction of particles with vacuum also gives an explanation of inertial motion of particles. We will consider one more special case of the general CS-model described in Chapter 2.

Let $M = T^*\mathbb{R}^3$ be a usual coordinate-momentum phase space. We will consider it as the complex space $M = \mathbb{C}^3$. Points of this space will be denoted by (setting $\hbar = 1$, $c = 1$)

(3.63) $$\mathbf{z} = \mathbf{q} - ia^2\mathbf{p} \in \mathbb{C}^3.$$

We take the following system of coherent states $f(\overline{\mathbf{z}}) \in L_2(\mathbb{R}^3)$:

(3.64) $$f(\overline{\mathbf{z}};\mathbf{x}) = \exp\left\{-\frac{(\mathbf{x}-\overline{\mathbf{z}})^2}{2a^2}\right\}.$$

Here $a$ is the fixed parameter defining the sizes of a wave packet $f(\overline{\mathbf{z}})$. Recognising that such packages should be something like to «frozen» solitons of the future fundamental field theory, it is necessary to accept, that package sizes are proportional to the fundamental length $c/\kappa$:

(3.65) $$a = \rho \cdot c / \kappa,$$

with the dimensionless multiplier $\rho$ which value must be defined by the field theory. This multiplier is the dimensionless parameter of the model.

Note, that for packages of the given form the transformation of a function $F$ to a holomorphic function $F(\mathbf{z}) = (F, f(\overline{\mathbf{z}}))$ is known as the Weierstrass transformation.

We choose the coherent states $f(\overline{\mathbf{z}})$ because the corresponding potential of Kahler metrics in $\mathbb{C}^3$:

(3.66) $$V = \log\|f(\overline{\mathbf{z}})\|^2 = -\frac{(\mathbf{z}-\overline{\mathbf{z}})^2}{4a^2} + \frac{3}{2}\log(\pi a^2) = a^2\mathbf{p}^2 + \frac{3}{2}\log(\pi a^2),$$

leads to the standard action form

(3.67) $$\alpha_M = \operatorname{Im}\partial V = \mathbf{p}d\mathbf{q}.$$

Let $G$ be the motion group of Euclidean space:

(3.68) $$g\mathbf{x} = O\mathbf{x} + \mathbf{h}, \quad \mathbf{x} \in \mathbb{R}^3.$$

The unitary representation $T$ of this group in $L_2(\mathbb{R}^3)$ acts by the formula

(3.69) $$(T(g)F)(\mathbf{x}) = F\bigl(O^{-1}(\mathbf{x}-\mathbf{h})\bigr),$$

and

(3.70) $$T(g)f(\overline{\mathbf{z}}) = f\left(\overline{g\mathbf{z}}\right).$$

The last, that is necessary for the CS-model deriving, is to choose a $G$-invariant Hamiltonian:

(3.71) $$H = H\left(|\mathbf{p}|\right).$$

Let's write out the dynamic equations of obtained system:

(3.72) $$\begin{cases} \dot{\mathbf{z}} = H'(|\mathbf{p}|)\dfrac{\mathbf{p}}{|\mathbf{p}|} + 2\varepsilon a^2 \left(i\mathbf{p} + \dfrac{\partial}{\partial \overline{\mathbf{z}}}\log\left(f(\overline{\mathbf{z}}), F\right)\right) \\ \dot{F} = -\varepsilon F + \varepsilon f(\overline{\mathbf{z}})/\left(f(\overline{\mathbf{z}}), F\right) \end{cases}$$

Our purpose consists in finding solutions to these equations, following the general scheme of Chapter 2. Let's consider at first *a subgroup of the shifts* with generators

(3.73) $$A(\mathbf{z}, \xi) = \langle P(\mathbf{z}), \xi \rangle = \mathbf{p} \cdot \mathbf{v},$$

where $\mathbf{v}$ is a velocity vector. Therefore in this case the element $\xi$ of a Lie algebra is reduced to a velocity vector, and the Lie solution looks like

(3.74) $$\mathbf{z}(t) = \mathbf{z} + \mathbf{v}t, \quad F(t; \mathbf{x}) = e^{i\omega t} F(\mathbf{x} - \mathbf{v}t).$$

Let's use the results of Chapter 2. The function $R$ in this case:

(3.75) $$R(\mathbf{z}, \xi, \omega - i\varepsilon) = -i\int_0^\infty \exp\left\{-\dfrac{\mathbf{v}^2 t^2}{4a^2} - \varepsilon t + i(\mathbf{p} \cdot \mathbf{v} - \omega)t\right\}dt, \quad \varepsilon = \kappa/2,$$

and the equations defining Lie solutions, are reduced to one vectorial relation

(3.76) $$\dfrac{\mathbf{p}}{|\mathbf{p}|} \cdot H'(|\mathbf{p}|) - \mathbf{v} + i\varepsilon R^{-1}\dfrac{\partial R}{\partial \mathbf{p}} = 0,$$

and to one scalar relation

(3.77) $$\mathbf{p} \cdot \mathbf{v} = \omega.$$

From (3.76) follows, that the momentum and velocity are proportional:

(3.78) $$\mathbf{p} = \mu \cdot \mathbf{v}.$$

We will consider the case, when $H'(|\mathbf{p}|) > 0$. Then, as it is easy to see, $\mu > 0$.

Let's choose a frame, in which $\mathbf{v} = (0, 0, v)$, and $p_3 > 0$. Then $p_1 = p_2 = 0$, the value $p_3 > 0$ can be arbitrary, and all equations are reduced to scalar equations for parameter $\mu$ and a particle velocity $v > 0$:

(3.79) $$H'(\mu v) = v \cdot (1 + \Phi(v/\rho)),$$

(3.80) $$\mu v^2 = \omega,$$

where we have introduced the function $\Phi(x)$, $x \geq 0$:

(3.81) $$\Phi(x) = \frac{\int_0^\infty t \exp\{-x^2 t^2 - t\}\, dt}{\int_0^\infty \exp\{-x^2 t^2 - t\}\, dt}.$$

This function is monotonically decreasing and has the asymptoticses:

(3.82) $$\Phi(x) = 1 - 4x^2 + O(x^4),\ x \to 0,$$
$$\Phi(x) = \frac{1}{x\sqrt{\pi}} + O(x^{-2}),\ x \to \infty.$$

Let's consider the massless Hamiltonian $H = |\mathbf{p}|$ which in the Hamilton theory describes particles, moving with the speed of light. Then equation (3.79) becomes the equation for definition of the velocity. Let $v(\rho)$ be a unique solution to this equation. Then $v(\rho)$ is a monotone function of parameter $\rho$, and $v(0) = 1$, $v(\infty) = 1/2$. Let's consider our particle well localised in the sense that $\rho \ll 1$. This supposition ensures a particle velocity, close to the speed of light:

(3.83) $$v(\rho) \approx 1 - \sqrt{\frac{\rho}{\pi}} + O(\rho).$$

The frequency $\omega > 0$ can be arbitrary, and the equality (3.80) gives the value of parameter $\mu$. Thereby, all Lie solutions (for a subgroup of shifts) are discovered.

If to ignore a vacuum component $F$, then these solutions exactly coincide with solutions to Hamilton system with the Hamiltonian $H_0 = v \cdot |\mathbf{p}|$. According the formulas (3.78) and (3.80), for the energy of a particle we obtain $E = \hbar\omega$. It is remarkable, that the particle velocity does not depend in any way on its energy and momentum, or from frequency and length of a wave. This property is characteristic for the strictly relativistic theory.

Let's check up, that, as well as in the model of an oscillator considered above, the vacuum component $F$ at a limit $\varepsilon = \kappa/2 \to 0$ looks as a corresponding one-particle wave function. In this case we have a free particle, and we must get a de Broglie's flat wave. According to the Chapter 2, we have, to within normalization,

(3.84) $$F = \mathbf{R}f(\overline{\mathbf{z}}) = -i\int_0^\infty e^{-(\varepsilon+i\omega)t} T(e^{-\xi t}) f(\overline{\mathbf{z}}) dt,$$

that gives, to within a constant factor:

(3.85) $$F(x) = \exp\left\{-\frac{(x_1-q_1)^2+(x_2-q_2)^2}{2a^2} + ip_3 x_3\right\} \cdot \int_0^\infty \exp\left\{-\frac{(x_3-q_3+vt)^2}{2a^2} - \varepsilon t\right\} dt.$$

Obviously, at $\kappa \to 0$ we will obtain formally $a \to \infty$, and really $F$ turns into a flat wave $F \sim \exp(ip_3 x_3)$.

At the same time, at $a \to 0$ a particle becomes a *point*, and we see, that behind a particle arises a narrow wave tail, with the breadth $a$, and with the length $v/\kappa$.

Note one peculiarity. The Hamiltonian $H_{cl} = v \cdot |\mathbf{p}|$ is not *physical* in the sense of §9 Chapters 2 because the vacuum component $F$ is not stationary. Basically, the similar phenomenon is available also in oscillator model, but there the non-stationarity of $F$ has the character of oscillations near the average value. However, in the given case a particle makes unbounded motion and averaging of $F$ is not possible.

Let's look still, for the completeness of a picture that turns out for the Lie solutions in more complicated case of the general *motion group* $G$ of Euclidean space. The purpose still consists in search of quantities $\mathbf{z}, \xi, \omega$. In this case $\xi$ will be identified with the pair $\{\mathbf{v}, \mathbf{u}\}$ including both a linear velocity $\mathbf{v}$ and an angular velocity $\mathbf{u}$. Indeed, the generator of three-dimensional motions $\exp(t\xi)$ on $M$ is the function

(3.86) $$A(\mathbf{z}, \xi) = \langle P(\mathbf{z}), \xi \rangle = \mathbf{p} \cdot \mathbf{v} + [\mathbf{q}, \mathbf{p}] \cdot \mathbf{u}.$$

The basic equations can be noted in the form

(3.87) $$\begin{cases} A(\mathbf{z}, \xi) = \omega, \\ \partial/\partial \mathbf{z}\left(|\mathbf{p}| - A(\mathbf{z}, \xi)\right) = O(\varepsilon), \quad \varepsilon = \kappa/2 \end{cases}$$

and, at first glance, the elementary way to solve them is to apply the theory of perturbations to a solution with $\varepsilon = 0$. It is easy to transform (3.87) to the form

(3.88) $$\begin{cases} |\mathbf{p}| = \omega + O(\varepsilon), \\ \mathbf{v} + [\mathbf{u}, \mathbf{q}] = \omega \frac{\mathbf{p}}{|\mathbf{p}|^2} + O(\varepsilon), \\ [\mathbf{u}, \mathbf{p}] = O(\varepsilon). \end{cases}$$

In the case $\mathbf{u} = \mathbf{0}$ the solution easily can be defined by the theory of perturbations and this solution coincides with already discovered. However, if $\mathbf{u} \neq \mathbf{0}$, the third of equations (3.88) is not obliged to have a solution $\mathbf{p}, \mathbf{q}$ at arbitrarily small function in a right member. It shows, that in the case $\mathbf{u} \neq \mathbf{0}$ the critical point of the function $f(\mathbf{z}) = |\mathbf{p}| - A(\mathbf{z}, \xi)$ is degenerate and method of perturbations *is inapplicable*. More detailed analysis shows, that equations (3.87), unfortunately, have no other solutions, except already discovered.

Nevertheless, it should be noted, that the equations (3.88) give the grounds to believe, that such more general a Lie solution could describe inertial motion of the massive relativistic particle possessing a spin. In solutions of such type the particle moves on a circular helix and it is possible to consider, that vectors $\mathbf{u}$, $\mathbf{v}$ are parallel. Then the following interpretation arises: $\omega \approx |\mathbf{p}|$ is an «energy», $\mathbf{v}$ is an average velocity of the particle, and $|\mathbf{v}|^2 \approx 1 - \|[\mathbf{u}, \mathbf{q}]\|^2 = const$. Further, the projection of a momentum $\mathbf{p}_\parallel$ to a drift direction $\mathbf{v}$ is a «physical» momentum, the module of an orthogonal component $|\mathbf{p}_\perp|$ is a mass, and $[\mathbf{q}, \mathbf{p}]$ is the particle spin.

Possibly the model considered here is too primitive to give a dynamic model of a massive relativistic particle. But, apparently, the principle obstacles for this are absent. Really, it is possible, using the toy example of Appendix 2 to construct the system describing a free motion of a massless particle in *cylindrical* space, obtained by adding of the closed additional space-like dimension. Also it is possible to consider the coherent states connected with *twistors*, etc. However construction of a realistic model of elementary particles is possible only within the frameworks of *local field theory*.

### 3.5. The spin model.

Let's formulate a special case of the CS-model describing the spin $m/2$. Here the phase space is a two-dimensional sphere $M = S^2 \approx \mathbb{C}P^1$. The holomorphic sections of a line fiber bundle are identified with homogeneous polynomials of degree $m$, and coordinates $x$ are considered to be the homogeneous coordinates in $\mathbb{C}P^1$:

$$p(x) = (x_1 : x_2) \in \mathbb{C}P^1.$$

We can consider them also as usual coordinates in a principal bundle $\mathbb{C}^2 \setminus \{0\} \to \mathbb{C}P^1$. The Hermitian metric in fiber bundle is set by the formula

(3.89) $$\langle F, F \rangle_{p(x)} = m! \cdot |F(x)|^2 \cdot (x^* x)^{-m},$$

and coherent states look like

(3.90) $$E(\overline{s}; x) = (s^* x)^m / m!,$$

where the column $s = \begin{pmatrix} s_1 \\ s_2 \end{pmatrix} \in \mathbb{C}^2$ is considered as usual spinor, because the symmetry group here is the spinor group of Euclidean space, $G = SU(2)$.

At a corresponding normalization of a scalar product in the space of sections we will have the formulas

(3.91) $$\|F\|^2 = \sum_{k=0}^{m}|F_k|^2, \text{ где } F(x) = \sum_{k=0}^{m} F_k \frac{x_1^k x_2^{m-k}}{\sqrt{k!(m-k)!}},$$

(3.92) $$\|E(\bar{s})\|^2 = (s^*s)^m/m!,$$

(3.93) $$F(s) = (F, E(\bar{s})).$$

The unitary representation of the group acts by the formula

(3.94) $$(T(g)F)(x) = F(g^{-1}x),$$

and

(3.95) $$T(g)E(\bar{s}) = E(\overline{gs}), \quad g \in SU(2).$$

Let's identify the elements $\xi$ of a Lie algebra with traceless Hermitian matrixes $A$ of the second order:

$$\xi = iA, \ A = A^*, \ \operatorname{Sp} A = 0,$$

Also we write out expression for generators of rotations $s \mapsto \exp(itA)s$ on a sphere $M = S^2 = \mathbb{C}P^1$:

(3.96) $$\langle P;\xi\rangle_{p(s)} = m \cdot \frac{s^*As}{s^*s}.$$

As the Hamiltonian we take the function, proportional to a $z$-component of a spin vector. The spine vector $\boldsymbol{S} = (S_1, S_2, S_3)$ is built by means of Pauli matrixes $\sigma_1, \sigma_2, \sigma_3$:

(3.97) $$S_\alpha = -\frac{m}{2}\frac{s^*\sigma_\alpha s}{s^*s}, \quad \alpha = 1,2,3.$$

Note also, that $|\boldsymbol{S}| = m/2$ and that the correspondence $p(s) \leftrightarrow \boldsymbol{S}$ realizes a known isomorphism:

(3.98) $$\mathbb{C}P^1 \approx S^2.$$

So, we take the Hamiltonian:

(3.99) $$H(p(s)) = \lambda S_3 = m \cdot \frac{s^*A_0 s}{s^*s}, \quad A_0 = \operatorname{diag}(-\lambda/2, \lambda/2).$$

Definitively, the system is set by the Kahler potential

(3.100) $$U(p(s), F) = \|F\|^2 - \log\langle F, F\rangle_{p(s)}$$

and by the Hamiltonian (3.99). This system is symmetric with respect to the group $G_H = U(1)$ coinciding with the group of rotations round an axis $z$. According to the general scheme of Chapter 2, the Lie solutions for a subgroup $G_H$ in this case look like

(3.101) $$p(s(t)) = p(e^{itA}s), \quad F(t;x) = e^{i\omega t}F(e^{-itA}x)$$

with a diagonal matrix $A = \text{diag}(a, -a)$, and their search is reduced to a solution to a finite set of equations concerning the quantities $\omega$, $a$ and spinor $s$ (to within proportionality):

(3.102) $$\begin{cases} \text{Re}\, R(p(s), a, \omega) = 0, \\ \dfrac{\partial}{\partial s_k}\left(m \cdot \dfrac{s^*(A_0 - A)s}{s^*s} + i\varepsilon \log R(p(s), a, \omega)\right) = 0. \end{cases}$$

The resolvent symbol $R(z; \omega, \xi) = R(p(s), a, \omega) = (\mathbf{R}E, E)/\|E\|^2$ is easy calculated directly from the resolvent definition

(3.103) $$\mathbf{R} = (\mathbf{A}(\xi) - \omega + i\varepsilon)^{-1} = (a \cdot (x_1 \partial / \partial x_1 - x_2 \partial / \partial x_2) - \omega + i\varepsilon)^{-1},$$

because the eigenvectors are known:

(3.104) $$\mathbf{R}(x_1^p x_2^q) = \frac{1}{a \cdot (p - q) - \omega + i\varepsilon} \cdot x_1^p x_2^q, \quad p + q = m.$$

Therefore

(3.105) $$R(p(s), a, \omega) = \sum_{k=0}^{m} \frac{C_m^k}{(2k - m)a - \omega + i\varepsilon} \cdot \frac{|s_1|^{2k} |s_2|^{2m-2k}}{(|s_1|^2 + |s_2|^2)^m}.$$

Note also, that the «Hamilton operator», in the sense of Chapter 2, here is a differential operator

(3.106) $$\mathbf{H} = \frac{1}{2} \cdot (x_1 \partial / \partial x_1 - x_2 \partial / \partial x_2).$$

Unfortunately, equations (3.102) analytically don't be solved, except the case $m = 1$.

To model the given dynamic system on computer, it is more convenient to start with *the expanded* system. It is set by the potential

(3.107) $$U(s, F) = \|F\|^2 + s^*s - \log|F(s)|^2,$$

and by the Hamiltonian

(3.108) $$H(s) = s^* A_0 s.$$

Reduction of this system by group of phase transformations of a spinor $s$ gives the spin system. This expanded system has the following quasi-integrals

(3.109) $$Q_1 = \|F\|^2 - 1, \quad Q_2 = s^*s - m.$$

Let's write out the dynamic equations of the expanded system:

(3.110) $$\begin{cases} \dot{F}(x) = -\varepsilon F(x) + \varepsilon \dfrac{(s^*x)^m}{m!\overline{F(s)}}, \\ \dot{s} = \left(-\varepsilon + im\lambda\sigma_3/2\right)s + \dfrac{\varepsilon}{\overline{F(s)}} \cdot \dfrac{\overline{\partial F(s)}}{\partial s}, \quad \dfrac{\partial}{\partial s} = \begin{pmatrix} \partial/\partial s_1 \\ \partial/\partial s_2 \end{pmatrix}. \end{cases}$$

The computer research of this system discovers only the limit cycles being the «spectral solutions» described in Chapter 2. For these solutions the spin projection $S_3$ accepts the values $S_3 \approx k - m/2, \; k = 0, 1, 2, \ldots, m$. Hence, the remaining solutions are unstable. To prove the existence of such «not spectral solutions» we can investigate analytically the case of spin $\frac{1}{2}$.

**Remarks.** 1) That usually is called by «quantization of spin projection», in the given model means quantizations *of the precession angle* of the spin moment. The precession *velocity* is defined by a numerical factor $\lambda$, with dimensionality of energy. From comparison with the Pauli equation for electron in a magnetic field follows, that $\lambda = 2\mu_B H / m$, where $\mu_B$ is the Bohr magneton, and $H$ is the intensity of a magnetic field which is considered parallel an axis $z$. The dimensionless parameter of the model is the ration $\lambda/\hbar\kappa$.

2) To zeroes of polynomial $F(x)$ there correspond the poles on an sphere, repelling the phase point $S$ representing a spin state. In a «stationary» state of system when the phase point moves over a cycle $S_3 \approx k - m/2$, the $k$ poles are disposed near the south pole, and $m-k$ poles near the northern.

3) The value $m/2$ of a spin is defined by polynomial $F(x)$ degree, i.e. only by a vacuum component. On the other hand, in the unified field theory the vacuum is some *field*, and it is natural to assume, that it breaks up into components with a definite spin. Thus, *if to rely on the model considered here*, the spin of these vacuum components should be exhibited in form of the spin of stable elementary particles observed in the nature: neutrinos, electrons, protons, photons. In that case the vacuum field must contain only components of spin $\frac{1}{2}$ and spin $1$.

# Appendix 1. Relativity as a dynamic symmetry.

Let's illustrate the approach to relativity on the basis of concept of dynamic symmetry. For this purpose it is enough to consider the system of particles on a straight line.

The state of a separate particle is described by the momentum and coordinate: $z = (p, q)$. A role of dynamic symmetry the motion group $G$ of pseudo-Euclidean plane in this case plays. Its generators are[47] the Hamiltonian $H = \sqrt{m^2 + p^2}$, momentum $p$, and the Lorentz momentum $N = -qH$. These variables generate the Lie algebra of the given group:

(1.1) $\qquad \{H, p\} = 0, \quad \{N, H\} = p, \quad \{N, p\} = H.$

The corresponding one-parameter groups of transformations of a phase space acts as follows (here $t, h, \varepsilon$ are parameters):

(1.2)
$$H_t : (p, q) \mapsto (p, q + t \cdot pH^{-1}),$$
$$P_h : (p, q) \mapsto (p, q + h),$$
$$N_\varepsilon : (p, q) \mapsto \left(\operatorname{ch}(\varepsilon)p + \operatorname{sh}(\varepsilon)H, \ q \cdot H \left(\operatorname{ch}(\varepsilon)H + \operatorname{sh}(\varepsilon)p\right)^{-1}\right).$$

The group $H_t$ action, of course, is simply the phase flow of the Hamilton equations of motion for a free particle.

The formulas (1.2) can seem complicated enough in comparison with simple linear transformations of Lorentz. It should be noted in this respect that the transformations of dynamic symmetry are linearized when passing to the momentums:

(1.3)
$$H_t : (H, p, N) \mapsto (H, p, N - t \cdot p),$$
$$P_h : (H, p, N) \mapsto (H, p, N - h \cdot H),$$
$$N_\varepsilon : (H, p, N) \mapsto (\operatorname{ch}(\varepsilon)H + \operatorname{ch}(\varepsilon)p, \ \operatorname{sh}(\varepsilon)H + \operatorname{ch}(\varepsilon)p, \ N).$$

In particular, we see, that the parameter $\varepsilon$ makes the sense of an angle of hyperbolic turn.

To consider effects of a relativity theory, the rigid body model, in the form of some one-dimensional crystal is required to us. Generally speaking, it demands the fields ensuring a coupling of particles. But, not to think of fields, it is convenient to consider this interaction negligible. It is important for us that, basically, this interaction ensures the existence of an equilibrium $Z$, maybe, very fragile, in the form of set of particles being placed through equal gaps of length $d$:

(1.4) $\qquad Z = (p_0 = 0, q_0 = 0; \ p_1 = 0, q_1 = d; \ p_2 = 0, q_2 = 2d; \ \ldots \ ; \ p_n = 0, q_n = nd).$

This equilibrium $Z \equiv H_t Z = Z(t)$ will be our «rigid body». As the length of a body we will consider the distance between its extreme points $L_0 = q_n - q_0$.

So, we will consider *the system of free particles*, neglecting fields, and using an equilibrium $Z$ as the model of a motionless rigid body. Then the Hamiltonian, momentum, and the Lorentz moment of the system are equal to sums of corresponding one-particle quantities, and for these sums the commutation relations (1.1) are still valid.

---

[47] The speed of light is accepted to equal unit.

Now it is possible to start a deducing of consequences. They are based on the following *definitions*.

**Definition 1**. The body $Z$ shifted on a distance $h$, is set by a phase trajectory

(1.5) $$Z_h(t) = H_t(P_h Z).$$

**Definition 2**. The same body, moving by inertia, is set by a phase trajectory

(1.6) $$Z_\varepsilon(t) = H_t(N_\varepsilon Z).$$

Despite their exterior simplicity, these definitions are very essential in the conceptual plan. Really, the first definition shows, in what sense usual Euclidean geometry is the part of dynamics. The second definition comprises all effects of a theory of relativity[48].

To show it, write out the values of phase coordinates of $a$-th particle of a moving body what they follow from (1.6), using the velocity parameter $v = \text{th}\,\varepsilon$ instead of a hyperbolic angle $\varepsilon$:

(1.7) $$p_a(t) = \frac{mv}{\sqrt{1-v^2}}, \quad q_a(t) = a \cdot d\sqrt{1-v^2} + vt, \quad a = 0,1,\ldots,n.$$

From here it is clear, that the momentum and energy of a moving body correspond to usual relativistic formulas. The same is possible to tell about the body length

$$L(t) = q_n(t) - q_0(t) = L_0\sqrt{1-v^2}.$$

Further, the definitions 1-2 are natural to expand from solid bodies to any *isolated processes* describing the dynamics of closed subsystems. It means that instead of $Z$ it is possible to take the initial conditions of these processes. Let's consider, for example, the «process», including two particle. Let the first particle is immobile, and the second particle is moving with absolute velocity $u$:

(1.8) $$p_1(t) = 0, \; q_1(t) = 0, \quad p_2(t) = \frac{mu}{\sqrt{1-u^2}}, \; q_2(t) = q_2^0 + ut.$$

Then, acting on initial conditions by transformation $N_\varepsilon$, we get the same «process», but moving «as a whole» with the absolute velocity $v = \text{th}\,\varepsilon$:

(1.9) $$p_1(t) = \frac{mv}{\sqrt{1-v^2}}, \; q_1(t) = vt, \quad p_2(t) = \frac{mw}{\sqrt{1-w^2}}, \; q_2(t) = q_2^0\sqrt{1-v^2} + wt.$$

The absolute velocity of the second particle in these formulas is equal to

---

[48] A clear split into geometry and inertial motion is possible not always. The example is given by the dynamic symmetry corresponding to the de Sitter world.

(1.10) $$w = \frac{u+v}{1+uv},$$

demonstrating the relativistic addition of velocities. A mathematical base of this law is the group property of transformations $N_\varepsilon$:

(1.11) $$N_{\varepsilon_1} \cdot N_{\varepsilon_2} = N_{\varepsilon_1 + \varepsilon_2}.$$

In the role of an isolated process also it is possible to take *a clock* of any construction. For example, we can take the «light clock», considering the oscillations of a massless particle between the body ends acting as mirrors. Then from the definition 2 the formula for the magnification of oscillation period of a moving clock is derived

(1.12) $$T = \frac{T_0}{\sqrt{1-v^2}},$$

expressing a relativistic effect of time dilation. Let's show also, that this effect *does not depend on the clock's mechanism*. By the clock with the period $T_0$ it is possible to call any periodic process $Z(t)$:

(1.13) $$Z(t) = H_t Z \equiv H_{t+T_0} Z.$$

Let's use now of a commutation relation[49]

(1.14) $$H_t \cdot N_\varepsilon = P_{t \cdot \text{th}\,\varepsilon} \cdot N_\varepsilon \cdot H_{t/\text{ch}\,\varepsilon}.$$

Applying it to a state $Z$ we obtain, that «moving process» $Z_\varepsilon(t) = H_t(N_\varepsilon Z)$ possesses the following quasi-periodic property

(1.15) $$Z_\varepsilon(t + T_0 \cdot \text{ch}\,\varepsilon) = P_{(t+T_0 \cdot \text{ch}\,\varepsilon)\text{th}\,\varepsilon} Z_\varepsilon(t).$$

In translation on usual language the relation (1.15) means, that process $Z_\varepsilon(t)$ represents a clock, moving with velocity $v = \text{th}\,\varepsilon$, and with period $T = T_0 \cdot \text{ch}\,\varepsilon = T_0/\sqrt{1-v^2}$, as it was required.

Let's show now, that, according to the definition 2, the passage from rest to motion «as a whole» can convert simultaneous events into the non-simultaneous. We will consider the two identical bodies, moving with equal absolute velocities towards each other. Obviously, at some moment there will be coincidence of respective points of both bodies, i.e. this coincidence will be simultaneous events. However if to add to this «process» such common velocity at which the first body stops, the second body becomes shorter, and coincidence of bodies is excluded.

We have formulated the basic relativistic effects in *absolute* frame of reference (in old terms – in the system connected with motionless ether). Considering these effects, it is possible to deduce a relativity principle, as the basic *empirical consequence* of existence of dynamic symmetry: by any measurements it is impossible to discover an absolute motion; all inertial observers are empirically equivalent. Such is the sense of the dynamic approach to relativity.

---
[49] To derive this relation, it is convenient to use the formulas (1.3).

**Historical remark**. There are no doubts what exactly in such spirit H. Poincare understood the relativity[50]. For him the concept of inertial motion, instead of a relativity principle was primary. That is why Poincare *did not deduce* the Lorentz contraction from a relativity principle, his purpose was *inverse*. Dynamic understanding also *excludes* that «decisive step» which has been made by Einstein and which sometimes incorrectly is attributed also to Poincare (the new concept of space and time means). Poincare's ideology was absolutely another[51]. This ideology, unlike Einstein's approach, is quite compatible with a hypothesis of the immobile ether.

## Appendix 2. Singularities in geometrical interpretation of dynamics.

Let's give an example to show that the Minkowski space in theory of relativity can play a role of the convenient method (but not obligatory and *not universal* method) for representation or *interpretation* of the dynamics of a system.

Let's consider the expression for the energy of a relativistic particle $E = \sqrt{\mathbf{p}^2 + m^2}$ as the *Hamiltonian of geodesic motion* of a mass point with the unit velocity on the four-dimensional cylinder $\mathbb{R}^3 \times S^1$. For this purpose we will identify the parameter $m$ with the component of the four-dimensional momentum conjugate to arc coordinate $s$ on a circle $S^1$ of the fixed radius[52], allowing the zero and negative values $m$. Then $\mathbf{p}$ and $m$, obviously, will be the constants of motion.

The sense of this model is that such a particle looks as usual relativistic particle with the mass $|m|$ and momentum $\mathbf{p}$, if to consider the additional closed space-like dimension *to be hidden* (for example, because of a smallness of circle radius). The negative values $m < 0$ should be understood so, that we deal *with an antiparticle*.

The equation of motion

$$ds/dt = m/E = \operatorname{sgn}(m)\sqrt{1 - \dot{\mathbf{x}}^2}$$

shows, that the «proper time» $\tau$ of the particles appearing in relativity, here is necessary *to define* by equality $d\tau = |ds|$. Differently, the proper time is proportional *to number of rotations* along the hidden dimension, and the circle $S^1$ plays a role of interior clock face of a particle. In particular, for zero mass the proper time is deprived of sense.

---

[50] And also H. Lorentz, till the moment when he has accepted the Einstein's point of view.

[51] The root of misunderstanding of Poincare's views consists in considering the special theory of relativity as *the unique correct and definitive point of view* onto relativistic effects.

[52] The value of radius does not play a role.

Further, in the system of such particles it is possible to introduce gravitational interaction, taking the Hamiltonian

$$H = \sum_a E_a - k \sum_{b<c} E_b E_c / |\mathbf{x}_b - \mathbf{x}_c|, \quad E_a = \sqrt{\mathbf{p}_a^2 + m_a^2},$$

where $k$ is the gravitation constant. Then a motion of a test particle of small energy is defined by the Hamiltonian $H_0(\mathbf{x}, \mathbf{p}, m) = E \cdot (1+U)$, where $U = -k \sum_a E_a / |\mathbf{x} - \mathbf{x}_a|$ is the potential of external gravitational field. Such Hamiltonian, by the way, describes a geodesic motion on the cylinder $\mathbb{R}^3 \times S^1$ endowed with the conformally-flat Riemannian metric

$$d\sigma^2 = (1+U)^{-2}(d\mathbf{x}^2 + ds^2)$$

In this metric a motion on the cylinder has a unit velocity[53], therefore we have a relation

$$(d\mathbf{x}/dt)^2 + (ds/dt)^2 = (1+U)^2,$$

from which *the effective pseudo-Euclidean metric of the Minkowski space* is derived

$$d\tau^2 = ds^2 = (1+U)^2 dt^2 - d\mathbf{x}^2 = g_{ik} dx^i dx^k, \quad x^0 = t.$$

We will get the same pseudo-Euclidean metric if, fixing the mass of a test particle, will pass from Hamiltonian $H_0$ to Lagrangian. This procedure gives to us the standard Lagrangian

$$L_0 = -m\sqrt{g_{ik}\dot{x}^i \dot{x}^k}$$

of the relativistic particle in special gravitational field. The corresponding metric $g_{ik}$ is known as Newton's approximation to the GR metric for slow motions. From the equations of motion follows, that the free falling test particle which is close to one of surfaces $1+U = 0$, infinitely approaches to this surface. Together with this process its gyration in the direction of the hidden dimension stops. The particle gets characteristics, impossible for a motion in the region of a weak gravitational field, where $U \approx 0$. Here the concept of proper time for it loses the sense, but the «mass» $m$ remains the same[54], and this mass can be non-zero. There is also a sharp increase of the energy and momentum, and *the particle ceases to be a test one*.

We can now interpret this motion within the limits of a Minkowski space. The calculations show, that singularity on a surface $1+U = 0$ will be the «true» space-time singularity. It follows from the fact that the Ricci curvature invariant $R^{ik}R_{ik}$ is equal to infinity on this singularity.

---

[53] It is the known fact for a geodesic motion in the metric $\gamma_{ab}$, with the Hamiltonian $H_0 = \sqrt{\gamma^{ab} p_a p_b}$.

[54] To change the mass $m$ the Hamiltonian must depend on the hidden variable $s$.

In this case we face with the *singularity in interpretation* of dynamics in geometrical terms of space-time. Moreover, if we include dependence on the hidden variable into the Hamiltonian (for example, replacing the three-dimensional Newton's potential by the solution to Poisson's equation on the cylinder $\mathbb{R}^3 \times S^1$) such interpretation in general becomes problematic. In this case the masses will not be preserved any more, up to transmutation of particles into antiparticles.